\documentclass[english]{article}
\usepackage{preamble}

\begin{document}

\title{A Story of Algebra and Categories I: Algebraic Theories}
\maketitle

\thispagestyle{empty}

\begin{abstract}
        This is the first part of a series of expository notes on homotopical and categorical algebra. It charts a course from ordinary algebraic structures such as monoids, groups, and rings to their various homotopy-coherent generalizations and equivariant versions. The emphasis is on the ideas that lead to and connect the constructions, as well as where and how they appear in nature.
\end{abstract}

\par\vspace{0.5\baselineskip}
\noindent
\makebox[\textwidth][c]{%
  \includegraphics[
    width=\textwidth,
    height=\dimexpr\pagegoal-\pagetotal-4\baselineskip\relax,
    keepaspectratio
  ]{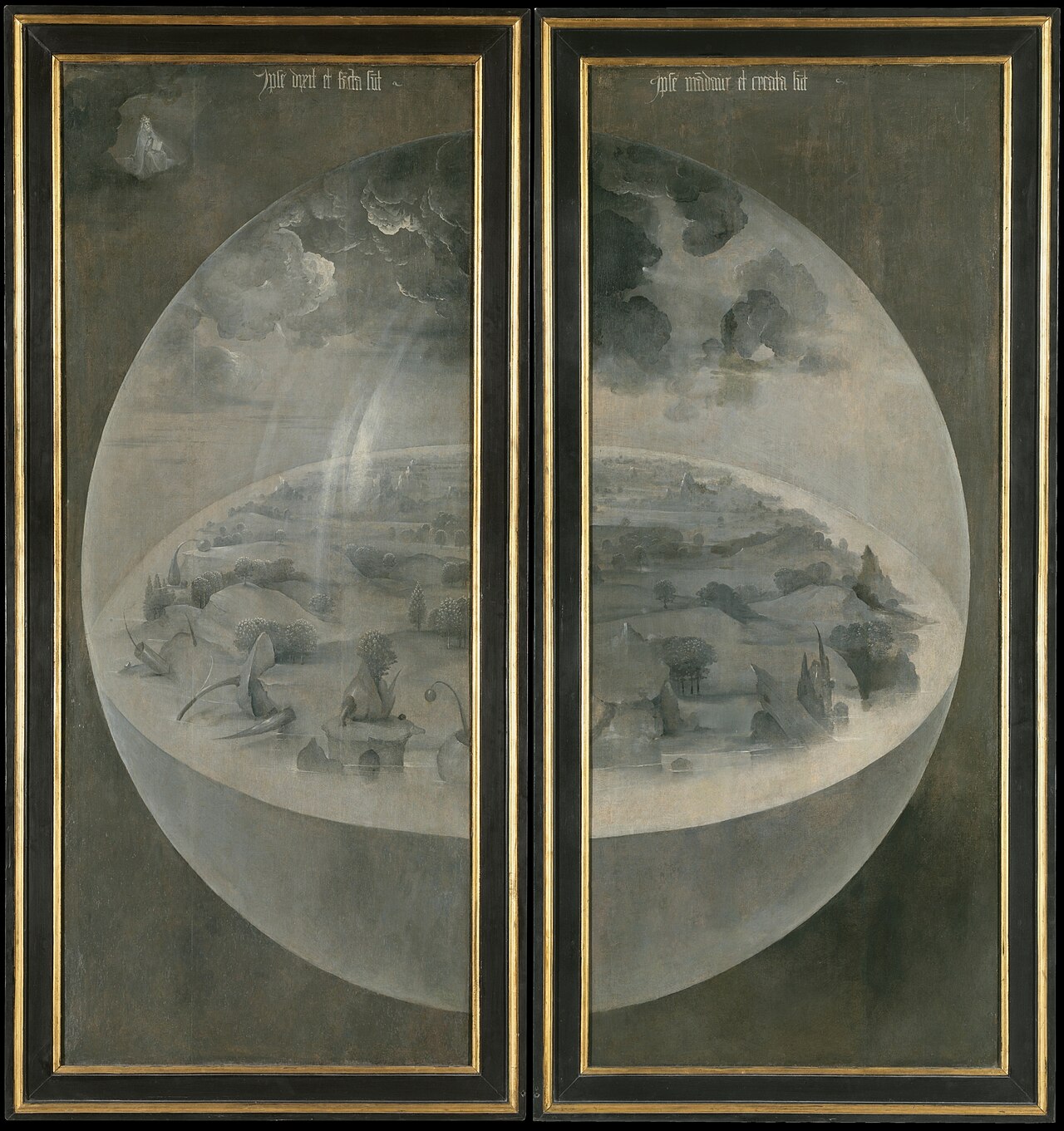}%
}
\par\vspace{0.4\baselineskip}
\noindent
{\centering\small
Hieronymus Bosch, \textit{The Garden of Earthly Delights}
(c.\ 1490--1500), exterior panels:
\textit{The Third Day of Creation}.
Museo del Prado, Madrid.\par}

\tableofcontents

\newpage

\section{Introduction}
\epigraph{For magic consists in this, the true naming of a thing.}{Ursula K. Le Guin, \textit{A Wizard of Earthsea}}

There are two possible meanings to the term \textit{categorical algebra}. 
The first is the use of abstract categorical frameworks --- monads, operads, algebraic theories, and their kin --- 
to bring out the common threads among algebraic structures arising across different areas of mathematics.
Recognition principles such as the monadicity theorem, bar-cobar dualities, as well as the many formulas for free constructions and envelopes belong to this meaning. The second is the study of algebraic structures on categories \textit{themselves}, viewed as objects of some higher category of categories. The extra categorical level brings in new modes of structure such as adjunctions, localizations, and laxity, lifting it beyond a mere special case of the first meaning. The characterization of \textit{stable categories} as modules over spectra and the resulting coherent existence and uniqueness of the smash product on the latter is a compelling demonstration of the power of this conceptual move. Another is the burgeoning study of \textit{dualizable categories} and their role in algebraic $K$-theory. 
 
Yet some of the most striking aspects of the subject lie in the interaction between the two meanings, broadly known as the \textit{microcosm principle}: every algebraic concept finds its native home in its own categorification.
As a guiding observation, algebras can be considered in arbitrary monoidal categories, while commutative algebras require the monoidal structure itself to have a symmetry. A deeper instance is \textit{Hochschild homology}. It is the universal receptacle for traces of endomorphisms of dualizable objects, precisely because it is itself the trace of the identity endofunctor of the category (viewed as a dualizable object of some category of categories). Thus, climbing the categorical ladder not only reveals deeper and more sophisticated phenomena, but also casts the lower rungs in a clearer light.  

\subsection*{About the origins}

As may be gathered from the examples above, the real force of this approach to algebra becomes apparent not in the study of monoids, groups, or rings in the traditional sense, but when such notions are considered in any derived, homotopical, or categorified settings. Accordingly, the historical roots of this theory are intertwined with the broader endeavor to lay workable foundations for such settings. This endeavor was for a long time of a largely disparate nature, as scholars of many schools, each with their own traditions and ambitions, set out to wrangle and harness the many-faced phenomenon we now call \textit{coherence}.

The category theorists tackled this in the most direct manner. In the categorification of a commutative monoid --- a symmetric monoidal category --- the associativity and commutativity conditions are replaced by specified isomorphisms obeying coherence constraints. As appealing as the explicitness of this approach is, categorifying it once more is, though possible, impractical, while categorifying it twice is practically impossible.\footnote{The reader is invited to look at the definition of a \textit{symmetric monoidal 2-category} in \cite[Appendix C]{schommer2009classification} and judge for themselves.}

At the same time, algebraic topologists were grappling with homotopy coherence in topological constructions and their chain complex analogues. The main tool devised for this purpose was the theory of \textit{model categories}, which stands at the opposite pole: instead of spelling out coherence, it directs one to work as strictly as possible while forcibly ensuring invariance under \textit{weak equivalences}. Its achievements command admiration, yet its form resists affection. 

Alongside, and often in close interaction with, the development of the many flavors and spin-offs of model categories, some continued to refine and elaborate the simplicial apparatus as a bookkeeping device for higher coherences. Others sought to trade the simplicial combinatorics for cubical (or dendroidal, or globular). Proponents of \textit{derivators},  \textit{computads},  \textit{crossed complexes}, and the like kept developing their own preferred formalisms, while the algebraists and representation theorists mostly opted for \textit{triangulated categories} and their \textit{dg enhancements} or explicit chain-level formulas. These theories kept proliferating, mutating and hybridizing, producing ungainly offspring such as \textit{monoidal pre-triangulated dg-enhanced opetopic multiderivators}.\footnote{Okay. This one is not real, but it could just as well have been.}

This tower of Babel was not merely an aesthetic inconvenience. It threatened to scatter its inhabitants into an archipelago of insulated communities discovering and rediscovering the same mathematical ideas in mutually unintelligible dialects. What kept the edifice from collapsing was the diligent and fastidious labor of comparison. Bridges were constantly built and re-built between the frameworks as they expanded to ensure that even when people could not quite understand one another, they could take comfort in knowing they were misunderstanding the same thing.

All this time, there was a shared belief that beneath the thicket of technical formalisms lay the \textit{homotopy theory} itself, a platonic form independent of any particular language used to present it. Different languages --- simplicial, triangulated, dg, model-theoretic --- gave access to different aspects of that homotopy theory. But conceptual maneuvers that ought to have been transparent at the level of homotopy theories remained technically formidable in all of them. None of the languages and models was apt for the task of describing homotopy theories in their pure form.

Even with the advent of $\oo$-category theory, the novelty was initially perceived as a collection of competing models for the ``homotopy theory of homotopy theories'': quasi-categories, simplicial categories, relative categories, Segal categories, and complete Segal spaces, to name a few. \textit{Their} comparison then became a central occupation of the field. Meanwhile, the vast body of work produced in the preceding languages did not simply translate itself into any of the new ones. Thus, to most outsiders the whole enterprise seemed too forbidding and volatile to invest in.

Eventually, the quasi-categorical model had become mature enough to provide all the required infrastructure in a self-contained manner. More importantly, it served as a proof of concept that $\infty$-category theory could become a working mathematical language, equipped with the $\infty$-categorical versions of all the basic machinery of ordinary category theory and robust enough to support substantial mathematics.\footnote{Since the standard models are connected by explicit Quillen equivalences, the difficulty in developing the theory in another model is arguably \textit{quantitative}. Though the precise meaning of that remains debatable.} 
With time, the general nervousness around what could be written off as a ``standard argument'' began to subside. A new mode of reasoning emerged, one that distilled the invariant principles of the theory and de-emphasized the particulars of \textit{any} chosen foundation. By now, a significant portion of the classical material has been reworked in this language and much of the new material is being developed natively in it. 

Nevertheless, some artifacts of implementation have remained so deeply ingrained in the collective mind as anchors of thought and sources of intuition that even experts occasionally find their heads screwed on backwards. Having said that, a word of caution is due to those too eager to dispose of the ``old models'' as mere relics of a confused transitional era, or as abandoned scaffolding left standing around a completed building. Many of the insights developed to access and manipulate homotopy-coherent constructions through these models are neither rendered obsolete nor trivialized by modern foundations. Thus, cofibrant replacement may no longer be needed to define a derived functor, but, recast as a particularly tractable kind of resolution, it may still be the most effective way to compute one. When one actually sets out to work with these objects, more often than not one finds oneself treading a path laid down decades ago by people who, by modern standards, did not yet fully grasp what a ``homotopy theory'' truly was.

At the same time, other, old and new, languages still provide points of contact with geometry, combinatorics, algebra, and analysis that an invariant $\infty$-categorical description may deliberately suppress. They also give direct access to structures for which no satisfactory intrinsic account is yet available, and serve as working models for still higher categorical phenomena where the general theory remains less consolidated. What has changed is therefore not quite that one common language has replaced the others, but that it lets us separate invariant content from presentation, which is often desirable even when the presentation itself carries valuable information. 

\subsection*{About these notes}

In this day and age, ordinary categorical notions such as functor categories, limits, and adjunctions have conceptually clear and well-understood $\infty$-categorical definitions and properties. This basic package, however, does not by itself prescribe how algebraic structures such as monoids, groups, abelian groups, rings, Hopf algebras, or Lie algebras should be transported into the $\infty$-categorical setting. Indeed, it is not even clear (or true) that these concepts have a single well-defined $\oo$-categorical generalization. In addition, completely new algebraic structures arise, which make sense only in the homotopy-coherent world.

There is no universal formalism for treating this full range of homotopy-coherent algebraic structures. Instead, there are several overlapping frameworks, whose relations are not always transparent even when everything is formulated natively in $\oo$-categorical language. This is not necessarily a defect. Different frameworks trade the generality of the situations they allow against how much one can actually say about them. They also bring different aspects of the same structure to the fore, and translating between them does real work in revealing how these aspects fit together.

The purpose of these notes is neither to survey the landscape of such frameworks nor to develop rigorously the theory underlying any one of them. 
Rather, it is to invite the reader, whatever their background, to enjoy a compelling story about some of the ideas involved and why one should care about them.
All the statements are ``well known'' even if some are more difficult to extract, in the exact form they are presented here, from the existing literature. Full proofs are included only when they are both simple and instructive. Likewise, references are chosen primarily for the explanations they offer, rather than to document the origins of every idea, result, or technique.

Freed from the demands of a strictly logical or chronological order, the story can unfold somewhat nonlinearly, using ``advanced technology'' to motivate, streamline, or contextualize more basic constructions. This mode of exposition is unlikely to suit the novice. It is intended rather for readers who have already heard many of these stories but still find them difficult to hold together. It assumes fluency in the language of $\oo$-categories --- or at least the ability to feign it convincingly --- while drawing on little beyond the basic facts. Above all, it is written for the young in spirit: those who appreciate being shown how to think about things and would rather be persuaded that a definition is \textit{inevitable} than merely assured that it is \textit{useful}. 


\subsection*{About this chapter}

We begin with the framework of \textit{algebraic theories}, and some of its close relatives, because it is the most accessible from first principles.\footnote{This point is also promoted in \cite{Cranch2011AlgebraicTS}.} Quite a lot can be achieved with it, though we shall also run up against some of its limitations that will eventually lead us to more general frameworks such as symmetric monoidal categories, operads, and monads in subsequent chapters.
The algebraic structures we shall start with are familiar ones, but their homotopy-coherent incarnations will take us rather far afield. The reader who prefers to know in advance where all this is going may consult the concluding section, which offers a retrospective view of the topics covered and how they interconnect.

Algebra is about equations as much as analysis is about inequalities. The question, then, is how to treat algebraic structures in a world where equations have been replaced by coherent systems of homotopies. 
The point of leverage is the innate ability of $\oo$-category theory to handle \textit{diagrams} with all the necessary coherence. Indeed,  homotopy-coherent associativity is already baked into how morphisms in an $\infty$-category compose. Thus, whenever an algebraic notion can be encoded by a diagram of a certain \textit{shape} satisfying certain \textit{properties}, it admits a natural $\infty$-categorical analogue. 
The challenge is therefore to find the right shape and the right properties. 

Yet there is no reason to expect a unique answer. Even classically, the same structure may be encoded by several quite different-looking diagrams. After passing to the $\oo$-categorical world, the ambiguity only grows: a single classical notion may admit more than one inequivalent generalization. Finally, the appropriate indexing shape may not be an ordinary category at all, but itself an $\oo$-category. The ultimate task is therefore to identify the $\oo$-categorical structures worth formalizing, guided both by the internal logic of the theory and by the examples that arise in nature, and construct efficient models for them. 

\section{Monoids}

The most natural starting point for this quest is the theory of \textit{monoids}.
The paradigmatic source of monoids in ordinary mathematics is the \textit{endomorphism monoid} of an object in a category. In fact, this source is universal: monoids are precisely \textit{one-object categories}.\footnote{This should be interpreted as having a \textit{pointed connected} groupoid of objects, rather than a \textit{contractible} one, which we would refer to as having a \textit{unique} object.} The same should be true, and indeed is true, for homotopy-coherent monoids as well. One might even be tempted to take one-object $\oo$-categories as the \textit{definition} of homotopy-coherent monoids. But then, it is not obvious how to recast this as an algebraic structure one can impose on an object of \textit{any} $\infty$-category. Nor is it immediate how to place it on equal footing with other, less privileged, structures, such as \textit{commutative} monoids. We shall thus begin at the beginning and proceed in small steps.

\subsection{Relations among relations}

A \tdef{monoid} is a set $M$ with a binary \textit{multiplication} operation which is associative and unital with respect to a given \textit{unit} element $e \in M$. This structure can be encoded diagrammatically as a pair of maps
\[
    \pt \oto{\ \ e\ \ } M \qquad,\qquad 
    M\times M \oto{\ \ m\ \ } M
\]
rendering the following diagrams commutative

\[\begin{tikzcd}[sep = 3em]
	M & {M\times M} & M & {M\times M\times M} & {M\times M} \\
	& M && {M\times M} & M
	\arrow["{(\Id,e)}", from=1-1, to=1-2]
	\arrow[equals, from=1-1, to=2-2]
	\arrow["m", from=1-2, to=2-2]
	\arrow["{(e,\Id)}"', from=1-3, to=1-2]
	\arrow[equals, from=1-3, to=2-2]
	\arrow["{m\times\Id}", from=1-4, to=1-5]
	\arrow["{\Id \times m}"', from=1-4, to=2-4]
	\arrow["m", from=1-5, to=2-5]
	\arrow["m", from=2-4, to=2-5]
\end{tikzcd}\]

As such, it can be interpreted verbatim in every category $\cC$ with finite products (with $\pt$ the terminal, ie empty product, object). This readily reproduces notions such as \textit{topological monoids}, \textit{ordered monoids}, etc. 

We pause to point out that this definition of a monoid structure is essentially by ``generators and relations''. Namely, we have one binary operation $m(x,y)$, one unary \textit{identity} operation $\Id(x) = x$, and one nullary operation $e$, from which further $n$-ary operations can be produced by composition. Associativity, for example, is then a \textit{relation} identifying two ternary operations thusly produced
\[
    m(m(x,y),z) = m(x,m(y,z)).
\]
In general $\infty$-categories, such relations become \textit{homotopies}, which can themselves be related via higher homotopies, thought of as \textit{relations among relations}.  

To be concrete, in ordinary algebra, the associativity law for \textit{threefold} products implies that for all $n$, any two groupings of an $n$-fold product produce the same result. In an $\oo$-category, the two ternary operations $m(m(x,y),z)$ and $m(x,m(y,z))$ are two points in the space of maps $M^3 \to M$, or equivalently, a single map
\[
    \gamma_3 \colon S^0 \too \Map(M^3,M).
\]
Associativity is then a choice of a \textit{path} between these two points, or equivalently, a null-homotopy for $\gamma_3$.
Given such a null-homotopy, the 5 different groupings of a \textit{fourfold} product fit into a diagram known as the \tdef{Stasheff pentagon}:
\[
\begin{tikzpicture}[
  x=0.85cm,y=0.85cm,
  every node/.style={fill=white, inner sep=1pt}
]
\coordinate (A) at (0,2.6);
\coordinate (B) at (2.4,0.9);
\coordinate (C) at (1.5,-1.8);
\coordinate (D) at (-1.5,-1.8);
\coordinate (E) at (-2.4,0.9);

\draw[shorten <=4pt, shorten >=4pt] (A)--(B)--(C)--(D)--(E)--(A);

\node[anchor=south] at (A) {$m(m(x,y),m(z,w))$};
\node[anchor=west]  at (B) {$m(x,m(y,m(z,w)))$};
\node[anchor=north west] at (C) {$m(x,m(m(y,z),w))$};
\node[anchor=north east] at (D) {$m(m(x,m(y,z)),w)$};
\node[anchor=east]  at (E) {$m(m(m(x,y),z),w)$};
\end{tikzpicture}
\]
One can view this diagram as a \textit{loop}
\[
    \gamma_4\colon S^1 \too \Map(M^4,M).
\]
If $\gamma_4$ fails to be null-homotopic, the fourfold product is in a homotopy-coherent sense \textit{nonunique}; and even when it is null-homotopic, the choice of the null-homotopy itself need not be unique up to homotopy.

Moving further, the 14 possible groupings of a \textit{fivefold} product form the vertices of a \href{https://en.wikipedia.org/wiki/File:Associahedron.gif}{3D polyhedron} with three square faces and six pentagonal ones. A null-homotopy for $\gamma_3$ (threefold associativity) provides the edges of this polyhedron and a null-homotopy for $\gamma_4$ (fourfold coherence) provides fillings for the pentagons. The squares have canonical fillings, so all this assembles into a map whose source is the boundary of said polyhedron
\[
    \gamma_5 \colon S^2 \too \Map(M^5,M),
\]
exhibiting a higher-homotopical ambiguity in the definition of the fivefold product.
This continues indefinitely to produce a sequence of elements 
\[
    \gamma_n\ \in\ \pi_{n-3}\Map(M^n,M) \quad,\quad n \ge 3,
\]
each depending on a chosen null-homotopy of its predecessor. We are thus in a classical obstruction-theory situation. It is only by nullifying iteratively all these obstructions that one can claim to have a truly homotopy-coherent associative multiplication on $M$.

\begin{rem}
    For those who indulge in \textit{negative thinking}, the multiplication map $m(x,y)$ itself can be viewed as a null-homotopy for the unique map
    \[
        \gamma_{2} \colon S^{-1} = \es \too \Map(M^2, M).
    \]
\end{rem}

This leads to our first attempt to ``define'' a homotopy-coherently associative multiplication. The fully spelled-out definition can be found in \cite{Stasheff1963}.

\begin{defn}
    A (nonunital) \tdef{$\AA_n$-structure} on an object $M$ of an $\oo$-category $\cC$ is an iterative sequence of null-homotopies for the maps $\gamma_2,\gamma_3,\dots,\gamma_n$. Such a sequence of null-homotopies for \textit{all} $\gamma_n$ is called a (nonunital) \tdef{$\AA_\oo$-structure}\footnote{Unitality turns out to be (almost) a property (see \Cref{Mon_nu}).}. 
\end{defn}

\begin{rem}
    The coherence diagrams appearing in this process are organized by a sequence of convex polytopes $K_n$ with $\partial K_n \simeq S^{n-3}$ known as \tdef{associahedra}, whose geometry and combinatorics are of independent interest. 
\end{rem}

The original and prototypical example of an $\AA_\oo$-space is the loop space $\Omega X$ of a pointed connected space $X$. Multiplication is given by concatenation of loops. When using an actual \textit{topological space} as a model, concatenation is not strictly associative, owing to discrepancies in parametrization, but its failure is controlled by explicit coherent higher homotopies. Importantly, an $\AA_\oo$-structure can be transported along any \textit{homotopy equivalence}. And conversely, such a structure on a space $Y$, under mild hypotheses, exhibits it, up to homotopy equivalence, as a loop space $Y \simeq \Omega X$ (see \Cref{Grp_Loop} for the precise statement).
For this reason, such models of homotopy-coherent associative multiplication were classically viewed as \textit{recognition principles} for loop spaces and came to be known as \textit{loop-space machines}.    

\begin{rem}
    The tensor product of a \textit{monoidal category} is required to satisfy the Stasheff pentagon identity for fourfold products. It is then a theorem that any two groupings of an $n$-fold tensor product are canonically isomorphic for all $n$. A monoidal category can be seen as a homotopy-coherent monoid in the $(2,1)$-category of 1-categories. While its mapping spaces are not discrete, they are \textit{1-truncated} and so the obstructions $\gamma_n$ vanish automatically for all $n\ge 5$. Similarly, for a \textit{monoidal 2-category}, one would have to consider only \textit{fivefold products}.\footnote{A systematic study of the connection between \textit{homotopy truncation} and \textit{arity truncation} for more general algebraic structures can be found in~\cite{barkan2022arity}.}
\end{rem}

\subsection{Composition tables}

While it is sometimes possible to build explicit hierarchies of obstructions governing homotopy-coherent refinements, for algebraic gadgets more general than monoids it tends to be a rather complicated business. The more essential point is that the noncanonical choice of generators, relations, relations among relations etc.\! makes this approach unwieldy even when it is available. For one thing, it is not immediately clear how to define the \textit{$\oo$-category} of such things.

A different route would be to consider the \textit{entirety} of $n$-ary operations generated under composition from the multiplication and unit together with their full \textit{composition table}.
Although it might seem forbiddingly complicated at first glance, it turns out we already have access to this data in the guise of morphisms between the \textit{finitely generated free} monoids. 

\newcommand{\ff}{\mathrm{ff}}
\begin{defn}
    We denote by $\mdef{F(n)}$ the free monoid on $n$ elements and by $\Mon^{\mdef{\ff}} \sseq \Mon$ the full subcategory on the $F(n)$-s.    
\end{defn}
 For a monoid $M$, a sequence $x_1,\dots,x_n \in M$ determines (and is determined by) a homomorphism 
\[
    x_1,\dots,x_n \colon 
    F(n)\too M \qin \Mon,
\]
under which a formal word $w\in F(n)$ evaluates to $w(x_1,\dots,x_n) \in M$. Dually, a word $w\in F(n)$ corresponds to a map
\[
    w^* \colon M^n \too M 
    \qin \Set,
\]
substituting $x_1,\dots,x_n \in M$ into $w$. Combining these two facts, we can identify maps of the form $F(1) \to F(n)$ with $n$-ary operations $M^n \to M$ by considering the composition
\[
    F(1) \oto{\quad w \quad} 
    F(n) \oto{x_1,\dots,x_n} 
    M.
\]
We can further contemplate how such operations compose
\[
    M^n \oto{w_1^*,\dots,w_m^*} M^m \oto{\quad u^*\quad }  M. 
\]
The composite should correspond to the single word obtained by formally substituting the words $w_1,\dots,w_m \in F(n)$ into the $m$ variables of $u \in F(m)$,
\[
    u(w_1,\dots,w_m) \qin F(n).
\]
The categorical reason that this is indeed so is the associativity of the composition
\[
    F(1) \oto{\quad u \quad} 
    F(m) \oto{w_1,\dots,w_m}
    F(n) \oto{x_1,\dots,x_n}
    M \qin \Mon.
\]
Thus, $\Mon^\ff$ records the complete calculus of monoid operations: the morphisms are (lists of) operations, and their composition is substitution. To realize this formal ``composition table'' as an actual diagram of $n$-ary operations on a particular monoid $M$, we simply apply the functor 
\[
    \hom(-,M) \colon 
    \Mon^{\ff,\op} \too 
    \Set.
\]
Such functors have the distinctive feature of taking the $n$-fold coproduct $F(n)=F(1)^{\sqcup n}$ to the $n$-fold product 
\[
    \hom(F(n),M) \simeq 
    \hom(F(1)^{\sqcup n},M) \simeq 
    \hom(F(1),M)^n \simeq
    M^n
    \qin \Set.
\]
This turns out to be a \textit{complete characterization} of them.

\begin{prop}
    The restricted Yoneda along $\Mon^\ff\into \Mon$ is an equivalence of categories onto the product-preserving presheaves:\footnote{This is more than characterizing monoids as certain presheaves on $\Mon^\ff$, as it also identifies the \textit{homomorphisms} between them as precisely the natural transformations.}
    \[
        \Mon \iso \Fun^\times(\Mon^{\ff,\op}, \Set).
    \]
\end{prop}

Indeed, for a product-preserving functor $X$, the set $M := X(F(1))$ determines all other values $X(F(n)) \simeq M^n$,
while the obvious maps (given by $x_1\mapsto 1$ and $y_1\mapsto x_1x_2$) 
\[
    F(1) \too F(0) \quad\text{and}\quad 
    F(1) \too F(2) 
\]
recover the unit and multiplication, whose axioms are then enforced by composition in~$\Mon^\ff$.

Crucial for our purposes is that this description remains valid with $\Set$ replaced by \textit{any} category $\cC$ with products: monoid objects in $\cC$ are equivalent to product-preserving functors from $\Mon^{\ff,\op}$. The point is that, for \textit{$\oo$-categories}, this can be taken as the \textit{definition}. 

\begin{defn}
    Let $\cC$ be an $\oo$-category with finite products.\footnote{It is not strictly speaking a necessary assumption, but it makes the construction more well behaved.}
    The $\oo$-category of \tdef{monoids} in $\cC$ is given by
    \[
        \mdef{\Mon}(\cC) := \Fun^\times(\Mon^{\ff,\op}, \cC).
    \]
    When $\cC = \Spc$ we simply refer to these as \tdef{homotopy-coherent monoids}.\footnote{We shall use this convention for other algebraic structures without further comment.} 
\end{defn}

At the risk of stating the obvious, we point out that even though $\Mon^{\ff,\op}$ is an \textit{ordinary} category, functors from it to an $\infty$-category preserve composition only up to \textit{coherent homotopy}. The fact that various groupings of $n$-fold products are identified in the source $\Mon^{\ff,\op}$ explicitly, rather than being mere consequences of threefold associativity, is what creates an elaborate system of higher coherence homotopies in the \textit{target} rendering the multiplication in a monoid object of $\cC$ associative in the full homotopy-coherent sense. 

\begin{rem}
    Objects of $\Mon(\Spc)$ are also called \tdef{$\EE_1$-spaces}. More generally, a product preserving functor $\Mon^{\ff,\op} \to \cC$ taking $F(1)$ to $X \in \cC$ is referred to as an \tdef{$\EE_1$-structure} on $X$. We shall touch upon $\EE_n$-structures in \Cref{En_Mon}. 
\end{rem}

\subsection{Segal simplicial spaces}

Our model for homotopy-coherent monoids using $\Mon^\ff$ enjoys a certain air of canonicity, but, as canonical models often do, it also suffers from \textit{redundancy}, making it rather ponderous in practice.
To encode associative and unital multiplication one needn't bother with \textit{all} possible operations generated by it. It suffices to encode only the \textit{partial groupings}.  

\begin{defn}\label{Gap}
    Let 
    \[
        \mdef{\Gap} \colon \Delta \too \Mon^\ff \qin \Cat
    \]
    be the functor taking a monotone map $\alpha \colon [n] \to [m]$ to the homomorphism $F(n) \to F(m)$ defined in terms of the generators $y_1,\dots,y_n \in F(n)$ and $x_1,\dots,x_m \in F(m)$ by 
    \[
        y_i \longmapsto \prod_{j=\alpha(i-1)+1}^{\alpha(i)}x_j \qin F(m)
        \qquad\text{for}\qquad i=1,\dots,n.
    \]
\end{defn}

\begin{example}
    The following illustrates the value of $\Gap$ on a certain map $\alpha \colon [3]\to [5]$ in~$\Delta$,
    \[\begin{tikzcd}[column sep = tiny]
    	\textcolor{rgb,255:red,36;green,143;blue,36}{{[3]}} &&&& \textcolor{rgb,255:red,36;green,143;blue,36}{0} & \textcolor{rgb,255:red,33;green,72;blue,171}{{y_1}} & \textcolor{rgb,255:red,36;green,143;blue,36}{1} & \textcolor{rgb,255:red,33;green,72;blue,171}{{y_2}} & \textcolor{rgb,255:red,36;green,143;blue,36}{2} & \textcolor{rgb,255:red,33;green,72;blue,171}{{y_3}} & \textcolor{rgb,255:red,36;green,143;blue,36}{3} &&&& \textcolor{rgb,255:red,33;green,72;blue,171}{{F(3)}} \\
    	\\
    	\textcolor{rgb,255:red,36;green,143;blue,36}{{[5]}} && \textcolor{rgb,255:red,36;green,143;blue,36}{0} & \textcolor{rgb,255:red,33;green,72;blue,171}{{x_1}} & \textcolor{rgb,255:red,36;green,143;blue,36}{1} & \textcolor{rgb,255:red,33;green,72;blue,171}{{x_2}} & \textcolor{rgb,255:red,36;green,143;blue,36}{2} & \textcolor{rgb,255:red,33;green,72;blue,171}{{x_3}} & \textcolor{rgb,255:red,36;green,143;blue,36}{3} & \textcolor{rgb,255:red,33;green,72;blue,171}{{x_4}} & \textcolor{rgb,255:red,36;green,143;blue,36}{4} & \textcolor{rgb,255:red,33;green,72;blue,171}{{x_5}} & \textcolor{rgb,255:red,36;green,143;blue,36}{5} && \textcolor{rgb,255:red,33;green,72;blue,171}{{F(5)}}
    	\arrow["\alpha"', color={rgb,255:red,36;green,143;blue,36}, from=1-1, to=3-1]
    	\arrow[draw={rgb,255:red,36;green,143;blue,36}, from=1-5, to=3-5]
    	\arrow[draw={rgb,255:red,36;green,143;blue,36}, from=1-7, to=3-7]
    	\arrow[draw={rgb,255:red,36;green,143;blue,36}, from=1-9, to=3-7]
    	\arrow[draw={rgb,255:red,36;green,143;blue,36}, from=1-11, to=3-13]
    	\arrow["{\Gap(\alpha)}", color={rgb,255:red,33;green,72;blue,171}, from=1-15, to=3-15]
    	\arrow[draw={rgb,255:red,33;green,72;blue,171}, dotted, no head, from=3-6, to=1-6]
    	\arrow[draw={rgb,255:red,33;green,72;blue,171}, dotted, no head, from=3-8, to=1-10]
    	\arrow[draw={rgb,255:red,33;green,72;blue,171}, dotted, no head, from=3-10, to=1-10]
    	\arrow[draw={rgb,255:red,33;green,72;blue,171}, dotted, no head, from=3-12, to=1-10]
    \end{tikzcd}\]
    So, $\Gap(\alpha)$ takes the generators $y_1,y_2,y_3$ of $F(3)$ to the words $x_2,1,x_3x_4x_5$ in $F(5)$. 
\end{example}

Importantly, the so-called \tdef{Segal maps} $F(1) \to F(n)$, picking a generator and jointly exhibiting $F(n)$ as the coproduct $F(1)^{\sqcup n}$, are in the image of~$\Gap$. They correspond to the spine edges 

\[
    [1] \simeq \{i,i+1\} \intoo [n] \qin \Delta.
\]
With this, we can write down a leaner and more agile model for homotopy-coherent monoids:

\begin{prop}\label{Segal_Mon}
    Let $\cC$ be an $\oo$-category with products. Restriction along 
    \[
        \Gap^\op \colon 
        \Delta^\op \too 
        \Mon^{\ff,\op}
    \]
    induces an equivalence of $\oo$-categories:
    \[
        \Gray{\Mon(\cC) :=}
        \Fun^\times(\Mon^{\ff,\op},\cC) \iso 
        \Fun^{\mdef{\seg}}(\Delta^\op,\cC) 
        \Gray{\sseq \Fun(\Delta^\op,\cC),}
    \]
    where `seg' indicates the full subcategory on the functors $M$ satisfying the \tdef{Segal condition}: 
    for each $[n] \in \Delta$, the Segal maps jointly induce an isomorphism
    \[
        M([n]) \iso M([1])^n \qin \cC.
    \]
\end{prop}

In simple terms, a monoid object in $\cC$ can be equivalently encoded by a \textit{simplicial diagram} of the form:
\[\begin{tikzcd}
	{\dots } & {M\times M} & M & \pt & \cC.
	\arrow[shift left=3, from=1-1, to=1-2]
	\arrow[shift right=3, from=1-1, to=1-2]
	\arrow[shift left, from=1-1, to=1-2]
	\arrow[shift right, from=1-1, to=1-2]
	\arrow[between={0.2}{0.8}, from=1-2, to=1-1]
	\arrow[shift left=2, between={0.2}{0.8}, from=1-2, to=1-1]
	\arrow[shift right=2, between={0.2}{0.8}, from=1-2, to=1-1]
	\arrow[from=1-2, to=1-3]
	\arrow[shift left=2, from=1-2, to=1-3]
	\arrow[shift right=2, from=1-2, to=1-3]
	\arrow[shift left, between={0.2}{0.8}, from=1-3, to=1-2]
	\arrow[shift right, between={0.2}{0.8}, from=1-3, to=1-2]
	\arrow[shift right, from=1-3, to=1-4]
	\arrow[shift left, from=1-3, to=1-4]
	\arrow[between={0.2}{0.8}, from=1-4, to=1-3]
	\arrow["{\in }"{marking, allow upside down}, draw=none, from=1-4, to=1-5]
\end{tikzcd}\]

Intuitively, at degree $n$ we have $n$-tuples of elements in $M$. The \textit{degeneracy} maps correspond to inserting a unit element. The \textit{inner face} maps ($0<i<n$) correspond to multiplying together an adjacent pair of elements, and the \textit{outer face} maps ($i=0, n$) correspond to discarding the leading or trailing element. Finally, the compositions in the diagram encode the associativity, unitality, and their higher coherence. This is the most common model for homotopy-coherent monoids.

\begin{defn}\label{Bar_Construction}
    This simplicial diagram is denoted~$\mdef{\BB} M$ and is known as the \tdef{bar construction} on~$M$.    
\end{defn}

\begin{proof}[Proof sketch]
    Let us write $j := \Gap^\op$. Restriction and right Kan extension along $j$ form an adjunction on the (unrestricted) functor $\oo$-categories\footnote{We are not too careful with \textit{existence of limits} here. We shall harmlessly write as though the relevant right Kan extensions exist; the computation below shows that they do whenever needed.}
    \[
        j^* \colon 
        \Fun(\Mon^{\ff,\op}, \cC) \adj
        \Fun(\Delta^\op, \cC)
        \noloc j_*.
    \]
    Given a \textit{finite-product preserving} functor $M\colon \Mon^{\ff,\op} \to \cC$, its restriction along $j$
    \[
        \Delta^\op \oto{\ \ j\ \ } 
        \Mon^{\ff,\op} \oto{\ \ M\ \ }
        \cC
    \]
    takes the cone formed by the Segal maps $F(1)\from F(n)$ in $\Mon^{\ff,\op}$ to a limit cone in $\cC$ and hence $j^*M$ satisfies the Segal condition. Thus, $j^*$ restricts to the full subcategories
    \[
        j^* \colon \Fun^{\times}(\Mon^{\ff,\op}, \cC) \too 
        \Fun^{\seg}(\Delta^\op, \cC).  
    \]
    The first pleasant surprise here is that $j_*$ likewise restricts to these full subcategories without any further ``sheafification''. In other words, for a Segal functor $M\colon \Delta^\op \to \cC$, the right Kan extension $j_*M \colon \Mon^{\ff,\op}\to \cC$ automatically preserves finite products. The second is that $j_*$ (as well as $j^*$) does not change the value of $M$ on \textit{objects}. This implies that the (co)unit maps of the restricted adjunction are \textit{isomorphisms}, hence we get an equivalence of $\oo$-categories.
    We shall sketch the argument. 

    The composition with the Yoneda embedding 
    \[
        \Mon^\ff \intoo 
        \Fun(\Mon^{\ff,\op},\Spc)  \oto{\ j^* \ }
        \Fun(\Delta^\op,\Spc)
    \]
    takes the monoid $F(n)$ to the \textit{simplicial nerve} of the corresponding one-object category $BF(n)$ (note this is a completely 1-categorical move),
    \[
        NBF(n) \colon 
        \Delta^\op \too 
        \Set \sseq \Spc.
    \]
    By the general yoga of Kan extensions, for every test object $Z \in \cC$ we have
    \[
        \Map_\cC(Z, (j_*M)(F(n))) \simeq
        \Map_{\Fun(\Delta^\op,\Spc)}(NBF(n), \Map(Z,M)).
    \]
    Both the source and the target in the right mapping space lie in the full subcategory on \textit{Segal} functors $\Delta^\op \to \Spc$. The nontrivial ingredient needed to finish the computation is that $NBF(n)$ is the \textit{Segalification} of the simplicial space $\bigvee_n \Delta^1/\partial\Delta^1$, which follows from the \textit{necklace  formula} for Segalification in \cite[Theorem A]{barkan2025segalification}. Namely, in the $\oo$-category of Segal simplicial spaces, $NBF(n)$ is freely generated by $n$ elements in degree $1$.
    Given that, we can proceed as follows:
    \[
        \Map_{\Fun(\Delta^\op,\Spc)}(NBF(n), \Map(Z,M)) \simeq 
        \Map(Z,M_1)^n \simeq 
        \Map(Z,M_1^n).
    \]
    The Yoneda lemma then implies
    \[
        (j_*M)(F(n)) \simeq M_1^n \qin \cC. 
    \]
    We have therefore obtained an adjunction between the indicated full subcategories. Evaluation at $F(1)$ on the left and at $[1]$ on the right are both conservative. Under the canonical equivalence
    \[
        (j_*M)(F(1)) \simeq M_1 \qin \cC,
    \]
    the underlying maps of both the unit and counit are identities. They are therefore isomorphisms.
\end{proof}

\subsection{Unitality}

There is an obvious \textit{nonunital} analogue of the above story. 

\begin{defn}
    Let $\Mon^{\mdef{\rm nu},\ff}$ be the category of finitely generated free nonunital monoids. The $\oo$-category of \tdef{nonunital monoids} in an $\oo$-category $\cC$ with finite products is given by
    \[
        \Mon^{\mdef{\rm nu}}(\cC) := 
        \Fun^\times(\Mon^{\rm nu,\ff,\op}, \cC).
    \]
\end{defn}

There is also an obvious \textit{Segal model}. As the unit of the monoid structure is encoded  by the degeneracy maps, there is an equivalence between \textit{nonunital} monoid objects in $\cC$ and \textit{semi}simplicial objects $\Delta_\inj^\op \to \cC$ satisfying the \textit{Segal condition},

\[\begin{tikzcd}
	{\dots } & {M\times M} & M & \pt & \cC.
	\arrow[shift right=3, from=1-1, to=1-2]
	\arrow[shift right, from=1-1, to=1-2]
	\arrow[shift left, from=1-1, to=1-2]
	\arrow[shift left=3, from=1-1, to=1-2]
	\arrow[from=1-2, to=1-3]
	\arrow[shift left=2, from=1-2, to=1-3]
	\arrow[shift right=2, from=1-2, to=1-3]
	\arrow[shift right, from=1-3, to=1-4]
	\arrow[shift left, from=1-3, to=1-4]
	\arrow["{\in }"{marking, allow upside down}, draw=none, from=1-4, to=1-5]
\end{tikzcd}\]
The semisimplex category $\Delta_\inj$ is a \textit{directed Reedy category},\footnote{A good source for Reedy categories in $\oo$-categories is \cite{haine2025fully}.} which roughly means that morphisms only ``go up''. Diagrams indexed by such categories are amenable to an inductive analysis. Concretely, for all $n\ge 2$, a solution to the extension problem 
\[\begin{tikzcd}[column sep = large, row sep = 2.5em]
	{\Delta_{\inj, \le n-1}^\op} & \cC \\
	{\Delta_{\inj,\le n}^\op}
	\arrow["X", from=1-1, to=1-2]
	\arrow[hook, from=1-1, to=2-1]
	\arrow[dashed, from=2-1, to=1-2]
\end{tikzcd}\]
is equivalent to providing an object $X_n \in \cC$ with a map to the \tdef{matching object} of the $(n-1)$-skeleton~$X$,
\[
    X_n \too \mdef{M_n} X := \invlim_{[k]\subsetneq [n]} X_k. 
\]
If we impose the Segal condition, we have $X_k = X_1^k$ for all $k=1,\dots,n$. Intuitively, as $M_nX$ is now a limit of $X_1$-s,  all the information lives in the space of maps $X_n=X_1^n \to X_1^{}$. With a considerable amount of work, the extension problem translates to none other than a null-homotopy for the map
\[
    \gamma_n\colon S^{n-3} \too \Map(X_1^n,X_1^{})
\]
from the definition of a (nonunital) $\AA_n$-structure. A detailed account of this translation can be found in \cite[\S 4.1]{HA}. 

It turns out that for a nonunital monoid $M$, the space of ways to lift it to a (unital) monoid is either \textit{empty} or \textit{contractible}. Moreover, for the latter, it suffices for $M$ to have a unit \textit{up to homotopy} in which case $M$ is called \textit{quasi-unital}. 
This is not quite enough to qualify unitality as a \textit{property} though. The issue is that a morphism of quasi-unital nonunital monoids\footnote{Given that `mono' in `monoids' stands for having a \textit{unit}, the above literally reads as objects with an associative multiplication ``with a unit without a unit with a unit.''} need not preserve the (essentially unique) unit. This, however, \textit{is} a property for the morphism. We can summarize the situation as follows:

\begin{prop}\label{Mon_nu}
    Let $\cC$ be an $\oo$-category with finite products.
    The forgetful functor
    \[
        \Mon(\cC) \too \Mon^{\rm nu}(\cC)
    \]
    is a (nonfull) subcategory inclusion.\footnote{See \cite[\S 5.4.3]{HA} for a detailed discussion.}  
\end{prop}

\subsection{Animation}
\newcommand{\Anim}{\mathrm{Ani}}
\newcommand{\UU}{\mathbb{U}}

Having discussed the case of monoids in some detail, we now step back and extract the general mechanism behind it.
As a matter of terminology, the category $\Mon^\ff$ is the \textit{algebraic theory} of monoids, and product-preserving functors $\Mon^{\ff,\op} \to \cC$ are the \textit{models} of that algebraic theory in $\cC$.
This method of adapting 1-categorical structures to $\oo$-categories is quite general and goes nowadays by the name \textit{animation}.\footnote{It is derived from the term \textit{anima} for spaces (the animation of sets). An alternative, more classical, term is \textit{nonabelian derived category}.} As for what exactly an algebraic theory is, we adopt the following absurdly general definition:

\begin{defn}
    An \tdef{algebraic theory} is a category $\TT$ with finite coproducts.\footnote{It is common to restrict to the case where $\TT$ is generated  under coproducts by a \textit{single object}, which is indeed the case in most (but not all) of the examples we'll consider.} A \tdef{model} of $\TT$ in an $\oo$-category $\cC$ is a product preserving functor $\TT^\op \to \cC$. 
\end{defn}
The category of models is contravariantly functorial in the theory and covariantly functorial in the target. A finite \textit{coproduct} preserving functor $\UU\to\TT$ and a finite \textit{product} preserving functor $\cC\to \cD$ induce, by pre- and postcomposition respectively,
\[
    \Fun^\times(\TT^{\op}, \cC) \too \Fun^\times(\UU^{\op}, \cC) \qquad,\qquad
    \Fun^\times(\TT^{\op}, \cC) \too \Fun^\times(\TT^{\op}, \cD).
\]
The forgetful functor from monoids to nonunital monoids illustrates the first, while the ordinary monoid structure on the connected components of a homotopy-coherent monoid is an example of the second. The first is always limit preserving and thus admits a left adjoint under general presentability assumptions. The same holds for the second, if $\cC \to \cD$ is limit (and not only finite product) preserving.  

Whenever we can realize a given category as models of some algebraic theory in \textit{sets}, we can consider models of this same theory in any $\oo$-category, and most importantly, in \textit{spaces}. 

\begin{defn}
    Let $\cC$ be the category of models in $\Set$ for some algebraic theory $\TT$. 
    The \tdef{animation} of $\cC$ is the $\oo$-category
    \[
        \mdef{\Anim}(\cC) := \Fun^\times(\TT^\op, \Spc) \qin \Cat_\oo.
    \]
\end{defn}

There are many interesting $\infty$-categories that arise in this fashion. We begin with one from classical homological algebra.
For a ring $R$, the finitely generated free (say, right) $R$-modules are of the form~$R^n$. The set of morphisms $R^m \to R^n$ can be identified with the set of matrices $M_{n\times m}(R)$. Thus, the algebraic theory of $R$-modules $\Mod_R^\ff$ is equivalent to the following category: 

\begin{defn}
    Let $\mdef{\Mat}(R)$ be the category whose objects are natural numbers, morphisms are matrices over $R$ and composition is given by matrix multiplication. 
\end{defn}

For a reader who might be somewhat put off by this whole categorical formalism, because it feels far removed from the way one usually thinks about algebraic structures, the example of $\Mat(R)$ can help bridge that gap. One usually defines a vector space $V$ over a field $\FF$ in terms of a pair of operations, vector addition and scalar multiplication, subject to a list of axioms. However, one then usually proceeds to develop the theory almost exclusively in terms of general \textit{linear combinations} of the form
\[
    a_1 v_1 + a_2v_2 +\dots +a_nv_n  \qin V.
\]    
Thus, thinking of a vector space structure in terms of $n$-ary operations $V^n \to V$ parameterized by~$\FF^n$, composing according to the rules of matrix multiplication, is perhaps not so foreign. 

\begin{example}
    The animation of $\Mod_R$ is the $\infty$-categorical (connective) \textit{derived category} of $R$,
    \[
        \Anim(\Mod_R) :=
        \Fun^\times(\Mat(R)^\op, \Spc) \simeq 
        \sD_{\ge 0}(R).
    \]    
    Note that this method of constructing the derived category avoids any explicit mention of chain complexes and their homology. Furthermore, a ring homomorphism $R \to S$ induces a coproduct preserving functor $\Mat(R) \to \Mat(S)$ restriction along which yields \textit{restriction of scalars} 
    \[
        \sD_{\ge_0}(S) \too \sD_{\ge 0}(R).
    \]
    Its left adjoint is given by the appropriately ``derived'' \textit{extension of scalars}. 
    
\end{example}

Adding \textit{multiplicative structure} to the previous example produces a more intricate object.

\begin{example}\label{Ani_CRing}
    The algebraic theory of commutative rings is the category of finitely generated polynomial rings $\ZZ[x_1,\dots,x_n]$. The animation of the category of commutative rings is called simply the $\infty$-category of \tdef{animated commutative rings}.\footnote{This notion appears, and plays a key role, in \cite{cesnavicius2019purity}.} 
\end{example}

Our description of animation as a process of the form
\[
    \cC \simeq \Fun^\times(\TT^\op, \Set) \longmapsto \Fun^\times(\TT^\op, \Spc) =: \Anim(\cC)
\]
glossed over the fact that $\TT$ is determined by $\cC$ only up to \textit{idempotent completion}. This ends up having no effect on $\Anim(\cC)$, which is insensitive to this ambiguity. It is nevertheless useful to describe the process starting from $\cC$ itself rather than from a choice of $\TT$. It is similarly useful to be able to recognize when a given $\oo$-category is the animation of some $\cC$. To do so, we must briefly shift gears and plunge into some rather abstract category theory.

Just as the presheaf category is the free cocompletion of a category, the Yoneda embedding
\[
    \TT \intoo \Fun^\times(\TT^\op, \Set) := \mdef{\sP_\Sigma}(\TT; \mdef{\Set})
\]
exhibits its target as the 1-categorical free \tdef{sifted cocompletion} of $\TT$: it is obtained by freely adjoining sifted colimits, which one may think of roughly as filtered colimits and geometric realizations\footnote{For 1-categories, these are equivalent to \textit{reflexive coequalizers}. }. The extent to which this process can be reversed relies on the following notion:

\begin{defn}
    An object $X$ of an ordinary category $\cC$ is called \tdef{compact-projective} if the functor it corepresents  
    \[
        \Hom(X,-)\colon \cC \too \Set
    \]
    preserves sifted colimits. 
    We denote by $\cC^{\mdef{\cp}} \sseq \cC$ the full subcategory on compact-projective objects.  
\end{defn}

It is easy to see that the Yoneda image of $\TT$ in $\sP_{\Sigma}(\TT;\Set)$ consists of compact-projective objects and that the latter are closed under retracts. It is similarly easy to deduce that $\sP_{\Sigma}(\TT;\Set)^\cp$ is precisely the \textit{idempotent completion} of $\TT$. 

\begin{example}
    For $\TT = \Mat(R)$, this is the difference between \textit{free} and \textit{projective} finitely generated $R$-modules. They both have the same free sifted cocompletion, the $R$-modules.
\end{example}

The $\oo$-categorical Yoneda embedding
\[
    \TT \intoo \Fun^\times(\TT^\op, \Spc) 
    := \mdef{\sP_\Sigma}(\TT)
\]
enjoys the analogous free sifted cocompletion universal property in the world of $\oo$-categories. Thus, one can think of the animation of an ordinary category $\cC$ as extracting the compact-projective objects $\cC^\cp$ of which everything else is built by formal sifted colimits, followed by re-adjoining the similarly formal sifted colimits but in the \textit{$\oo$-categorical} world. 

\begin{rem}
    The inclusion $\Set \into \Spc$ preserves products, so it induces a fully faithful inclusion
    \[
        \cC := 
        \Fun^\times(\TT^\op, \Set) \intoo 
        \Fun^\times(\TT^\op, \Spc)
        =: \Anim(\cC),
    \]
    identifying $\cC$ with the \textit{discrete objects} of $\Anim(\cC)$. It can be shown that every object of $\Anim(\cC)$ is a geometric realization of a simplicial diagram in $\cC$. 
    This explains why ``injecting homotopy'' into a category $\cC$ was historically modeled by the category of \textit{simplicial objects} in $\cC$.\footnote{For the same reason, animated rings are  still occasionally referred to, somewhat reductively, as \textit{simplicial rings}.}
    In the additive setting, the Dold-Kan correspondence identifies simplicial objects with \textit{chain complexes}, recovering the classical route to derived categories.
\end{rem}

The above discussion can be complemented by a very useful \textit{recognition principle} for $\oo$-categories of the form~$\sP_\Sigma(\TT)$. In short, the evident necessary conditions are also sufficient.

\begin{prop}\label{Ani_Recog}
    Let $\cC$ be an $\oo$-category and let $\TT \sseq \cC$ be a full subcategory of compact-projective objects,\footnote{Here we adapt the definition to the $\oo$-categorical setting by considering the representing functors into $\Spc$. This is \textit{not} compatible with the 1-categorical notion even for 1-categories.} closed under finite coproducts and generating $\cC$ under (sifted) colimits. Then, there is a canonical equivalence $\sP_\Sigma(\TT) \simeq \cC$. 
\end{prop}

\begin{proof}[Proof sketch]
    Consider the embedding $f_0\colon \TT \into \cC$. By the universal property of $\TT \into \sP_\Sigma(\TT)$, there is a unique colimit preserving extension of $f_0$,
    \[
        f\colon \sP_\Sigma(\TT) \too \cC. 
    \]
    Essential surjectivity of $f$ follows from the fact that $\TT$ generates $\cC$ under colimits. For fully faithfulness, observe that $f$ admits a right adjoint $g$ and it suffices to show that the unit $\Id \to gf$ is an isomorphism. The fact that $\TT$ consists of compact-projective generators implies that $g$ preserves sifted colimits so it suffices to check on $\TT$ where it is a consequence of $f_0$ being fully faithful.\footnote{The omitted details can be found in \cite[Proposition 5.5.8.22]{HTT}.} 
\end{proof}

    

\newcommand{\coMon}{\mathrm{coMon}}

\section{Groups}

Among the many kinds of monoids arising in mathematics, by far the most popular are groups, which appear as symmetries of objects. Homotopy-coherent groups form a full subcategory of homotopy-coherent monoids, but also admit a surprisingly geometric model --- they are exactly the loop spaces. This is, however, a special feature of the $\oo$-category $\Spc$ and cannot be used as an alternative route to group objects in arbitrary $\oo$-categories. Even in spaces, the monoid perspective remains useful: the interaction between groups and monoids which are not groups, most notably through group completion, is an important part of the story.

\subsection{Shear power}

Groups are simply monoids where every element is invertible, but one might worry that this notion too has to be handled in some subtle way in the homotopy-coherent world. Fortunately, this is not the case. An ordinary monoid $M$ is a group precisely when its \tdef{shear map}
\[
    s \colon M\times M \too M\times M,
\]
\[
    s(x,y)  = (x, x\cdot y)
\]
is a bijection. Its inverse then gives \textit{division}
\[
    s^{-1}(x,y) = (x,x^{-1} y).
\]
This shear map exists for every $M \in \Mon(\cC)$ as it is part of the algebraic theory. It is encoded by the map $F(2) \to F(2)$ that picks the elements $x, xy \in F(2)$ where $x$ and $y$ are the generators of~$F(2)$. 

\begin{defn}
    Let $\cC$ be an $\oo$-category with finite products.
    The $\infty$-category of \tdef{groups} in $\cC$ is the full subcategory 
    \[
        \mdef{\Grp}(\cC)\ \sseq\ \Mon(\cC)
    \]
    spanned by those monoids $M$ for which the shear map is an isomorphism:
    \[
        s \colon M\times M \iso M \times M.
    \]
\end{defn}

In the case of spaces, there is an even simpler criterion. 
\begin{prop}\label{CGrp_pi0}
    A monoid $M \in \Mon(\Spc)$ is a group if and only if $\pi_0M \in \Mon(\Set)$ is an ordinary group.\footnote{These are classically called \textit{grouplike }.} 
\end{prop}
In particular, if $M$ is \textit{connected}, it is always a group. 
\begin{proof}
    One direction is easy. The shear map of $\pi_0 M$ is literally the $\pi_0$ of the shear map of $M$ ($\pi_0$ commutes with products) and a homotopy equivalence induces an isomorphism on $\pi_0$. Conversely, assume $\pi_0 M$ is a group. The shear map lies over $M$ via the projection onto the first factor:
    \[\begin{tikzcd}
    	{M\times M} & {M\times M} \\
    	M & {M.}
    	\arrow["s", from=1-1, to=1-2]
    	\arrow["{p_1}", from=1-1, to=2-1]
    	\arrow["{p_1}", from=1-2, to=2-2]
    	\arrow[equals, from=2-1, to=2-2]
    \end{tikzcd}\]
    Over each $x\in M$, the fibers of both vertical maps are $M$ and the induced map $L_x\colon M\to M$ is left multiplication by $x$. Since $x$ is invertible in $\pi_0 M$, there is $y \in M$ such that $xy \sim yx \sim 1$. It follows that 
    \[
        L_xL_y \sim L_{xy} \sim \Id \quad,\quad
        L_y L_x \sim L_{yx} \sim \Id. 
    \]
    Thus, $s$ is a fiber-wise equivalence and hence an equivalence. 
\end{proof}

\subsection{Loop spaces}

It is not hard to see that $\Grp(\Spc)$ is itself the animation of the ordinary category of groups,
\[
    \Grp(\Spc) \simeq \Fun^\times(\Grp^{\ff,\op}, \Spc).
\]
From this we get another, deceptively simple, presentation of this $\oo$-category. 

\begin{thm}\label{Grp_Loop}
    We have an equivalence of $\oo$-categories
    \[
        \Grp(\Spc) \simeq \Spc_{*,\ge 1} \qin \Cat_\oo
    \]
    where on the right we have the $\oo$-category of pointed connected spaces.\footnote{This is an $\oo$-categorical take on the classical theme of modeling pointed connected spaces by \textit{simplicial groups}.} 
\end{thm}
\begin{proof}
    We use the recognition principle for algebraic theories from \Cref{Ani_Recog}. Let $\mdef{\mathcal{W}} \sseq \Spc_{*,\ge 1}$ be the full subcategory on finite \textit{wedges} of circles $\bigvee_n S^1$. 
    Since the latter is the classifying space of the free group on $n$ elements, the \textit{fundamental group} functor restricts to an equivalence of categories:
    \[
        \pi_1 \colon \mathcal{W} \iso \Grp^\ff.
    \]  
    
    The circle $S^1$ corepresents the loop-space functor $\Omega\colon \Spc_{*,\ge1} \to \Spc$. It is conservative (by looking at homotopy groups), which by standard abstract nonsense implies that $S^1$, and hence $\mathcal{W}$, generate $\Spc_{*,\ge1}$ under colimits. Furthermore, $S^1$ is compact-projective. Preservation of \textit{filtered} colimits follows from $\Omega$ itself being a finite limit. Preservation of \textit{geometric realizations} is more subtle and relies crucially on the \textit{connectedness} assumption (see \cite[Lemma 5.5.6.17]{HA}). 
    Compact-projective objects are closed under finite coproducts, hence $\mathcal{W}$ consists of compact-projective objects. 
    The general recognition principle yields 
    \[
        \Spc_{*,\ge 1} \simeq \sP_\Sigma(\mathcal{W}) \simeq \sP_\Sigma(\Grp^\ff) \simeq \Grp(\Spc).
    \]
\end{proof}

By construction, the circle $S^1 \in \Spc_{*,\ge 1}$ corresponds to the free group on one element $F(1) \in \Grp(\Spc)$ and so the functors they corepresent make the following diagram commute:
\[\begin{tikzcd}
	{\Spc_{*,\ge 1}} && {\Grp(\Spc)} \\
	& \Spc.
	\arrow[from=1-1, to=1-3]
	\arrow["{\sim }", shift right, draw=none, from=1-1, to=1-3]
	\arrow["\Omega"', from=1-1, to=2-2]
	\arrow["U", from=1-3, to=2-2]
\end{tikzcd}\]
Unpacking this a bit further, for a pointed connected space $X$ the multiplication map
\[
    m\colon \Omega X \times \Omega X \too \Omega X
\]
is given by \textit{concatenation} of loops. Indeed, it is corepresented by the pinch map $S^1 \to S^1 \vee S^1$ which on $\pi_1$ induces the map $F(1)\to F(2)$ given by $y_1 \mapsto x_1x_2$. We can summarize the situation loosely by saying that the homotopy-coherent groups in spaces are exactly the \textit{loop spaces}.  

There is a beautiful (and useful) geometric interpretation of the above equivalence in terms of Segal's \textit{simplicial space} model. Any map of spaces $Y\to X$ admits a \textit{\v{C}ech resolution}:
\[\begin{tikzcd}
	{\dots } & {Y\times_XY\times_XY} & {Y\times_XY} & Y & {\Spc.}
	\arrow[shift left=3, from=1-1, to=1-2]
	\arrow[shift right=3, from=1-1, to=1-2]
	\arrow[shift left, from=1-1, to=1-2]
	\arrow[shift right, from=1-1, to=1-2]
	\arrow[between={0.2}{0.8}, from=1-2, to=1-1]
	\arrow[shift left=2, between={0.2}{0.8}, from=1-2, to=1-1]
	\arrow[shift right=2, between={0.2}{0.8}, from=1-2, to=1-1]
	\arrow[from=1-2, to=1-3]
	\arrow[shift left=2, from=1-2, to=1-3]
	\arrow[shift right=2, from=1-2, to=1-3]
	\arrow[shift left, between={0.2}{0.8}, from=1-3, to=1-2]
	\arrow[shift right, between={0.2}{0.8}, from=1-3, to=1-2]
	\arrow[shift right, from=1-3, to=1-4]
	\arrow[shift left, from=1-3, to=1-4]
	\arrow[between={0.2}{0.8}, from=1-4, to=1-3]
	\arrow["{\in }"{marking, allow upside down}, draw=none, from=1-4, to=1-5]
\end{tikzcd}\]
This is loosely speaking the \textit{equivalence relation} induced on $Y$ by the map $Y \to X$. Formally, it is the \textit{right adjoint} of the functor 
\[
    \ev_{[0]} \colon 
    \Fun(\Delta^\op, \Spc_{/X}) \too 
    \Spc_{/X}.
\]

\begin{prop}\label{Segal_Cech}
    Let $X$ be a pointed connected space. The bar construction $\BB (\Omega X)$ identifies with the \v{C}ech resolution of the basepoint inclusion $\pt \to X$, 
    \[\begin{tikzcd}
    	{\dots } & {\Omega X\times \Omega X} & {\Omega X} & \pt & \Spc.
    	\arrow[shift left=3, from=1-1, to=1-2]
    	\arrow[shift right=3, from=1-1, to=1-2]
    	\arrow[shift left, from=1-1, to=1-2]
    	\arrow[shift right, from=1-1, to=1-2]
    	\arrow[between={0.2}{0.8}, from=1-2, to=1-1]
    	\arrow[shift left=2, between={0.2}{0.8}, from=1-2, to=1-1]
    	\arrow[shift right=2, between={0.2}{0.8}, from=1-2, to=1-1]
    	\arrow[from=1-2, to=1-3]
    	\arrow[shift left=2, from=1-2, to=1-3]
    	\arrow[shift right=2, from=1-2, to=1-3]
    	\arrow[shift left, between={0.2}{0.8}, from=1-3, to=1-2]
    	\arrow[shift right, between={0.2}{0.8}, from=1-3, to=1-2]
    	\arrow[shift right, from=1-3, to=1-4]
    	\arrow[shift left, from=1-3, to=1-4]
    	\arrow[between={0.2}{0.8}, from=1-4, to=1-3]
    	\arrow["{\in }"{marking, allow upside down}, draw=none, from=1-4, to=1-5]
    \end{tikzcd}\]
    Conversely, one recovers $X$ by taking the colimit of this simplicial diagram,\footnote{We use the traditional \textit{geometric realization} notation $|-|$ for simplicial colimits.}
    \[\begin{tikzcd}
    	{X\simeq \big|\dots } & {\Omega X\times \Omega X} & {\Omega X} & {\pt\big|.}
    	\arrow[shift left=3, from=1-1, to=1-2]
    	\arrow[shift right=3, from=1-1, to=1-2]
    	\arrow[shift left, from=1-1, to=1-2]
    	\arrow[shift right, from=1-1, to=1-2]
    	\arrow[between={0.2}{0.8}, from=1-2, to=1-1]
    	\arrow[shift left=2, between={0.2}{0.8}, from=1-2, to=1-1]
    	\arrow[shift right=2, between={0.2}{0.8}, from=1-2, to=1-1]
    	\arrow[from=1-2, to=1-3]
    	\arrow[shift left=2, from=1-2, to=1-3]
    	\arrow[shift right=2, from=1-2, to=1-3]
    	\arrow[shift left, between={0.2}{0.8}, from=1-3, to=1-2]
    	\arrow[shift right, between={0.2}{0.8}, from=1-3, to=1-2]
    	\arrow[shift right, from=1-3, to=1-4]
    	\arrow[shift left, from=1-3, to=1-4]
    	\arrow[between={0.2}{0.8}, from=1-4, to=1-3]
    \end{tikzcd}\]
\end{prop}

\begin{proof}[Proof idea]
    The fact that the \v{C}ech resolution of $\pt \to X$ is the Segal space corresponding to $\Omega X$ follows ``by inspection''. That the colimit recovers $X$ is a general fact about \textit{epimorphisms} in $\oo$-topoi. That's where the \textit{connectedness} assumption on $X$ is crucially used. 
\end{proof}

\begin{rem}
    Geometrically speaking, the group $\Omega X$ acts on the one-point space $\pt$ with a quotient map to the homotopy orbits space 
    \[
        \pt \too \pt/\!\!/\Omega X \ \simeq\  X \qin \Spc_*.
    \]
    The group $\Omega X$ can then be recovered from the ``equivalence relation'' imposed by this quotient map on $\pt$. This approach to group objects in $\Spc$ relies heavily on the latter being an \textit{$\oo$-topos} and would not work in most $\oo$-categories.     
\end{rem}

While perhaps striking at first glance, under the homotopy hypothesis, \Cref{Grp_Loop} reads as an equivalence between homotopy-coherent groups and one-object $\oo$-groupoids, which is exactly what one would expect. Indeed, for a pointed space $X$ viewed as an $\oo$-groupoid, the loop space $\Omega X$ corresponds precisely to the group of \textit{automorphisms} of the basepoint in $X$. With more effort, the same line of reasoning can be applied to show the following: 

\begin{prop}
    There is an equivalence between $\Mon(\Spc)$ and the $\oo$-category of one-object $\oo$-categories. 
\end{prop}

\subsection{Group theory}
The pointed connected space model for groups in $\Spc$ translates ``homotopy-coherent group theory'' into ordinary homotopy theory, yielding interesting reinterpretations and generalizations of classical group-theoretic phenomena. We give some examples. 

To begin with, it allows one to compute \textit{colimits} in $\Grp(\Spc)$.\footnote{Limits are easy as the forgetful functor to spaces preserves (and creates) limits.} The full subcategory $\Spc_{*,\ge 1}\sseq \Spc_*$ is closed under colimits, and in $\Spc_*$ colimits are computed as in $\Spc$ followed by collapsing the colimit of the basepoints into one. The prototypical example is the \textit{coproduct}, which is given by the \tdef{wedge sum}
\[\begin{tikzcd}
	{\pt \sqcup \pt } & {X\sqcup Y} \\
	\pt & {X\mathbin{\mdef{\vee}} Y.}
	\arrow[from=1-1, to=1-2]
	\arrow[from=1-1, to=2-1]
	\arrow[from=1-2, to=2-2]
	\arrow[from=2-1, to=2-2]
	\arrow["\lrcorner"{anchor=center, pos=0.125, rotate=180}, draw=none, from=2-2, to=1-1]
\end{tikzcd}\]
For colimits of contractible shape, such as \textit{pushouts}, the situation is even simpler: the colimit is computed as in $\Spc$.
This apparent simplicity clashes with the fact that the homotopy truncation functor
\[
    \pi_0 \colon \Grp(\Spc) \too \Grp(\Set)
\]
is colimit preserving and colimits in $\Grp$ are already rather complicated. The issue is of course that the underlying space of the corresponding group object is given by the \textit{loop space} functor
\[
    \Omega \colon \Spc_{*, \ge 1 } \too \Spc_*,
\]
whose interaction with colimits is complicated. In particular, under the equivalence $\Grp(\Spc) \simeq \Spc_{*,\ge1}$, the truncation functor corresponds to the \textit{fundamental group},
\[
    \pi_1 \colon \Spc_{*,\ge 1} \too \Grp(\Set). 
\]
Thus, for pushouts, we recover the classical Seifert--van Kampen theorem,
\[
    \pi_1(X \cup_Z Y) \simeq \pi_1(X) *_{\pi_1(Z)} \pi_1(Y) \qin \Grp(\Set).
\]
But there is more. For ordinary groups, the kernel of the canonical map 
\[
    G * H \too G\times H
\]
is always \textit{freely generated} by the basic commutators 
\[
    \{[g,h] \mid g \in G\smallsetminus\{e\}\ ,\ h \in H\smallsetminus\{e\} \}.
\]
As the left adjoint of $\Omega$ is $\Sigma$, the free group on a pointed space $X$ is modeled by the pointed connected space $\Sigma X$ (which, as a group, therefore has \textit{underlying space} $\Omega \Sigma X$\footnote{This is an $\oo$-categorical formulation of Milnor's classical simplicial model for $\Omega\Sigma$ (see \cite{milnor1972construction}).}).
The following special case of the Hilton-Milnor theorem is therefore a direct homotopy-coherent generalization of the above:
\begin{prop}
    For all $X,Y \in \Spc_{*,\ge 1}$ there is a fiber sequence,
    \[
        \Sigma(\Omega X\wedge \Omega Y) \too 
        X \vee Y \too 
        X\times Y 
        \qin \Spc_{*,\ge 1}.
    \]
\end{prop}

\begin{proof}
    The only thing required for the argument is the \textit{universality} of colimits in spaces. Namely, that for every map $f\colon A \to B$ in $\Spc$, the pullback functor
    \[
        f^* \colon \Spc_{/B} \too \Spc_{/A}
    \]
    preserves colimits. Denote by $Z$ the fiber of $X\vee Y\to X\times Y$. The whole move can be summarized in the following diagram:
    \[\begin{tikzcd}[sep = 1.3em]
    	{\Omega X\times \Omega Y} && {\Omega Y} & \\
    	& \pt && Y \\
    	{\Omega X} && {Z} \\
    	& X && {X\vee Y}
    	\arrow[from=1-1, to=1-3]
    	\arrow[dashed, from=1-1, to=2-2]
    	\arrow[from=1-1, to=3-1]
    	\arrow[dashed, from=1-3, to=2-4]
    	\arrow[from=1-3, to=3-3]
    	\arrow[from=2-2, to=2-4]
    	\arrow[from=2-2, to=4-2]
    	\arrow[from=2-4, to=4-4]
    	\arrow[from=3-1, to=3-3]
    	\arrow[dashed, from=3-1, to=4-2]
    	\arrow["\lrcorner"{anchor=center, pos=0.125, rotate=180}, draw=none, from=3-3, to=1-1]
    	\arrow[dashed, from=3-3, to=4-4]
    	\arrow[from=4-2, to=4-4]
    \end{tikzcd}\]
    The front face lives over $X\times Y$ and the back face is obtained from it by pulling back along $\pt \to X\times Y$. Since the front one is a pushout, so is the back one, but the pushout of 
    \[
        \Omega X \longleftarrow \Omega X \times \Omega Y \too \Omega Y
    \]
    is the \textit{join} of $\Omega X$ and $\Omega Y$, which is equivalent to $\Sigma(\Omega X \wedge \Omega Y)$.
\end{proof}

The pointed connected space model is also particularly well suited for discussing \textit{group actions}.  For a pointed connected space $X$, an action of the group $\Omega X$ on some object of an $\oo$-category $\cC$ is simply a functor $X \to \cC$. Restriction along the basepoint $\pt \to X$ corresponds to forgetting the action. 

The case $\cC = \Spc$ is particularly nice as we have the \textit{straightening--unstraightening} equivalence 
\[
    \Spc_{/X} \ \simeq \ \Fun(X,\Spc) \qin \Cat_\oo. 
\]
To $Y \to X$ on the left-hand side, we associate the functor $X\to \Spc$ given by the family of fibers  $x \mapsto Y_x$. Conversely, given the latter, we have
\[
    Y = \colim_{x\in X} Y_x \too \colim_{x \in X}\pt = X.
\]
The bottom line is that a homotopy-coherent action of $\Omega X$ on a space $Z$ is completely encoded in the fiber sequence  
\[
    Z \too Y \too X \quad\text{where}\quad Y = Z/\!\!/\Omega X.
\]
Furthermore, the space of \textit{homotopy fixed points} of the action of $\Omega X$ on $Z$ is then the space of \textit{sections} of the map $Y \to X$.  This is a remarkably \textit{low-coherence} encoding of the entire situation. 

\begin{example}\label{Conj_LX}
    Every homotopy-coherent group $\Omega X$ acts on itself by \textit{conjugation}. The associated fiber sequence is 
    \[
        \Omega X \too LX \too X \quad\text{where}\quad \mdef{L} X := \Map(S^1,X)
    \]
    is the \tdef{free loop space} and the map is evaluation at the basepoint of $S^1$. Incidentally, there is also an action of the group $S^1 = \Omega \CC P^\oo$ on $LX$ by \textit{rotating the loops}. 
\end{example}

A related topic is the classification of \textit{extensions}. For ordinary groups $G$ and $H$ it is a classical problem in group theory to classify all extensions of $G$ by $H$, which are short exact sequences of the form:
\[
    1 \too H \too K \too G \too 1 \qin \Grp.
\]
The corresponding notion of an extension of $\Omega X$ by $\Omega Y$ is a fiber sequence 
\[
    Y \too E \too  X \qin \Spc_{*,\ge 1},
\]
or in other words, an action of $\Omega X$ on the (unpointed) space $Y$. Such fiber sequences are classified by pointed maps 
\[
    X \too B\Aut(Y) \ \subset\ \Spc. 
\]
In the case $X = BG$ and $Y = BH$, the classifying space is a \textit{2-type} with 
\[
    \pi_1(B\Aut(BH)) =\Out(H) \quad\text{and}\quad
    \pi_2(B\Aut(BH)) = Z(H).
\]
The \textit{Postnikov invariant} is a canonical class 
\[
    \mdef{\kappa}_H \ \in\  H^3(\Out(H);Z(H)). 
\]
This leads immediately to the classical Schreier--Eilenberg--Mac Lane classification of the extensions of $G$ by $H$.
\begin{prop}
    An extension of $G$ by $H$ induces a group homomorphism
    \[
        \rho \colon G \too \Out(H),
    \]
    such that for the induced action of $G$ on $Z(H)$ we have
    \[
        \rho^*(\kappa_H) \equiv 0 \qin H^3(G;Z(H)).
    \]
    Conversely, for every homomorphism $\rho$ satisfying this condition, the equivalence classes of extensions inducing $\rho$ form a torsor under $H^2(G;Z(H))$.
\end{prop}
For general spaces $X$ and $Y$, the Postnikov tower of $B\Aut(Y)$ similarly yields an obstruction theory for such extensions and, under the usual convergence hypotheses, an unstable spectral sequence for the homotopy groups of the corresponding mapping space.

\subsection{Group completion}

The condition defining groups among monoids is that a certain natural map is an equivalence. Somewhat unexpectedly, the full subcategory it defines is closed under both \textit{limits} and \textit{colimits}.

\begin{prop}
    The inclusion 
    \[
        \Grp(\Spc) \intoo \Mon(\Spc)
    \]
    preserves all limits and colimits.
\end{prop}
\begin{proof}
    Preservation of limits is immediate. Indeed, we have a commutative diagram
    \[\begin{tikzcd}
    	{\Grp(\Spc)} && {\Mon(\Spc)} \\
    	& \Spc
    	\arrow[hook, from=1-1, to=1-3]
    	\arrow["{U_{\Grp}}"', from=1-1, to=2-2]
    	\arrow["{U_{\Mon}}", from=1-3, to=2-2]
    \end{tikzcd}\]
    where both forgetful functors are limit preserving and conservative. Preservation of colimits is less immediate, at least from abstract considerations. However, by \Cref{CGrp_pi0}, being a group can be detected on $\pi_0$. Namely we have a pullback diagram of categories,
    \[\begin{tikzcd}[sep = 2.3em]
    	{\Grp(\Spc)} & {\Mon(\Spc)} \\
    	\Grp & \Mon.
    	\arrow[hook, from=1-1, to=1-2]
    	\arrow["{\pi_0}"', from=1-1, to=2-1]
    	\arrow["\lrcorner"{anchor=center, pos=0.125}, draw=none, from=1-1, to=2-2]
    	\arrow["{\pi_0}", from=1-2, to=2-2]
    	\arrow[hook, from=2-1, to=2-2]
    \end{tikzcd}\]
    The functor $\pi_0$ preserves colimits and the ordinary subcategory inclusion $\Grp \sseq \Mon$ is known to be closed under colimits. 
\end{proof}

Consequently, by the adjoint functor theorem, the inclusion admits both adjoints. As with ordinary monoids and groups, the right adjoint 
\[
    (-)^{\mdef{\times}} \colon 
    \Mon(\Spc) \too 
    \Grp(\Spc)
\]
is straightforward; one simply discards the (connected components of) noninvertible elements, producing the \tdef{group of units} of $M$. We thus have a pullback diagram in $\Mon(\Spc)$:
\[\begin{tikzcd}[sep = 2.3em]
	{M^\times} & M \\
	{(\pi_0 M)^\times } & {\pi_0 M.}
	\arrow[hook, from=1-1, to=1-2]
	\arrow[from=1-1, to=2-1]
	\arrow["\lrcorner"{anchor=center, pos=0.125}, draw=none, from=1-1, to=2-2]
	\arrow[from=1-2, to=2-2]
	\arrow[hook, from=2-1, to=2-2]
\end{tikzcd}\]
The left adjoint, the \tdef{group completion},   
\[
    (-)^{\mdef{\gp}} \colon 
    \Mon(\Spc) \too 
    \Grp(\Spc),
\]
is more interesting. On the level of connected components we recover the ordinary group completion 
\[
    \pi_0(M^\gp) \simeq (\pi_0 M)^\gp \qin \Grp.
\]
To describe the full space-level group completion, we again make use of the equivalence of \Cref{Grp_Loop}, or rather its \v{C}ech resolution interpretation in \Cref{Segal_Cech}. Recall that every $X \in \Spc_{*,\ge 1}$ is the colimit of the bar construction on the corresponding loop space $\Omega X$. But note that we can take the colimit of the bar construction on \textit{any} homotopy-coherent monoid~$M$, 
\[
    \mdef{B}M  := |\BB M| \qin \Spc. 
\]

\begin{prop}\label{Bar_Grp_Cpl}
    The composition
    \[
        \Omega B  \colon 
        \Mon(\Spc) \too  \Grp(\Spc)
    \]
    is naturally equivalent to the group completion functor $(-)^\gp$.
\end{prop}

\newcommand{\redu}{\mathrm{red}}

\begin{proof}[Proof sketch]
    By standard arguments, the \v{C}ech resolution of $\pt \to X$ is the right adjoint to the simplicial colimit functor 
    \[
        |-|\colon \Fun^{\mdef{\redu}}(\Delta^\op, \Spc) \too \Spc_*,
    \]
    where the superscript `red' refers to the full subcategory on \textit{reduced} simplicial objects $X$ (ie those with $X_0 =\pt$). This adjunction restricts to the full subcategory $\Fun^{\seg}(\Delta^\op,\Spc)$ on the left and $\Spc_{*,\ge1}$ on the right. It follows that 
    \[
        \Mon(\Spc) \isob[\ \ \BB\ \ ]
        \Fun^{\seg}(\Delta^\op,\Spc) \oto{\ |-|\ }
        \Spc_{*,\ge 1} \isob[\ \ \Omega\ \ ]
        \Grp(\Spc)
    \]
    is the left adjoint of the inclusion $\Grp(\Spc) \into \Mon(\Spc)$, and is therefore the group-completion functor.
\end{proof}

\begin{rem}
    For an ordinary group $G$, the space $BG$ is an Eilenberg--Mac Lane space $K(G,1)$. However, for an ordinary \textit{monoid} $M$, the space $BM$ can be very complicated. In fact, \textit{every} pointed connected space can be realized as $BM$ for some ordinary monoid $M$.\footnote{This is classical \cite{mcduff1979classifying}.} Hence, while the fundamental group of $BM$ is the ordinary group completion of $M$, the higher homotopy groups of $BM$ can be viewed as its \textit{derived functors}.\footnote{Indeed, they can be defined/computed in a more classical fashion using levelwise free simplicial resolutions.} 
\end{rem}

\begin{example}
    Consider the ordinary 5-element monoid\footnote{This nice concrete example can be found with details in \cite{fiedorowicz2002counterexample}.}
    \[
        M = \{1,x_{11},x_{12},x_{21},x_{22}\} \quad\text{with}\quad
        x_{ij}\cdot x_{k\ell} = x_{i\ell}.
    \]
    Since every element of $M$ is \textit{idempotent}, the ordinary group completion of $M$ is the trivial group. However, it turns out that $BM \simeq S^2$ and therefore
    \[
        M^\gp \simeq \Omega S^2 \qin \Grp(\Spc).
    \]
    In particular, $\pi_n(M^\gp)\simeq \pi_{n+1}S^2$ is nonzero for all $n\ge 1$ (see \cite{ivanov2015nontriviality}). 
\end{example}
 
\section{Commutative monoids}

Animation produces many useful homotopy-coherent structures, but misses some equally fundamental ones. Notably, (connective) spectra --- the homotopy theorist's \textit{commutative groups} --- do not arise as the animation of any 1-category. Simply put, animation produces $\oo$-categorical models of $1$-categorical algebraic theories, but sometimes the algebraic theory itself calls for a natural $\infty$-categorical refinement. 
With an eye towards understanding homotopy-coherent commutative groups, we begin with homotopy-coherent \textit{commutative monoids}, which are also quite interesting on their own, not least for their intimate relationship with \textit{semiadditivity}.

\subsection{The matural numbers}

As with abelian groups, the algebraic theory of commutative monoids is given by $\Mat(\NN)$. We can thus consider, for any $\infty$-category $\cC$, product preserving functors\footnote{Even though $\Mat(\NN)$ is \textit{self-opposite}, mental hygiene dictates retaining the $\op$.}
\[
    \Mat(\NN)^\op \too \cC.
\]
For a $1$-category $\cC$ this indeed recovers the notion of a commutative monoid object, yet it fails to capture the most natural of higher categorical examples.

The semiring of natural numbers $\NN$ is famously a \textit{decategorification} of the category of finite sets~$\Fin$ with its categorical product $\times$ and coproduct $\sqcup$. The universal properties governing these operations ensure they are coherently associative, commutative and distributive. Thus, the underlying groupoid (or space) of $\Fin$ inherits this structure. 

\newcommand{\isom}{\mathrm{iso}}
\begin{defn}
    We define the \tdef{matural numbers} to be the groupoid of finite sets
    \[
        \mdef{\MM} := 
        \Fin_\isom \simeq 
        \bigsqcup_{n \in \NN} B\Sigma_n, 
    \]
    with the disjoint union and cartesian product\footnote{These are no longer categorical (co)products.} symmetric monoidal structures. When convenient, we think of $\MM$ in terms of its skeletal version with a single object $\underline{n}=\{1,\dots,n\}$ for each $n\in \NN$. 
\end{defn}

Among other things, $\MM$ should be a commutative monoid in $\Spc$ with respect to $\sqcup$. This structure however cannot be encoded by a product preserving functor $\Mat(\NN)^\op \to \Spc$ as it would provide a symmetric monoidal functor (where the source is viewed as a discrete groupoid)
\[
    (\underline{-}) \colon \NN \oto{\qquad} \MM
\]
taking an integer $n \in \NN$ to the $n$-element set $\underline{n}  \in \MM$. This is impossible as the commutativity isomorphism in $\MM$,
\[
    \sigma_{n,m} \colon 
    \underline{n+m} = 
    \underline{n} + \underline{m} \iso
    \underline{m} + \underline{n} = 
    \underline{m+n}
    \qin \Sigma_{n+m}
\]
is \textit{not} the identity, but rather the permutation:
\[\begin{tikzcd}[column sep= 1em]
	\textcolor{rgb,255:red,43;green,112;blue,16}{{\ \quad 1\quad \ }} & \textcolor{rgb,255:red,43;green,112;blue,16}{{\ \quad 2\quad\ }} & \textcolor{rgb,255:red,43;green,112;blue,16}{\cdots} & \textcolor{rgb,255:red,43;green,112;blue,16}{{\ \quad n\quad \ }} & \textcolor{rgb,255:red,40;green,69;blue,164}{{n+1}} & \textcolor{rgb,255:red,40;green,69;blue,164}{{n+2}} & \textcolor{rgb,255:red,40;green,69;blue,164}{\cdots} & \textcolor{rgb,255:red,40;green,69;blue,164}{{n+m}} \\
	\\
	\\
	\textcolor{rgb,255:red,40;green,69;blue,164}{{\ \quad 1\quad \ }} & \textcolor{rgb,255:red,40;green,69;blue,164}{{\ \quad 2\quad\ }} & \textcolor{rgb,255:red,40;green,69;blue,164}{\cdots} & \textcolor{rgb,255:red,40;green,69;blue,164}{{\quad m\quad}} & \textcolor{rgb,255:red,43;green,112;blue,16}{{1+m}} & \textcolor{rgb,255:red,43;green,112;blue,16}{{2+m}} & \textcolor{rgb,255:red,60;green,156;blue,22}{\cdots} & \textcolor{rgb,255:red,43;green,112;blue,16}{{n+m}}
	\arrow[color={rgb,255:red,43;green,112;blue,16}, from=1-1, to=4-5]
	\arrow[color={rgb,255:red,43;green,112;blue,16}, from=1-2, to=4-6]
	\arrow[color={rgb,255:red,43;green,112;blue,16}, from=1-4, to=4-8]
	\arrow[color={rgb,255:red,40;green,69;blue,164}, from=1-5, to=4-1]
	\arrow[color={rgb,255:red,40;green,69;blue,164}, from=1-6, to=4-2]
	\arrow[color={rgb,255:red,40;green,69;blue,164}, from=1-8, to=4-4]
\end{tikzcd}\]

The failure is instructive: the true semiring was never $\NN$ itself, but rather its untruncated origin~$\MM$. The way forward is therefore also clear. We simply replace $\NN$ by $\MM$ throughout, defining the algebraic theory $\Mat(\MM)$ and considering \textit{its} models instead. Although $\Mat(\MM)$ is not a 1-category, it is still a fairly tangible object.  

\begin{defn}
    There is a $(2,1)$-category $\mdef{\Mat(\mathbb M)}$ whose objects are natural numbers, $1$-morphisms $m\to n$ are $n\times m$ matrices of finite sets $X=(X_{ij})$, composed by matrix multiplication (with $\sqcup$ and $\times$ in place of $+$ and $\cdot$), and $2$-morphisms are entrywise bijections $X_{ij}\iso Y_{ij}$.
\end{defn}

It is not difficult to fill in the details of this definition by supplying the associators for composition, induced by the associativity and
distributivity of $\sqcup$ and $\times$, and verify the Stasheff pentagon identity. There is, however,
a neater presentation of this $(2,1)$-category, more natural in the nonskeletal version.

For $A,B\in \Fin$, a $B\times A$ matrix of finite sets $(X_{ba})$ can be encoded by a single map 
\[
    X := \bigsqcup X_{ba} \too
    B\times A, 
\]
or equivalently, by a \tdef{span}
\[
    B \longleftarrow X \too A \qin \Fin.
\]
This curious encoding now performs the following trick: \textit{matrix product} is simply \textit{span composition}. Indeed, if $C\from Y \to B$ encodes another such matrix, then their span composition $C\from Y\cdot X \to A$ is defined by the following pullback:
\[\begin{tikzcd}[sep = 1.6em]
	&& {Y\cdot X} && \\
	& \ Y && X\ \\
	C && B && A
	\arrow[from=1-3, to=2-2]
	\arrow[from=1-3, to=2-4]
	\arrow[curve={height=18pt}, from=1-3, to=3-1]
	\arrow["\lrcorner"{anchor=center, pos=0.125, rotate=-45}, draw=none, from=1-3, to=3-3]
	\arrow[curve={height=-18pt}, from=1-3, to=3-5]
	\arrow[from=2-2, to=3-1]
	\arrow[from=2-2, to=3-3]
	\arrow[from=2-4, to=3-3]
	\arrow[from=2-4, to=3-5]
\end{tikzcd}\]
A pullback is also known as a \textit{fiber product}, since it induces \textit{products} on the \textit{fibers}. So the fiber of $Y\cdot X$ over the triple $(c,b,a)$ in $C\times B\times A$ is $Y_{cb}\times X_{ba}$ and forgetting the projection to $B$ is precisely taking the union over all $b \in B$, yielding
\[
    (Y\cdot X)_{ca} \simeq \bigsqcup_{b\in B} Y_{cb} \times X_{ba} \qin \MM.
\]

That is, $\Mat(\MM)$ can be identified, up to equivalence, with a special case of the following general construction of span categories:

\begin{defn}
    Let $\cC$ be an $\infty$-category with \textit{pullbacks}. There is an $\infty$-category $\mdef{\Span}(\cC)$ with the same objects as $\cC$, morphisms from $A$ to $B$ given by spans $A \from X \to B$ in $\cC$ and composition given by the above ``span composition'' formula.\footnote{A thorough, rigorous, and pleasant treatment of spans in $\oo$-categories can be found in \cite{haugseng2023two}.}
\end{defn}

\begin{rem}
    For an $\oo$-category $\cC$, the spans from $A$ to $B$ form a \textit{space}. Even when $\cC$ is an ordinary category, such as $\Fin$, this space is not discrete but rather a \textit{groupoid} whose (iso)morphisms are given by a commutative diagram of the form
    \[\begin{tikzcd}[row sep = 1.5em, column sep = 2em]
    	& X &&& \\
    	A && B && \cC \\
    	& Y
    	\arrow[from=1-2, to=2-1]
    	\arrow[from=1-2, to=2-3]
    	\arrow["\wr"', from=1-2, to=3-2]
    	\arrow["{\in }"{description}, draw=none, from=2-3, to=2-5]
    	\arrow[from=3-2, to=2-1]
    	\arrow[from=3-2, to=2-3]
    \end{tikzcd}\]
    This is essential. The pullback used in the composition of spans is only \textit{defined} up to an isomorphism, so associativity of composition requires nonidentity 2-isomorphisms to witness it. 
\end{rem}

Recapitulating, we have an equivalence of $(2,1)$-categories\footnote{Again, the `op' is optional, but it makes things nicer.} 
\[
    \Mat(\MM)^\op \simeq \Span(\Fin),
\]
which is a perfectly legitimate algebraic theory in the $\oo$-categorical world. Namely, it is an $\infty$-category with finite products --- and, by the evident self-duality, finite coproducts as well. These structures can be seen explicitly from the natural wide-subcategory embeddings
\[
    \Fin \intoo \Span(\Fin) \longhookleftarrow \Fin^\op
\]
by restricting to spans whose backward or forward leg is an \textit{isomorphism}. 
In $\Fin$, every finite set is the coproduct of its points. This remains true in $\Span(\Fin)$. Namely, for every $A\in \Fin$, the maps
\[
    \pt \xleftarrow{\ =\ } \pt \oto{\ a\ } A \qin \Span(\Fin)
\]
exhibit $A$ as the colimit of the constant $A$-shaped diagram on $\pt$. By self-duality, the \textit{transposes} of the above maps, known (again) as the \tdef{Segal maps}, 
\[
    A \ofrom{\ a\ } \pt \xrightarrow{\ =\ } \pt \qin \Span(\Fin)
\]
exhibit $A$ as the \textit{product} of its points.

\begin{defn}
    Let $\cC$ be an $\oo$-category with finite products.
    The $\oo$-category of \tdef{commutative monoids} in $\cC$ is given by
    \[
        \mdef{\CMon}(\cC) := \Fun^\times(\Span(\Fin),\cC).
    \]    
\end{defn}

Forgetting the commutative monoid structure is given by evaluation at $\pt \in \Span(\Fin)$, which gives the \tdef{underlying object} functor and its left adjoint, the \tdef{free commutative monoid} functor (which exists under standard niceness assumptions on $\cC$)
\[
    \mdef{U_{\CMon}}\colon 
    \CMon(\cC) \adj \cC \noloc
    \mdef{F_{\CMon}}.
\]
The precise nature of $F_{\CMon}$ depends on the $\oo$-category $\cC$, but there is a familiar formula when the cartesian product distributes over colimits. 

\begin{prop}\label{F_CMon_Cart}
    When $\cC$ is \textit{cartesian closed}, $F_{\CMon}$ is given by\footnote{In claims involving (co)limits we always tacitly assume they \textit{exist}.}
    \[
        F_{\CMon}(X) \simeq \bigsqcup_{n\in \NN} X^n /\!\!/\Sigma_n \qin \cC.
    \]
\end{prop}
Informally, the free commutative monoid on $X$ in this case consists of multisets of elements in $X$ with disjoint union.
\begin{proof}[Proof idea]
     Let $M \in \CMon(\cC)$ with a map $X \to M$ in $\cC$. For every $n\in \NN$ we get a $\Sigma_n$-equivariant map $X^n \to M^n$. Composing with the $n$-fold multiplication of $M$ we obtain a $\Sigma_n$-\textit{invariant} map, which hence factors through the homotopy orbits:
    \[\begin{tikzcd}
    	{X^n} & {M^n} \\
    	{X^n/\!\!/\Sigma_n} & M.
    	\arrow[from=1-1, to=1-2]
    	\arrow[from=1-1, to=2-1]
    	\arrow["{\text{mult}}", from=1-2, to=2-2]
    	\arrow[dashed, from=2-1, to=2-2]
    \end{tikzcd}\]
    These assemble into a single map
    \[
        \bigsqcup_{n\in \NN} X^n/\!\!/\Sigma_n \too M
        \qin \cC.
    \]
    The hypothesis on $\cC$ comes into play when one wants to promote this to a map of \textit{commutative monoids}. On the pointwise level, the commutative multiplication on the source is given by the canonical maps
    \[
        X^n/\!\!/\Sigma_n \times X^m/\!\!/\Sigma_m \iso
        (X^n \times X^m)/\!\!/(\Sigma_n \times \Sigma_m) \too 
        X^{n+m}/\!\!/\Sigma_{n+m}.
    \]
    The isomorphism above, and the assembly of these maps into a multiplication on the coproduct, rely on binary products distributing over colimits; this follows from cartesian closedness.
\end{proof}

The above formula is compatible with $F_{\CMon}(\pt)$ being $\MM$ in $\CMon(\Spc)$ and $\NN$ in $\CMon(\Set)$. In the noncommutative case, the \tdef{free monoid} functor is given by 
\[
    \mdef{F_{\Mon}}(X) \simeq \bigsqcup_{n\in \NN} X^n \qin \cC.
\]
In particular, the finite free monoids in sets and spaces are the same. We find that $\NN$ is the free homotopy-coherent \textit{monoid} on one element, which happens to be \textit{commutative} (and being discrete, this is a property), yet it is not the free homotopy-coherent \textit{commutative monoid} on one element (which is $\MM$).

\begin{rem}
    The form of $F_{\CMon}$ in general depends a lot on the category $\cC$. By the Eckmann-Hilton argument, for example, every abelian group admits a unique commutative monoid structure so $\CMon(\Ab) \simeq \Ab$ and $F_{\CMon}$ is essentially the \textit{identity}. This boils down to the fact that finite \textit{products} in $\Ab$ are also \textit{coproducts} (see \Cref{CMon_Sadd_Iff} for a fuller discussion). 
\end{rem}

\begin{rem}
    When $\cC$ is a $1$-category, every functor 
    \[
        \Span(\Fin) \simeq \Mat(\MM)^\op \too \cC
    \]
    factors through the \textit{homotopy category} of the source which is $\Mat(\NN)^\op$ and product-preservation is unaffected. Consequently, homotopy-coherent commutative monoids in $\cC$ are just ordinary commutative monoids and the notation $\CMon(\cC)$ is unambiguous. In particular, the ordinary commutative monoids are precisely the homotopy-coherent ones that happen to be homotopically discrete. 
\end{rem}

\subsection{Symmetric monoidal categories}
In ordinary category theory, there is a fundamental difference between a commutative monoid structure on a \textit{set} and a symmetric monoidal structure on a \textit{category}. In the former, unitality, associativity and commutativity identities are imposed on the nose, whereas in the latter they are replaced by isomorphisms, which form additional data and are themselves required to satisfy the triangle, pentagon, and hexagon identities. From our perspective, this reflects the fact that $\Cat$ is not an ordinary category with discrete mapping sets, but a $(2,1)$-category with $1$-truncated mapping groupoids of functors and natural isomorphisms. Since only natural \textit{isomorphisms} are used in the definition, the distinction is not so much a matter of \textit{categorical level} as of \textit{homotopy coherence}. The classical coherence theorem for symmetric monoidal categories can be reinterpreted as the following comparison result:

\begin{prop}
    The $(2,1)$-category of symmetric monoidal categories, symmetric monoidal functors and symmetric monoidal natural isomorphisms is equivalent to $\CMon(\Cat)$.
\end{prop}

Unlike with $\Cat$ and $\Set$, at infinity there is no longer a difference in homotopy truncatedness between $\Cat_\oo$ and $\Spc$. Consequently, the definition of a symmetric monoidal \textit{$\oo$-category} is the same as that of a commutative monoid object in any ambient $\oo$-category. 

\begin{defn}
    The $\oo$-category of \tdef{symmetric monoidal $\oo$-categories} is given by
    \[
        \Cat^{\mdef{\otimes}}_\oo := \CMon(\Cat_\oo).
    \]
\end{defn}

While the noninvertible natural transformations do not participate in the \textit{definition} of a symmetric monoidal structure, they do interact with it in a nontrivial way, making the \textit{theory} of symmetric monoidal $\oo$-categories much richer than that of homotopy-coherent commutative monoids. Most importantly, by the microcosm principle, every symmetric monoidal $\oo$-category carries its own internal notion of a \textit{commutative algebra}. This, however, requires one to contemplate (in one way or another) the $(\oo,2)$-categorical structure of $\Cat_\oo^\otimes$.

The systematic study of this theory will be the subject of the next chapter. For now, we shall treat it mostly as a source of examples. For instance, when an $\oo$-category $\cC$ 
admits finite \textit{products} or \textit{coproducts}, these define a symmetric monoidal structure, called \tdef{cartesian} or \tdef{cocartesian}, respectively. The coherence homotopies for associativity, commutativity and unitality are provided canonically by the corresponding universal properties. This illustrates a general principle: categorification allows for certain canonical \textit{structure} to emerge from \textit{properties}.

Since the functor taking an $\oo$-category to its underlying space preserves limits, we have an induced functor on commutative monoid objects,
\[
    \Cat^\otimes_\oo \simeq 
    \CMon(\Cat_\oo) \too 
    \CMon(\Spc).
\]
Thus, we obtain a very cheap way to produce homotopy-coherent commutative monoids. Indeed, that's how we got $\MM$. 

\begin{example}\label{Ex_Proj_CMon}
    For any ring $R$, one can form a commutative monoid $\mdef{M}\!{}_R$ in (1-truncated) spaces by taking the underlying groupoid of finite free $R$-modules under $\oplus$. 
    For $R=\CC$, we use the special notation
    \[
        M^{\mdef{\delta}}_\CC \simeq \bigsqcup_{n \in \NN} B\GL_n^{\mdef{\delta}}(\CC) \qin \CMon(\Spc),
    \]
    where $\delta$ emphasizes that we treat $\GL_n(\CC)$ as a \textit{discrete} group. We reserve $M_\CC$ for a different commutative monoid in (nontruncated) spaces where one uses instead the \textit{natural topology} on $\GL_n(\CC)$, or equivalently (and somewhat more standardly), 
    \[
        \mdef{M_\CC} := \bigsqcup_{n \in \NN} BU(n) \qin \CMon(\Spc).
    \]
    It does not seem to be inherited from the (co)product of any ``$\infty$-category of $\CC$-vector spaces''.
\end{example}

\subsection{Segal $\Gamma$-spaces}

Having made a point of $\Span(\Fin)$ not being a 1-category, there actually \textit{is} a way to encode commutative monoid objects in an $\oo$-category using only diagrams of 1-categorical shape, which goes back to Segal's seminal work on \textit{$\Gamma$-spaces} \cite{segal1974categories}. 

The core idea has nothing particularly ``homotopy-coherent'' about it. It is essentially a commutative analogue of how Segal simplicial spaces are used to model (noncommutative) monoids (see \Cref{Segal_Mon}). To encode a commutative, associative and unital addition one doesn't need the full category $\Mat(\NN)$. One only needs the strikingly sparse wide subcategory consisting of $0$-$1$ matrices with at most one $1$ in each column. For instance, threefold associativity is encoded by the identity:
\[
    \begin{bmatrix}
        1 & 1
    \end{bmatrix}
    \begin{bmatrix}
        1 & 1 & 0 \\
        0 & 0 & 1
    \end{bmatrix} =
    \begin{bmatrix}
        1 & 1 & 1
    \end{bmatrix} =
    \begin{bmatrix}
        1 & 1
    \end{bmatrix}
    \begin{bmatrix}
        0 & 1 & 1 \\
        1 & 0 & 0
    \end{bmatrix}.
\]
We can just as well consider the wide subcategory of $\Mat(\MM)$ consisting of $\underline{0}$-$\underline{1}$ matrices with at most one $\underline{1}$ in each column. Since $\underline{0}=\es$ and $\underline{1}=\{1\}$ have no nonidentity automorphisms, this wide subcategory is an honest-to-God \textit{1-category}. 

Translating to $\Span(\Fin)$, this wide subcategory corresponds to restricting to those spans whose backward leg is an \textit{injection}:
\[
    A \longhookleftarrow X \too B. 
\]
One can think of such a span as a \textit{partially defined} map $A \rightharpoonup B$ with $X$ the domain of definition. Notably, this is an ordinary category, which is equivalent to the more familiar category $\Fin_*$ of pointed finite sets. Indeed, one can treat the basepoint  as a rubbish bin for all the elements the function is undefined on. Crucially, $\Fin_*$ contains the Segal maps $A \ofrom{a} \pt = \pt$ for all $a \in A$, which allows it to control the product preserving functors from $\Span(\Fin)$. 

\begin{prop}
    Let $\cC$ be an $\oo$-category with finite products. Restriction along the inclusion $\Fin_* \to \Span(\Fin)$ induces an equivalence of categories:
    \[
        \Gray{\CMon(\cC) :=}
        \Fun^\times(\Span(\Fin),\cC) \iso 
        \Fun^{\mdef{\seg}}(\Fin_*,\cC) 
        \Gray{\sseq \Fun(\Fin_*,\cC)},
    \]
    where `seg' stands for the full subcategory of functors $M$ satisfying the \tdef{Segal condition}: for every finite set $A$, the Segal maps jointly induce an equivalence\footnote{A modern treatment of (a vast generalization of) this story can be found in \cite[\S 5.2]{harpaz2020ambidexterity}. }
    \[
         M(A) \iso M(\pt)^A \qin \cC.
    \]
\end{prop}

Modeling commutative monoids with Segal functors has the questionable advantage of $\Fin_*$ being a $1$-category rather than a $(2,1)$-category like $\Span(\Fin)$. From the classical perspective this allowed one to model homotopy-coherent commutative monoids by 1-categorical functors $\Fin_* \to {\rm Top}$. The flip side is that it somewhat obfuscates the ``algebraic theory'' picture. 

\begin{rem}\label{Gamma_Op}
    The op-conventions are liable to cause confusion about the role played here by the category $\Fin_*$. It encodes not so much partial functions as \textit{partial sums}.
    This may also explain why Segal defined the category $\Gamma \simeq \Fin_*^\op$ directly and formulated the theory in terms of functors $\Gamma^\op \to \cC$.
\end{rem}

\begin{rem}
    The two Segal models:
    \[
        \Mon(\cC) \simeq \Fun^{\seg}(\Delta^\op, \cC) \quad,\quad
        \CMon(\cC) \simeq \Fun^{\seg}(\Fin_*, \cC)
    \]
    belong to a bigger story about \textit{operator categories} developed in $\cite{barwick2018operator}$. Pushing this Segal philosophy beyond algebraic theories leads to the framework of \textit{algebraic patterns} of \cite{chu2021homotopy}, which accommodates more general homotopy-coherent structures, notably \textit{category objects}. Iterating and varying this construction yields models for $(\oo,n)$-categories, their flagged and $n$-uple variants, and related structures, all internal to $\oo$-category theory.
\end{rem}

Ordinary commutative monoids are monoids satisfying an additional property. In the homotopy-coherent world, by contrast, it is not even clear what the underlying monoid of a commutative monoid should mean. To define one, we need a finite coproduct preserving functor along which to restrict:
\[
    \Mon^\ff \too \Mat(\MM) \simeq \Span(\Fin)^\op. 
\]
Intuitively, this functor should send a homomorphism of free monoids $F(m) \to F(n)$ specified by words $w_1,\dots,w_m$ in generators $x_1,\dots,x_n$ to its \textit{abelianization} $\MM^m \to \MM^n$: the $n\times m$ matrix whose $(i,j)$-entry is the \textit{set of occurrences} of $x_i$ in $w_j$. Making this precise is somewhat tedious, since the target is a $(2,1)$-category. We may instead trade some conceptual transparency for technical ease by passing to the equivalent Segal models.

\begin{defn}\label{CGrp_Prestable}
    We define a functor
    \[
        \mdef{\Cut} \colon \Delta^\op \too \Fin_*
    \]
    by 
    \[
        \Cut([n]) = \hom_{\Delta}([n],[1])\ \big/\ \mathrm{const}\sim * \qin \Fin_*,
    \]
    where the quotient collapses the constant maps (on either 0 or 1) to a single element, which is taken to be the basepoint.\footnote{Most concisely, this is the standard simplicial circle $\Delta^1/\partial\Delta^1$ viewed as a pointed levelwise finite simplicial set.}  
\end{defn}

The encoding specifics make this look worse than it is. Behind them lies, by no coincidence, essentially the same combinatorics as for the functor
$\Gap\colon \Delta \to \Mon^{\ff}$ (see \Cref{Gap}).

\begin{example}
    We write $\angles{n} = \{*,1,2,\dots,n\}$ for the objects of $\Fin_*$. The following illustrates the value of $\Cut$, viewed as a contravariant functor, on a particular map $\alpha \colon [3]\to [5]$ in $\Delta$  (the basepoint $*$ appears twice for purely typographical reasons). 
    \[\begin{tikzcd}[column sep = tiny]
    	\textcolor{rgb,255:red,36;green,143;blue,36}{{[3]}} &&& \textcolor{rgb,255:red,33;green,72;blue,171}{{*}} & \textcolor{rgb,255:red,36;green,143;blue,36}{0} & \textcolor{rgb,255:red,33;green,72;blue,171}{1} & \textcolor{rgb,255:red,36;green,143;blue,36}{1} & \textcolor{rgb,255:red,33;green,72;blue,171}{2} & \textcolor{rgb,255:red,36;green,143;blue,36}{2} & \textcolor{rgb,255:red,33;green,72;blue,171}{3} & \textcolor{rgb,255:red,36;green,143;blue,36}{3} & \textcolor{rgb,255:red,33;green,72;blue,171}{{*}} &&& \textcolor{rgb,255:red,33;green,72;blue,171}{{\angles{3}}} \\
    	\\
    	\textcolor{rgb,255:red,36;green,143;blue,36}{{[5]}} & \textcolor{rgb,255:red,33;green,72;blue,171}{{*}} & \textcolor{rgb,255:red,36;green,143;blue,36}{0} & \textcolor{rgb,255:red,33;green,72;blue,171}{1} & \textcolor{rgb,255:red,36;green,143;blue,36}{1} & \textcolor{rgb,255:red,33;green,72;blue,171}{2} & \textcolor{rgb,255:red,36;green,143;blue,36}{2} & \textcolor{rgb,255:red,33;green,72;blue,171}{3} & \textcolor{rgb,255:red,36;green,143;blue,36}{3} & \textcolor{rgb,255:red,33;green,72;blue,171}{4} & \textcolor{rgb,255:red,36;green,143;blue,36}{4} & \textcolor{rgb,255:red,33;green,72;blue,171}{5} & \textcolor{rgb,255:red,36;green,143;blue,36}{5} & \textcolor{rgb,255:red,33;green,72;blue,171}{{*}} & \textcolor{rgb,255:red,33;green,72;blue,171}{{\angles{5}}}
    	\arrow["\alpha"', color={rgb,255:red,36;green,143;blue,36}, from=1-1, to=3-1]
    	\arrow[draw={rgb,255:red,36;green,143;blue,36}, from=1-5, to=3-5]
    	\arrow[draw={rgb,255:red,36;green,143;blue,36}, from=1-7, to=3-7]
    	\arrow[draw={rgb,255:red,36;green,143;blue,36}, from=1-9, to=3-7]
    	\arrow[draw={rgb,255:red,36;green,143;blue,36}, from=1-11, to=3-13]
    	\arrow[draw={rgb,255:red,33;green,72;blue,171}, dashed, from=3-2, to=1-4]
    	\arrow[draw={rgb,255:red,33;green,72;blue,171}, dashed, from=3-4, to=1-4]
    	\arrow[draw={rgb,255:red,33;green,72;blue,171}, dashed, from=3-6, to=1-6]
    	\arrow[draw={rgb,255:red,33;green,72;blue,171}, dashed, from=3-8, to=1-10]
    	\arrow[draw={rgb,255:red,33;green,72;blue,171}, dashed, from=3-10, to=1-10]
    	\arrow[draw={rgb,255:red,33;green,72;blue,171}, dashed, from=3-12, to=1-10]
    	\arrow[draw={rgb,255:red,33;green,72;blue,171}, dashed, from=3-14, to=1-12]
    	\arrow["{\Cut(\alpha)}"', color={rgb,255:red,33;green,72;blue,171}, from=3-15, to=1-15]
    \end{tikzcd}\]
    In the associated Segal models, this replaces the sequence of noncommutative words by its abelianization:
    \[
        (x_2,1,x_3x_4x_5) \longmapsto (x_2,0,x_3+x_4+x_5). 
    \]
\end{example}

The restriction along $\Cut$ yields a functor
\[
    \CMon(\cC) = 
    \Fun^{\seg}(\Fin_* , \cC) \too 
    \Fun^{\seg}(\Delta^{\op} , \cC) =
    \Mon(\cC).
\]
Unlike with ordinary monoids, this functor is not fully faithful. In the homotopy-coherent world, commutativity is not a \textit{property} but \textit{structure}. 

\subsection{Limits and colimits}

Understanding a category is, in particular, understanding how limits and colimits are formed in it. The forgetful functor 
\[
    U_{\CMon} \colon \CMon(\cC) \too \cC
\]
preserves all limits and, being conservative, also \textit{creates} them. Intuitively, given a diagram of commutative monoids in $\cC$, the limit of the underlying objects in $\cC$ admits a unique compatible commutative monoid structure making it the limit of the diagram in $\CMon(\cC)$. 

\begin{example}\label{Vect_X}
    Recall from \Cref{Ex_Proj_CMon} the  homotopy-coherent commutative monoid $M_\CC$ of $\CC$-vector spaces. For a space $X$, the constant $X$-shaped limit
    \[
        \mdef{\Vect}_\CC(X) := 
        \invlim_X M_\CC \simeq 
        \Map\big(X, \bigsqcup\nolimits_{n\in \NN} BU(n)\big)
    \]
    is the space of \textit{$\CC$-vector bundles} on $X$ equipped with pointwise $\oplus$.
\end{example}

Colimits in $\CMon(\cC)$ are more interesting. Even in the case $\cC = \Spc$, the functor $U_{\CMon}$ preserves (and creates) \textit{filtered} and even \textit{sifted} colimits,\footnote{This holds for many $\oo$-categories $\cC$, but not all. Counterexamples include $\mathrm{CRing}^\op$ and $\rm Top$. The latter can be fixed by working with \textit{compactly generated weakly Hausdorff} spaces.} but not all colimits. The fundamental example is the \textit{zero} commutative monoid of $\cC$ which is both \textit{terminal} and \textit{initial}. Its underlying object is terminal in $\cC$, but usually \textit{not} initial. Furthermore, 

\begin{prop}\label{CMon_Sadd}
    Finite products in $\CMon(\cC)$, which are computed as products in $\cC$, are also \textit{coproducts} in $\CMon(\cC)$, hence \textit{biproducts}, which are therefore denoted by $\mdef{\oplus}$.
\end{prop}

\begin{proof}[Proof sketch]
    The crucial step is to show that $\Span(\Fin)$ itself has biproducts. As this is a fairly explicit $(2,1)$-category, everything can be checked by hand. In fact, we already noted that every finite set $A\in \Span(\Fin)$ is both the \textit{product} and \textit{coproduct} of its elements and this is just an elaboration of that.

    Now, it is completely formal to show that the $\oo$-category of product preserving functors from \textit{any} $\oo$-category with biproducts, such as $\Span(\Fin)$, to \textit{any} $\oo$-category $\cC$ with finite products has itself biproducts. Indeed, given product preserving functors
    \[
        F,G,H\colon 
        \Span(\Fin) \too 
        \cC,
    \]
    and a pair of natural transformations 
    \[
        \alpha \colon F \too H \quad\text{and}\quad \beta\colon G \too H, 
    \]
    the dashed arrow in the following diagram:
    \[\begin{tikzcd}
    	{F(X)\times G(X)} && {H(X)} \\
    	{H(X)\times H(X)} & {H(X\times X)} & {H(X\sqcup X)}
    	\arrow[dashed, from=1-1, to=1-3]
    	\arrow["{\alpha\times \beta}"', from=1-1, to=2-1]
    	\arrow[from=2-1, to=2-2]
    	\arrow["\sim", shift right, draw=none, from=2-1, to=2-2]
    	\arrow[from=2-2, to=2-3]
    	\arrow["\sim", shift right, draw=none, from=2-2, to=2-3]
    	\arrow["\nabla"', from=2-3, to=1-3]
    \end{tikzcd}\]
    is the essentially unique extension of the given data to a natural transformation $F\times G\to H$, exhibiting the source as the \textit{coproduct} of $F$ and $G$.\footnote{A rigorous argument can be found in \cite[Corollary 2.4]{gepner2016universality}.} 
\end{proof}

By the general yoga of \textit{simplicial replacement}, a colimit of any shape can be written as a geometric realization of coproducts, providing a fairly concrete recipe for computing it.\footnote{This trick is particularly popular in classical formulas for ``homotopy colimits'' (see \cite{dugger2008primer}).} 
The case of a \textit{pushout} is particularly important and illuminating. The colimit of a span diagram 
\[
    N \froom M \too L \qin \CMon(\cC)    
\]
is isomorphic to the colimit of the simplicial diagram known as the  \tdef{two-sided bar construction}:
\[\begin{tikzcd}
	{\mdef{\BB}(N,M,L)} & {\dots } & {N\oplus M\oplus M\oplus L} & {N\oplus M\oplus L} & {N\oplus L}
	\arrow["{:=}"{marking, allow upside down}, draw=none, from=1-1, to=1-2]
	\arrow[shift left=3, from=1-2, to=1-3]
	\arrow[shift right=3, from=1-2, to=1-3]
	\arrow[shift left, from=1-2, to=1-3]
	\arrow[shift right, from=1-2, to=1-3]
	\arrow[between={0.2}{0.8}, from=1-3, to=1-2]
	\arrow[shift left=2, between={0.2}{0.8}, from=1-3, to=1-2]
	\arrow[shift right=2, between={0.2}{0.8}, from=1-3, to=1-2]
	\arrow[from=1-3, to=1-4]
	\arrow[shift left=2, from=1-3, to=1-4]
	\arrow[shift right=2, from=1-3, to=1-4]
	\arrow[shift left, between={0.2}{0.8}, from=1-4, to=1-3]
	\arrow[shift right, between={0.2}{0.8}, from=1-4, to=1-3]
	\arrow[shift right, from=1-4, to=1-5]
	\arrow[shift left, from=1-4, to=1-5]
	\arrow[between={0.2}{0.8}, from=1-5, to=1-4]
\end{tikzcd}\]

We obtain the following formula for pushouts:

\begin{prop}
    The pushout of a span diagram $N \from M \to L$ in $\CMon(\cC)$ is given by the colimit of the two-sided bar construction
    \[
        N\oplus_M L \simeq 
        |\BB(N,M,L)| \qin \CMon(\cC). 
    \]
\end{prop}

The special case of $N=L=0$,
\[\begin{tikzcd}
	{\BB(0,M,0)} & {\cdots M\oplus M\oplus M} & {M\oplus M} & M & 0
	\arrow["{:=}"{marking, allow upside down}, draw=none, from=1-1, to=1-2]
	\arrow[shift left=3, from=1-2, to=1-3]
	\arrow[shift right=3, from=1-2, to=1-3]
	\arrow[shift left, from=1-2, to=1-3]
	\arrow[shift right, from=1-2, to=1-3]
	\arrow[between={0.2}{0.8}, from=1-3, to=1-2]
	\arrow[shift left=2, between={0.2}{0.8}, from=1-3, to=1-2]
	\arrow[shift right=2, between={0.2}{0.8}, from=1-3, to=1-2]
	\arrow[from=1-3, to=1-4]
	\arrow[shift left=2, from=1-3, to=1-4]
	\arrow[shift right=2, from=1-3, to=1-4]
	\arrow[shift left, between={0.2}{0.8}, from=1-4, to=1-3]
	\arrow[shift right, between={0.2}{0.8}, from=1-4, to=1-3]
	\arrow[shift right, from=1-4, to=1-5]
	\arrow[shift left, from=1-4, to=1-5]
	\arrow[between={0.2}{0.8}, from=1-5, to=1-4]
\end{tikzcd}\]
is the usual bar construction $\BB M$ from \Cref{Segal_Mon} lifted to~$\CMon(\cC)$, which now gets reinterpreted as the internal \textit{suspension} in $\CMon(\Spc)$.

\begin{cor}
    For $M \in \CMon(\Spc)$ the suspension of $M$ internally to the $\oo$-category of commutative monoids has underlying space $BM = |\BB M|$.
\end{cor}

\begin{rem}
    For $M\in \CMon(\Spc)$, the notation $\Omega M$ is unambiguous as the underlying space functor commutes with taking $\Omega$. However, for $\Sigma M$ one has to specify whether we are considering the suspension operation internal to $\CMon(\Spc)$ or of the underlying space of $M$, which rarely even admits a commutative monoid structure. To avoid confusion, we shall use $BM$ to refer to the suspension operation in $\CMon(\Spc)$, indicating the construction on the level of spaces. 
\end{rem}

\begin{rem}
    Even for \textit{noncommutative} monoids, given a span diagram 
    \[
        N \froom M \too L \qin \Mon(\cC),    
    \]
    one can produce the two-sided bar construction simplicial diagram $\BB(N,M,L)$,
    \[\begin{tikzcd}
    	{\dots } & {N\times M\times M\times L} & {N\times M\times L} & {N\times L} & \cC
    	\arrow[shift left=3, from=1-1, to=1-2]
    	\arrow[shift right=3, from=1-1, to=1-2]
    	\arrow[shift left, from=1-1, to=1-2]
    	\arrow[shift right, from=1-1, to=1-2]
    	\arrow[between={0.2}{0.8}, from=1-2, to=1-1]
    	\arrow[shift left=2, between={0.2}{0.8}, from=1-2, to=1-1]
    	\arrow[shift right=2, between={0.2}{0.8}, from=1-2, to=1-1]
    	\arrow[from=1-2, to=1-3]
    	\arrow[shift left=2, from=1-2, to=1-3]
    	\arrow[shift right=2, from=1-2, to=1-3]
    	\arrow[shift left, between={0.2}{0.8}, from=1-3, to=1-2]
    	\arrow[shift right, between={0.2}{0.8}, from=1-3, to=1-2]
    	\arrow[shift right, from=1-3, to=1-4]
    	\arrow[shift left, from=1-3, to=1-4]
    	\arrow[between={0.2}{0.8}, from=1-4, to=1-3]
    	\arrow["{\in }"{marking, allow upside down}, draw=none, from=1-4, to=1-5]
    \end{tikzcd}\]
    However, in general, it can't be lifted to $\Mon(\cC)$.\footnote{This is because the multiplication map $M\times M\to M$ is not a monoid homomorphism unless $M$ is commutative.} Intuitively, the colimit of this diagram computes the quotient of $N\times L$ balancing the right action of $M$ on $N$ with the left action on $L$. It does not in general admit a monoid structure and, in particular, does not compute the pushout in $\Mon(\cC)$.
\end{rem}

\subsection{Semiadditivity}

The property of $\CMon(\cC)$ proved in \Cref{CMon_Sadd}, that finite products and coproducts coincide, cuts to the very heart of what commutative monoids are about. The reason, put succinctly, is that \Cref{CMon_Sadd} admits a \textit{converse}. 
We begin by describing the situation a bit more carefully. 

\begin{defn}
    A \tdef{zero object} in an $\oo$-category $\cC$ is an object $\mdef{0}$ that is both \textit{initial} and \textit{terminal}. If such an object $0$ exists it is essentially unique and $\cC$ is called \tdef{pointed}.\footnote{Or \textit{$(-1)$-semiadditive}.}     
    In such a case, for every pair of objects $X,Y\in \cC$ there is a canonical \tdef{zero map} defined by the composite of the unique maps
    \[
        \mdef{0_{X,Y}} \colon X \too 0 \too Y \qin \cC.
    \]
\end{defn}

If $\cC$ has a terminal object $\pt$, there is a canonical way to make it pointed by considering \textit{pointed objects} in $\cC$. 

\begin{prop}
    For an $\oo$-category $\cC$ with a terminal object, the forgetful functor
    \[
        U_{\mdef{\pt}} \colon \cC_{\pt/} \too \cC
    \]
    is an equivalence if and only if $\pt$ is also \textit{initial}, and hence \textit{zero}.     
\end{prop}

Thus, $\cC$ is pointed if and only if its objects are uniquely pointed, in which case the mapping spaces of $\cC$ are canonically pointed making it \textit{enriched} in $\Spc_*$.

\begin{rem}
    The structure of a \textit{basepoint} is also governed by an algebraic theory, the category $\Fin_*$ of pointed finite sets. If $\cC$ admits finite products, there is an equivalence
    \[
        \cC_{\pt/} \ \simeq \
        \Fun^\times(\Fin_*^\op,\cC).
    \]
    Note that the category $\Fin_*$ plays a different role here from the one in the Segal model for commutative monoids. For the latter, the fold map $\underline{2} \to \underline{1}$ encodes \textit{addition} $+\colon M\times M\to M$, while for a pointed object it simply encodes the \textit{diagonal} $\Delta \colon M \to M\times M$ (cf.\! \Cref{Gamma_Op}).
\end{rem}

If $\cC$ is pointed and furthermore admits finite (co)products, for every pair of objects $X,Y\in \cC$ there is a canonical map
\[
    \begin{bmatrix}
        1_X& 0_{X,Y} \\
        0_{Y,X} & 1_Y
    \end{bmatrix} \colon 
    X \sqcup Y \too
    X \times Y \qin \cC.
\]

\begin{defn}
    An $\oo$-category $\cC$ is called \tdef{semiadditive}\footnote{Or \textit{$0$-semiadditive}.} if it admits a zero object, finite (co)products, and for all $X,Y\in \cC$ the canonical map above is an isomorphism 
    \[
        X\sqcup Y\iso X\times Y
        \qin \cC,
    \]
    in which case we denote both sides by $X\oplus Y$.
\end{defn}

The fundamental theorem about semiadditivity is that it provides a canonical commutative monoid structure on all the objects of the $\oo$-category.  

\begin{thm}\label{CMon_Sadd_Iff}
    For an $\oo$-category $\cC$ with finite products, the forgetful functor is an \textit{equivalence} 
    \[
        U_{\CMon}\colon \CMon(\cC) \iso \cC,
    \]
    if and only if $\cC$ is \textit{semiadditive} (if and only if $h\cC$ is semiadditive). 
\end{thm}
\begin{proof}[Proof idea]
    Since $\cC \to h\cC$ preserves finite (co)products, it is rather easy to see that $\cC$ is semiadditive if and only if $h\cC$ is (provided $\cC$ was assumed in advance to admit finite products). By \Cref{CMon_Sadd}, the $\oo$-category $\CMon(\cC)$ is semiadditive, so if $U_{\CMon}$ is an equivalence, $\cC$ must be semiadditive too. Conversely, if $\cC$ is semiadditive, we need to show that every object $M\in \cC$ admits a unique (functorial) commutative monoid structure. First of all, being initial, there is a unique map from the zero object $0 \to M$ which is thus the unique unit element for a commutative monoid structure on $M$. A multiplication map
    \[
        m \colon M \oplus M \too M
    \]
    has to restrict to the identity on both copies of $M$ by unitality, but as $M\oplus M$ is also the \textit{coproduct} of $M$ and $M$ it forces the multiplication map to be the \textit{fold} map $\nabla$. More systematically, a span
    \[
        A \froom X \too B \qin \Span(\Fin)
    \]
    has to be mapped to the matrix 
    \[
        |X_{ab}| \colon 
        M^{\oplus A} \too 
        M^{\oplus B}
        \qin \cC,
    \]
    which gives the essentially unique commutative monoid structure on $M$.\footnote{A rigorous argument can be found in \cite[Proposition 2.3]{gepner2016universality} and in \cite[Corollary 5.15]{harpaz2020ambidexterity}.}
\end{proof}

The canonical commutative monoid structure on \textit{objects} induces a canonical commutative monoid structure on \textit{morphisms}. 

\begin{cor}
    Let $\cC$ be a semiadditive $\oo$-category. For every pair of objects $X,Y \in \cC$, the mapping space $\Map(X,Y)$ admits a canonical commutative monoid structure. 
\end{cor}
\begin{proof}
    For every $X\in \cC$, the associated corepresentable functor 
    \[
        \Map(X,-)\colon \cC \too \Spc
    \]
    is product preserving and this induces a functor on commutative monoid objects which fits into the following commutative diagram:
    \[\begin{tikzcd}[sep = 2.7em]
    	{\CMon(\cC)} & {\CMon(\Spc)} \\
    	\cC & \Spc
    	\arrow[from=1-1, to=1-2]
    	\arrow["\wr"', from=1-1, to=2-1]
    	\arrow["U", from=1-2, to=2-2]
    	\arrow[dashed, from=2-1, to=1-2]
    	\arrow["{\Map(X,-)}"', from=2-1, to=2-2]
    \end{tikzcd}\]
    The dashed arrow lifts the mapping space functor into commutative monoids in spaces. 
\end{proof}
In particular, under the assumption of semiadditivity, any pair of morphisms $f,g\colon X\to Y$ can be summed $f+g\colon X\to Y$ via the composition
\[
    X \oto{\ \Delta\ }
    X\oplus X \oto{f\oplus g}
    Y\oplus Y \oto{\ \nabla\ }
    Y.
\]
The more structural statement is that a semiadditive $\oo$-category $\cC$ is canonically \textit{enriched} in $\CMon(\Spc)$. So not only are the mapping spaces commutative monoids, but composition is (coherently) \textit{bilinear}.\footnote{We shall define what this means exactly later on.} 

\begin{rem}
    Semiadditivity is a \textit{property} for an $\oo$-category, which induces a great deal of coherent \textit{structure} on the mapping spaces. This is especially valuable in the $\oo$-categorical world where constructing algebraic structures ``by hand'' is very hard. If endomorphism monoids are the paradigmatic source of homotopy-coherent monoids, then mapping spaces in semiadditive $\infty$-categories play the same role for homotopy-coherent \textit{commutative} monoids. It is similarly (tautologically) \textit{universal}: every such $M$ arises as a mapping space $\Map_{\CMon(\Spc)}(\MM,M)$ in the semiadditive $\oo$-category~$\CMon(\Spc)$.
\end{rem}

With the notion of addition of morphisms one can start doing some elementary \textit{arithmetic}. 

\begin{defn}
    Every homotopy-coherent commutative monoid $M$ admits a natural \textit{multiplication by $m$} map
    \[
        \overset{\mdef{\cdot m}}{\overbrace{\Id+\Id+\cdots +\Id}} :
        M \too M \qin \CMon(\Spc),
    \]
    of which one can take the \textit{cofiber} (or \textit{quotient}) given by the following pushout square in $\CMon(\Spc)$,
    \[\begin{tikzcd}
    	M & M \\
    	0 & {M\mdef{/m}.}
    	\arrow["{\cdot m}", from=1-1, to=1-2]
    	\arrow[from=1-1, to=2-1]
    	\arrow[from=1-2, to=2-2]
    	\arrow[from=2-1, to=2-2]
    	\arrow["\lrcorner"{anchor=center, pos=0.125, rotate=180}, draw=none, from=2-2, to=1-1]
    \end{tikzcd}\]
\end{defn}

For an ordinary commutative monoid $M$, the multiplication by $m$ map on the quotient $M/m$ is \textit{zero}, as we have just killed all the multiples of $m$ in $M$. In the homotopy-coherent world, such identifications are \textit{structure} which has to interact coherently with the monoid operation. As it happens, $\pi_2(\MM/2) = \ZZ/4$ and so multiplication by $2$ is \textit{not} null-homotopic on $\MM/2$.\footnote{This is in fact the same as the \textit{mod $m$ Moore spectrum} $\Sph/m$.} One can show, however, that multiplication by $m^2$ \textit{is} null-homotopic on $M/m$.

\begin{rem}\label{Coherence_Puzzle}
    As a test of \textit{coherence kung fu}, we offer the following puzzle. Consider the diagram:
    \[\begin{tikzcd}
    	\textcolor{rgb,255:red,204;green,51;blue,51}{M} && \textcolor{rgb,255:red,204;green,51;blue,51}{M} && \textcolor{rgb,255:red,204;green,51;blue,51}{{M/m}} \\
    	\textcolor{rgb,255:red,204;green,51;blue,51}{M} && \textcolor{rgb,255:red,204;green,51;blue,51}{M} && \textcolor{rgb,255:red,204;green,51;blue,51}{0} \\
    	\textcolor{rgb,255:red,204;green,51;blue,51}{M} && \textcolor{rgb,255:red,204;green,51;blue,51}{M} && \textcolor{rgb,255:red,204;green,51;blue,51}{{M/m}}
    	\arrow["{\cdot m}", color={rgb,255:red,204;green,51;blue,51}, from=1-1, to=1-3]
    	\arrow["{\cdot m}"', color={rgb,255:red,204;green,51;blue,51}, from=1-1, to=2-1]
    	\arrow[draw={rgb,255:red,204;green,51;blue,51}, from=1-3, to=1-5]
    	\arrow[draw={rgb,255:red,204;green,51;blue,51}, equals, from=1-3, to=2-3]
    	\arrow[draw={rgb,255:red,204;green,51;blue,51}, from=1-5, to=2-5]
    	\arrow["{\cdot m}", color={rgb,255:red,204;green,51;blue,51}, curve={height=-18pt}, from=1-5, to=3-5]
    	\arrow[draw={rgb,255:red,204;green,51;blue,51}, equals, from=2-1, to=2-3]
    	\arrow[draw={rgb,255:red,204;green,51;blue,51}, equals, from=2-1, to=3-1]
    	\arrow[draw={rgb,255:red,204;green,51;blue,51}, from=2-3, to=2-5]
    	\arrow["{\cdot m}", color={rgb,255:red,204;green,51;blue,51}, from=2-3, to=3-3]
    	\arrow[draw={rgb,255:red,204;green,51;blue,51}, from=2-5, to=3-5]
    	\arrow["{\cdot m}"', color={rgb,255:red,204;green,51;blue,51}, from=3-1, to=3-3]
    	\arrow[draw={rgb,255:red,204;green,51;blue,51}, from=3-3, to=3-5]
    \end{tikzcd}\]
    The right column consists of the cofibers of the horizontal maps between the left and middle columns. The outer left rectangle is multiplication by $m$ on all sides, so by \textit{naturality} of multiplication by $m$ its induced map on cofibers should again be multiplication by $m$ on $M/m$, which seems to imply it is always \textit{null}. What is wrong with this argument?
\end{rem}

\section{Commutative groups}

We finally turn to the most fundamental homotopy-coherent algebraic structure: the commutative group. In spaces, these turn out to be precisely \textit{connective spectra}, which give rise to generalized (co)homology theories and help explain the ubiquity of spectra throughout mathematics. Just as abelian groups underlie ordinary algebra, spectra form the bedrock for homotopy-coherent rings, modules, Lie algebras, and much else besides. This equivalence, however, is not tautological. Before explaining it, we first consider homotopy-coherent commutative groups as a natural subclass of homotopy-coherent commutative monoids. 

Commutative groups are, again, simply commutative monoids where every element is \textit{invertible}. And again, there is no subtlety in how the homotopy-coherent commutativity interacts with this notion. There is a \textit{commutative shear map} for every $M \in \CMon(\cC)$, encoded by the matrix
\[
    \begin{bmatrix}
        \underline{1} & \underline{1} \\
        \underline{0} & \underline{1}
    \end{bmatrix} \colon 
    \MM^2 \too \MM^2.
\]

\begin{defn}
    Let $\cC$ be an $\oo$-category with finite products.
    The $\infty$-category of \tdef{commutative groups} in $\cC$ is the full subcategory 
    \[
        \mdef{\CGrp}(\cC)\ \sseq\ \CMon(\cC)
    \]
    spanned by those commutative monoids $M$ for which the shear map is an isomorphism.
\end{defn}

This also leads to the following natural strengthening of the \textit{semiadditivity} property. 

\begin{defn}
    An $\oo$-category $\cC$ is called \tdef{additive} if it is \textit{semiadditive} and for all $X,Y\in \cC$ the commutative monoid $\Map_\cC(X,Y)$ is in $\CGrp(\Spc)$.
\end{defn}

Thus, \textit{additive} $\oo$-categories are canonically enriched in $\CGrp(\Spc)$ just like \textit{semiadditive} ones are enriched in $\CMon(\Spc)$. We also have the following additive analogue of \Cref{CMon_Sadd_Iff}:

\begin{thm}\label{CGrp_add_Iff}
    For an $\oo$-category $\cC$ with finite products, the forgetful functor is an \textit{equivalence} 
    \[
        U_{\mdef{\CGrp}}\colon \CGrp(\cC) \iso \cC,
    \]
    if and only if $\cC$ is \textit{additive} (if and only if $h\cC$ is additive). 
\end{thm}

The $\oo$-category $\CGrp(\Spc)$ \textit{itself} is additive while $\CMon(\Spc)$ is not. The connective and unbounded derived $\oo$-categories $\sD_{\ge0}(R)$ and $\sD(R)$, for every ring $R$, are additive as well. The latter is even \textit{stable}, which is a stronger property not enjoyed by $\CGrp(\Spc)$.  

\subsection{Commutative group completion}

One could equivalently define a commutative group as a commutative monoid whose underlying monoid is a group. This rests on the easily verified fact that the \textit{commutative shear map} is mapped to the \textit{noncommutative shear map} under the forgetful functor $\CMon(\cC) \to \Mon(\cC)$. Namely, we have the following pullback diagram: 
\[\begin{tikzcd}[sep = 2.2em]
	{\CGrp(\cC)} & {\CMon(\cC)} \\
	{\Grp(\cC)} & {\Mon(\cC).}
	\arrow[hook, from=1-1, to=1-2]
	\arrow[from=1-1, to=2-1]
	\arrow["\lrcorner"{anchor=center, pos=0.125}, draw=none, from=1-1, to=2-2]
	\arrow[from=1-2, to=2-2]
	\arrow[hook, from=2-1, to=2-2]
\end{tikzcd}\]

The left adjoint of the bottom arrow is the (noncommutative) \textit{group completion} functor. The left adjoint of the top arrow is likewise the \tdef{commutative group completion}. Combining our observations about the bar construction for $\Grp(\Spc)$ and $\CMon(\Spc)$ we get a neat reinterpretation of it.

\begin{prop}\label{B_CGrp}
    For $M \in \CMon(\Spc)$, the unit of the internal loop--suspension adjunction
    \[
        M \too \Omega BM \qin \CMon(\Spc)
    \]
    exhibits the target as the reflection of $M$ into the full subcategory $\CGrp(\Spc) \sseq \CMon(\Spc)$.
\end{prop}

By abuse of notation, we denote this reflection by $M \mapsto M^{\mdef{\gp}}$. 
This is justified because its underlying (noncommutative) group object is the (noncommutative) group completion of the underlying (noncommutative) monoid object of $M$.\footnote{It is not obvious from abstract formal considerations and amounts to the above pullback diagram of $\oo$-categories being horizontally left \textit{adjointable} (aka satisfying the Beck-Chevalley condition).} 

Specializing \Cref{B_CGrp} to the case of $G \in \CGrp(\Spc)$, we have an \textit{equivalence}
\[
    G \iso \Omega BG \qin \CGrp(\Spc).
\]
Since $\CGrp(\Spc)$ is closed under (co)limits in $\CMon(\Spc)$, the above is also the unit of the loop--suspension adjunction in $\CGrp(\Spc)$, and so $B$ is \textit{fully faithful} with left-inverse-right-adjoint $\Omega$. The essential image of $B$ is also not a big mystery.

\begin{prop}
    The essential image of the internal suspension functor $B$ on $\CGrp(\Spc)$ consists precisely of those commutative groups whose underlying space is \textit{connected}.
\end{prop}
\begin{proof}
    We need to analyze when the \textit{counit} map 
    \[
        B\Omega G \too G \qin \CGrp(\Spc)
    \]
    is an equivalence, which can be checked on underlying spaces. As already noted in the discussion leading to \Cref{Bar_Grp_Cpl}, the colimit $|\BB \Omega G|$ recovers the connected component of the unit element~$1\colon \pt \to G$. 
\end{proof}

The natural extension of this analysis shows that the bar construction on $\CGrp(\Spc)$ behaves like the shift operator on the connective derived category of a ring.

\begin{cor}\label{Bar_Omega_CGrp}
    The internal loop--suspension adjunction 
    \[
        B \colon \CGrp(\Spc) \adj \CGrp(\Spc) \noloc \Omega
    \]
    restricts to an \textit{equivalence} between the full subcategories on $n$-connective and $(n+1)$-connective groups for all $n\ge 0$.\footnote{Here and elsewhere, $n$-connective means $(n-1)$-connected.}
\end{cor}

Group completion of homotopy-coherent commutative monoids provides a rich source of homotopy-coherent commutative groups. In some cases, it can be explicitly computed. In others, not quite. We give examples of both.  

Recall from \Cref{Ex_Proj_CMon} the homotopy-coherent commutative monoid of $\CC$-vector spaces with~$\oplus$,
\[
    M_\CC \simeq \bigsqcup_{n\in \NN} BU(n) \qin \CMon(\Spc).
\]
Its group completion is a well-known classical object. 
\begin{defn}\label{Ex_Top_K_Thy}
    Complex \textit{connective} \tdef{topological $K$-theory} is given by
    \[
        \mdef{\ku} := M_\CC^\gp \qin \CGrp(\Spc).
    \]
\end{defn}

The underlying space of $\ku$ admits a rather explicit description. 

\begin{prop}
    There is an equivalence of spaces 
    \[
        \ku \simeq \ZZ \times BU(\oo) \qquad\text{with}\qquad 
        U(\mdef{\oo}) = \colim U(n).
    \]
\end{prop}
\begin{proof}[Proof idea]
    We begin by observing that there is a canonical self-map of $M_\CC$ given by direct sum with a one-dimensional space $\CC$. On components it is induced by the embedding $U(n)\into U(n+1)$ as an upper $n\times n$ block. The sequential colimit of the iterative application of this operation is denoted
    \[
        M_{\CC, \mdef{\oo}} := 
        \colim_{\NN} \big( 
            M_\CC \oto{\ \oplus \CC\ }
            M_\CC \oto{\ \oplus \CC\ }
            M_\CC \oto{\ \oplus \CC\ }
            \dots
        \big)
    \]
    and is easily seen to satisfy
    \[
        M_{\CC,\oo} \simeq 
        \ZZ \times BU(\oo)  \qin \Spc.
    \]
    After passing to the group completion, the direct sum with $\CC$ becomes an \textit{equivalence}. Thus, the canonical map from $M_\CC$ to its group completion factors as
    \[
        M_\CC \too
        M_{\CC,\oo} \too
        M_\CC^\gp.
    \]
    The most nontrivial part is to show that, in this case, the second map is an \textit{equivalence}. This can be found in \cite{nikolaus2017group}.
\end{proof}

\begin{rem}
    From a modern perspective, referring to a group as a `theory' is rather odd. This terminology reflects the historical origins of $K$-theory, where the group
    \[
        \mdef{K^0}(X) := \pi_0\Map(X,\ku) \qin \Ab
    \]
    occurs as the degree-zero term of a certain generalized \textit{cohomology theory} applied to the space~$X$.\footnote{Not quite the one represented by $\ku$, but by the closely related $\KU$ (see \Cref{Def_KU}).} For a \textit{compact} $X$, this is precisely the ordinary group completion of the commutative monoid $\pi_0\Vect_\CC(X)$ of isomorphism classes of vector bundles on $X$ under $\oplus$ (see \Cref{Vect_X}).\footnote{Compactness is necessary to commute the mapping space with the sequential colimit in the computation of the group completion.}
\end{rem}

The \textit{algebraic} analogue of $K$-theory turns out to be considerably more involved. Recall from \Cref{Ex_Proj_CMon} the homotopy-coherent commutative monoid $M_R$ of finite free $R$-modules. 

\begin{defn}\label{Ex_Alg_K_Thy}
    The \tdef{algebraic $K$-theory} of a ring $R$ is given by\footnote{It is more common to use all finite \textit{projective} modules rather than just the \textit{free} ones, but this changes only~$K_0(R)$.}
    \[
        \mdef{K}(R) := M_R^\gp \qin \CGrp(\Spc).
    \]
    Its $n$th homotopy group is denoted $K_{\mdef{n}}(R)$. 
\end{defn}

One may attempt to compute its underlying space by the same procedure used for $\ku$ and indeed, the map $M_R \to K(R)$ factors through the space
\[
    (M_R)_\oo \simeq \ZZ \times B\GL_\oo(R) \qquad\text{with}\qquad \GL_{\mdef{\oo}}(R) := \colim \GL_n(R).
\]
But the latter fails to identify with $M_R^\gp$. In fact, as its fundamental group $\GL_\oo(R)$ is not even \textit{abelian}, it cannot be the loop space of anything. Hence, a further correction is required: Quillen's \textit{plus construction}.\footnote{For a modern discussion of this notion see \cite{hoyois2019quillen}.} Roughly speaking, it is characterized by abelianizing the fundamental group without altering the \textit{homology} groups. Its effect on \textit{homotopy} groups, however, is quite mysterious.   

\begin{example}
    Since $\ZZ$ is the initial ring, $K(\ZZ)$ is mapped canonically to $K(R)$ for any ring $R$.
    The groups $K_n(\ZZ)$ have not been fully computed; even $K_{12}(\ZZ)$ is still unknown. What \textit{is} known draws on an extraordinary range of mathematics, including the geometry and homology of arithmetic groups, \'{e}tale and motivic cohomology, Borel regulators, Iwasawa theory, and special values of zeta functions, with techniques spanning from stable homotopy theory to harmonic analysis. 
\end{example}

The algebraic $K$-theory of $\CC$ admits a natural comparison with \textit{topological} $K$-theory. The \textit{identity} maps 
\[
    \GL_n^\delta(\CC) \too \GL_n(\CC)
\]
induce a map of homotopy-coherent commutative monoids
\[
    M_\CC^\delta = 
    \bigsqcup B\GL_n^\delta(\CC) \too 
    \bigsqcup B\GL_n(\CC) = 
    M_\CC.
\]
The source and target of this map could hardly look more different. The first is 1-truncated with gigantic (uncountable) noncommutative fundamental groups. The second has simply connected components with nonzero homotopy groups in arbitrarily high degrees. Remarkably, upon \textit{group completion} the map becomes almost an equivalence.  

\begin{thm}[Suslin]
    The canonical map 
    \[
        K(\CC) \too \ku \qin \CGrp(\Spc)
    \]
    induces an equivalence on \textit{profinite} completions.\footnote{Equivalently, we get an equivalence $K(\CC)/m\iso \ku/m$ for every $m \ge 1$.} 
\end{thm}

Even more remarkably, we get the \textit{same result} when we replace $\CC$ with $\cl{\FF}_p$, as long as we complete away from $p$. 
Nor do the wonders end there: under this comparison, the \textit{Frobenius} automorphism of $K(\cl{\FF}_p)$ corresponds to the \textit{Adams} operation $\psi^p$ on $\ku$ and, away from $p$, the homotopy fixed points identify with~$K(\FF_p)$. All this is incredibly explicit yielding the following Fields Medal-worthy computation:

\begin{thm}[Quillen]
    For all primes $p$ and $n\in \NN$,
    \[
        K_n(\FF_p) \simeq \begin{cases}
            \ZZ & n=0 \\
            0 & n= 2m \neq 0 \\
            \ZZ/(p^m-1) & n = 2m-1.
        \end{cases}
    \]
\end{thm}
Virtually every known computation of (higher) algebraic $K$-theory ultimately rests on the one above. 

There is one more homotopy-coherent commutative group of particular importance, arguably the most important of all. It plays the role of $\ZZ$ in ordinary mathematics. 

\begin{defn}
    The free homotopy-coherent commutative group on one element
    \[
        \mdef{\Sph} := F_{\CGrp}(\pt)  \qin \CGrp(\Spc)
    \]
    is called the \tdef{sphere spectrum} (for reasons that will be discussed shortly).
\end{defn}

The adjoint functor theorem guarantees the existence of $F_\CGrp$ but offers no useful formula for it. Thus this definition gives little indication of what the underlying space of $\Sph$ actually looks like. 
There is nevertheless one enlightening description: just as $\ZZ$ is obtained from $\NN$ by group completion, chaining the relevant adjunctions yields the following homotopy-coherent analogue:\footnote{It is tempting to call this the \textit{Barratt--Priddy--Quillen theorem}. This would, however, be somewhat anachronistic. It is more accurate to view the latter as providing the \textit{homological criterion} for recognizing this group completion.} 
\[
    \Sph := F_{\CGrp}(\pt) \simeq 
        F_{\CMon}(\pt)^\gp \simeq 
        \MM^\gp 
        \qin \CGrp(\Spc).
\]
As with $K(R)$, the canonical map $\MM \to \Sph$ factors as
\[
    \MM \too 
    \MM_\oo \simeq \ZZ \times B\Sigma_\oo \too 
    \Sph
    \qquad\text{with}\qquad
    \Sigma_\oo := \colim \Sigma_n.
\]
However, $\Sigma_\oo$ is again noncommutative, so the second map is not an equivalence. The homotopy groups of $\Sph$ have been studied intensively, but no orderly closed-form description of them is known or expected (see \Cref{Stable_Stems}).

Reflecting back on our definition of $\CGrp(\Spc)$, we make the following observation. 
The category $\Ab$ can be defined directly, without passing through $\CMon$, as the category of models in sets of the algebraic theory $\Mat(\ZZ)$.
In a sense, the homotopy-coherent analogue also holds: $\CGrp(\Sph)$ is the category of models in $\Spc$ of the full subcategory $\Mat(\Sph) \sseq \CGrp(\Spc)$ on the finite free objects 
\[
    F_{\CGrp}(\underline{n}) \ \simeq \
    \Sph^{\oplus n} \quad,\quad n \in \NN.
\]
Unlike its ordinary counterpart, however, this does not provide an independent construction, since the mapping spaces
\[
    \Map_{\Mat(\Sph)}(n,m) \simeq 
    \Map_{\CGrp}(\Sph^{\oplus n}, \Sph^{\oplus m})
\]
are difficult to describe a priori without constructing $\CGrp(\Spc)$ first in some other way. 
Yet it does offer an interesting perspective on the connection between $\CGrp(\Spc)$ and the derived $\oo$-category~$\sD_{\ge 0}(\ZZ)$.
Passing to the \textit{connected components} on mapping spaces gives a finite coproduct preserving functor
\[
    \pi_0 \colon \Mat(\Sph) \too \Mat(\ZZ)
\]
and hence, by restriction of models, a functor
\[
    \sD_{\ge 0}(\ZZ) \too \CGrp(\Spc).
\]

As we shall see later, $\ZZ$ carries a natural \textit{commutative ring} structure in $\CGrp(\Spc)$, and the connective derived $\infty$-category $\sD_{\geq 0}(\ZZ)$ identifies with the $\infty$-category of $\ZZ$-modules therein. Under this identification, the functor above is \textit{restriction of scalars} along the commutative ring morphism $\Sph \to \ZZ$ while its left adjoint is \textit{extension of scalars} $\ZZ\otimes_\Sph M$. 
Thus, quite tautologically, homotopy-coherent commutative groups may be regarded as $\Sph$-modules, and the resulting theory as commutative algebra with $\ZZ$ replaced by the deeper base $\Sph$.

\subsection{Levels of commutativity}

We have seen that homotopy-coherent groups are modeled geometrically by loop spaces. It is natural to ask whether \textit{commutative} groups admit a similar description. Put differently, if we forget a homotopy-coherent commutative group to its underlying group and then identify that group with a loop space, how does the extra commutative structure manifest itself geometrically?
The answer comes from what might initially seem like a completely unrelated approach to homotopy-coherent commutativity. 

The celebrated \textit{Eckmann-Hilton argument} implies the following classical fact: group objects in groups (in sets) are \textit{abelian} groups,
\[
    \Grp(\Grp(\Set)) \simeq \CGrp(\Set).
\]
By analogy, one could naively try and define homotopy-coherent commutative groups as group objects in~$\Grp(\Spc)$. By iterating the identification $\Grp(\Spc) \simeq \Spc_{*,\ge 1}$ from \Cref{Grp_Loop} we get the following:

\begin{prop}
    There is an equivalence of $\oo$-categories
    \[
        \Grp(\Grp(\Spc)) \simeq \Spc_{*,\ge 2}.
    \]
\end{prop}
On viewing $X\in \mathcal S_{*,\ge 2}$ as an object of $\Grp(\Grp(\mathcal S))$ and forgetting both group structures, its underlying space is $\Omega^2X$.  
In other words, these \textit{double homotopy-coherent groups} are precisely \textit{double loop spaces}. The multiplication on $\Omega^2 X$ is still given by concatenation (in either direction) and it is commutative \textit{up to homotopy}. To wit, 
\[
    \pi_0 (\Omega^2 X) \ \simeq\ \pi_2 (X) \qin \Ab.
\]

\begin{example}
    One way to get an ordinary commutative group is to start with an arbitrary group $G$ and take its center $Z(G)$. The homotopy-coherent generalization of this operation rests on the observation that $Z(G)$ is the fixed point set of the conjugation action of $G$ on itself. For a pointed connected space $X$, the conjugation action of $\Omega X$ on itself is encoded by the free loop space fiber sequence (see \Cref{Conj_LX})
    \[
        \Omega X \too LX \too X.
    \]
    The \tdef{homotopy-coherent center} of $\Omega X$ is defined as the \textit{homotopy fixed points} of this action, which can be computed as the sections of the projection $LX \to X$. It is a simple exercise to realize this space as a \textit{double loop space}:\footnote{It is important that the automorphisms of $X$ are as an \textit{unpointed} space.}
    \[
        \mdef{Z}(\Omega X) := 
        (\Omega X)^{h(\Omega X)} \simeq 
        \Gamma_X(LX) \simeq 
        \Omega \Aut(X) \simeq 
        \Omega^2 B\Aut(X). 
    \]
    When $X= BG$ for an ordinary group $G$, the space $B\Aut(BG)$ is a 2-type whose $\pi_2$ is $Z(G)$. Consequently, its double loop space is the discrete group $Z(G)$, recovering the ordinary notion of the center.
\end{example}

Although commutative up to homotopy, $\Omega^2 X$ is not commutative up to full \textit{coherent homotopy} as higher ambiguities remain in the identifications between products differing only by a permutation of their factors.
To see what remains after the first level of commutativity has been imposed, begin one stage lower. In an ordinary loop space $\Omega X$, the two orders of multiplication $m(x,y)=xy$ and $m(y,x)=yx$ usually produce different elements. Their quotient, the \textit{commutator} $xyx^{-1}y^{-1}$, is a new element in $\Omega X$, witnessing the (potential) noncommutativity of $m$. This commutator construction
\[
    \Omega X \times \Omega X \too \Omega X \qin \Spc_*
\]
is corepresented by the (nontrivial) attaching map 
\[
    \mdef{\beta_1} \colon S^1 \too S^1\vee S^1 \qin \Spc_*
\]
of the 2-cell in the torus $S^1\times S^1$.\footnote{As $\Omega(S^1 \times S^1) = \ZZ\times \ZZ$, it corepresents pairs of \textit{commuting} group elements.} 

The \textit{suspension} of this map however \textit{is} null-homotopic. It translates to the fact that in a \textit{double} loop space $\Omega^2 X$, the commutator of $x$ and $y$ is homotopically trivial because $m(x,y) \sim m(y,x)$. This is the classical Eckmann--Hilton argument; the extra dimension allows one to ``rotate'' the picture and swap the positions of $x$ and $y$. However, it is trivial in \textit{two} different ways, corresponding to the two different ways to identify $m(x,y)$ with $m(y,x)$ --- either by rotating the picture $180^\circ$ \textit{clockwise}, or by rotating it $180^\circ$ \textit{counter-clockwise}. Concatenating one with the inverse of the other produces a \textit{loop} which is graphically represented by a full circle rotation: 

\begin{center}
    \includegraphics[width=0.8\textwidth]{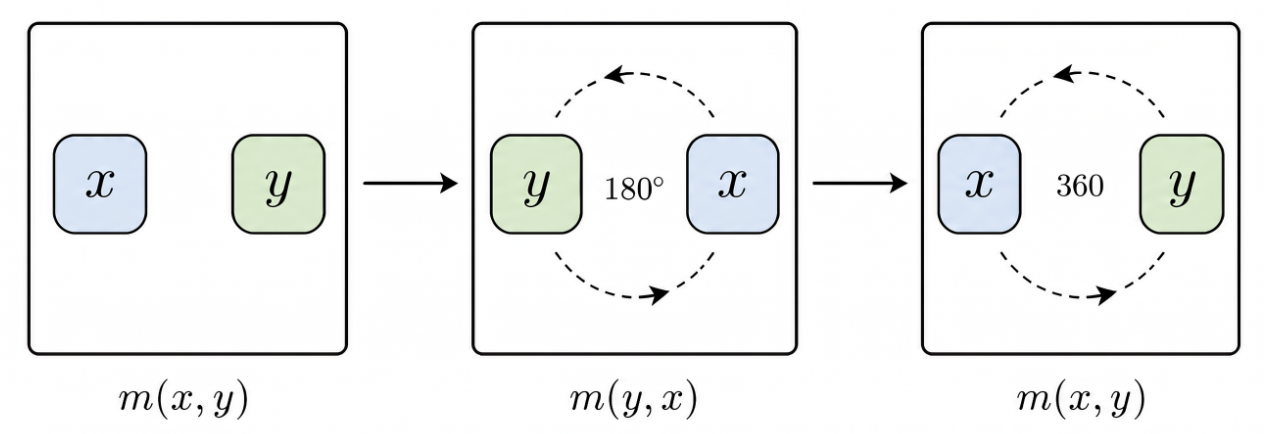}
\end{center}
\vspace{-10pt}

Multiplying by the inverse of $m(x,y)$ shifts this loop to the basepoint of $\Omega^2 X$, producing an element of $\Omega^3 X$. The resulting \tdef{higher commutator} map 
\[
    \Omega^2X \times \Omega^2 X \too \Omega^3 X \qin \Spc_*
\]
is corepresented by the attaching map 
\[
    \mdef{\beta_2\colon} S^3 \too S^2\vee S^2 \qin \Spc_*
\]
of the 4-cell in $S^2 \times S^2$. It measures the ambiguity in \textit{how} $m(x,y)$ and $m(y,x)$ are identified. 

Now, the \textit{suspension} of $\beta_2$ is, again, null-homotopic. This translates to the fact that in a \textit{triple} loop space, there is an extra dimension to produce a homotopy between the two ways to identify $m(x,y)$ and $m(y,x)$. 
One could rotate the entire path through the northern or southern hemisphere; comparing them produces the \textit{next ambiguity}. 

The pattern is now visible. At level $n$, the preceding ambiguity in the homotopy-coherence of commutativity has become null-homotopic, but the \textit{space} of its null-homotopies carries a new ambiguity, measured by the attaching map\footnote{The maps $\beta_n$, and the related operations they induce on homotopy and homology, appear under various names, including Whitehead products, Samelson products, and Browder brackets.}
\[
    \mdef{\beta_{n}} \colon S^{2n-1} \too S^n \vee S^n \qin \Spc_*
\]
of the top cell in $S^n\times S^n$.

To summarize, an $n$-fold loop space carries not merely a multiplication, but $n$ mutually interchanging multiplications:
\[
    \Spc_{*,\ge n} \simeq
    \overset{n\text{  times}}{\overbrace{\Grp(\Grp(\dots\Grp}}(\Spc)\dots)).
\]
Each additional loop coordinate supplies a new multiplication operation, which identifies via interchange with the previous ones, but in doing so resolves one further layer of ambiguity in the preceding commutativity data. Thus, $\Spc_{*,\ge n}$ may be viewed as the $\oo$-category of homotopy-coherent groups with $n$ levels of coherent commutativity, and
\[
    \Omega \colon \Spc_{*,\ge n} \too \Spc_{*,\ge n-1}
\]
as forgetting the highest level. 

\subsection{Infinite loop spaces}

The natural next step is to push this process all the way to its limit. 

\begin{defn}
    The $\oo$-category of \tdef{connective spectra} is given by the limit of the following tower of $\oo$-categories:
    \[
        \mdef{\Sp_{\ge 0}} := 
        \invlim_{n \in \NN}\big(
            \dots 
            \oto{\Omega} \Spc_{*, \ge 2}
            \oto{\Omega} \Spc_{*, \ge 1}
            \oto{\Omega} \Spc_{*, \Gray{\ge 0}}
        \big)
        \qin \Cat_{\oo}.
    \]
\end{defn}

Such \textit{infinite loop spaces} are to be thought of as \textit{infinitely commutative groups} in spaces. In such a group, the space of operations realizing a given finite product is \textit{contractible} (but see \Cref{Eta_Comm}). Quite remarkably, this notion turns out to \textit{coincide} with the one obtained from commutative monoids where the discrete semiring $\NN$ has been replaced by its groupoid refinement $\MM$.

Indeed, for $G$ in $\CGrp(\Spc)$ we can iterate the bar construction, viewed as internal suspension, to produce a sequence of groups
\[
    B^nG \qin \CGrp(\Spc) \qquad\text{for}\qquad n\ge 0,
\] 
whose underlying spaces are \textit{$n$-connective}. By \Cref{B_CGrp}, each is the loop space of the next one $B^nG =\Omega B^{n+1} G$. 
\begin{thm}\label{CGrp_Sp}
    The iterations of $B$ on $\CGrp(\Spc)$ followed by the forgetful functor to $\Spc_*$ assemble into an equivalence of $\oo$-categories
    \[
        B^{\mdef{\oo}} \colon \CGrp(\Spc) \iso \Sp_{\ge 0}.
    \]
\end{thm}

The argument is surprisingly simple given what we already know.

\begin{proof}
    Our strategy is to consider both structures --- that of a commutative group and that of an infinite loop space --- \textit{at once}, and show that both forgetful functors
     \[
        \Sp_{\ge 0} \froom 
        \CGrp(\Sp_{\ge 0}) \too
        \CGrp(\Spc)
     \]
     are \textit{equivalences}. 

    We first show that $\Sp_{\ge 0}$ is \textit{additive}, and so by \Cref{CGrp_add_Iff} absorbs the external commutative group structure.
    The levelwise point is a \textit{zero object} and products are computed levelwise. The main trick is the \textit{internal delooping}. For every $X \in \Sp_{\ge 0}$,
    \[
        X \simeq \Omega BX 
        \quad\text{where}\quad
        (BX)_n = X_{n+1}.
    \]
    Thus,
    \[
        \pi_0\Map_{\Sp_{\ge 0}}(X,Y) \simeq
        \pi_0\Map_{\Sp_{\ge 0}}(X,\Omega^2 B^2Y) \simeq
        \pi_2\Map_{\Sp_{\ge 0}}(X,B^2Y).
    \]
    These mapping sets therefore admit a natural abelian group structure with bilinear composition; the homotopy category of $\Sp_{\ge0}$ is hence $\Ab$-enriched. By ordinary category theory this implies that it, and therefore $\Sp_{\ge0}$, are additive.  

    Now, \Cref{Bar_Omega_CGrp} implies that $\Omega$ induces an equivalence
    \[
        \Omega \colon 
        \CGrp(\Spc_{*,\ge n+1}) \iso 
        \CGrp(\Spc_{*,\ge n})
        \quad\text{for all}\quad n\ge0.
    \]
    The construction $\CGrp$ preserves limits of $\oo$-categories along product preserving functors so
    \[
        \CGrp(\Sp_{\ge 0}) \iso 
        \CGrp(\invlim \Spc_{*,\ge n}) \iso
        \invlim \CGrp(\Spc_{*,\ge n}) \iso
        \CGrp(\Spc_{*}) \iso
        \CGrp(\Spc).
    \]
    For the last equivalence, note that the \textit{basepoint} is forced to be the \textit{unit} of the commutative group structure and is therefore automatically preserved by morphisms. Thus, $\CGrp(\Spc)$ absorbs the infinite loop structure. 
    The inverse of the above equivalence, followed by the forgetful functor to $\Sp_{\ge0}$, sends a homotopy-coherent commutative group $G$ to the sequence of iterated deloopings $B^nG$. This is precisely~$B^\oo G$.

\end{proof}

The construction makes it clear that forgetting to the last stage of the tower
\[
    \Omega^{\mdef{\oo}} \colon \Sp_{\ge 0} \too \Spc_{*,\ge 0} = \Spc_*
\]
corresponds, under composition with $B^\oo$, to forgetting the entire commutative group structure save for the unit element. Furthermore, the projection to the \textit{next-to-last} stage of the tower
\[
    \Omega^{\mdef{\infty-1}} \colon \Sp_{\ge 0} \too \Spc_{*,\ge 1} \simeq \Grp(\Spc)
\]
corresponds to forgetting the \textit{commutativity data} of the group structure. The higher stages of the tower track the gradual accumulation of higher levels of commutativity  ultimately assembling into a lift to $\CGrp(\Spc)$. We conclude this discussion with some examples, variations, and observations. 

\begin{example}
    It is instructive to consider the case where the underlying space is truncated. 
    Group structures on a \textit{set} $G$ correspond bijectively to deloopings $G \simeq \Omega X$ for a pointed connected space $X$.
    Indeed, such a delooping identifies $G$ with $\pi_1X$, while a group structure on $G$ gives $X \simeq BG$.
    The existence of a \textit{second} delooping $B^2G$ is then the \textit{property} of the group structure on $G$ being \textit{commutative}. 
    In that case all further deloopings $B^nG$ exist uniquely and assemble into the \tdef{Eilenberg--Mac Lane spectrum}~$B^\oo G$.\footnote{Classically denoted $HG$.} This recovers the classical fact with which we began our discussion:
    \[
        \Grp(\Grp(\Set)) \simeq \CGrp(\Set).
    \]
    If $\cG$ is a \textit{groupoid}, viewed as a $1$-truncated space, a delooping $B\cG$ corresponds to a \textit{groupoidal}\footnote{Ie, \textit{monoidal}, where every object is tensor invertible.} structure on $\cG$. A \textit{second} delooping $B^2\cG$ corresponds to a \textit{braiding}, while a \textit{third} exists precisely when the braiding is \textit{symmetric}. From that point onwards all further deloopings exist uniquely. For \textit{2-groupoids} the analogous hierarchy has one extra rung, bearing the intriguing name \textit{sylleptic}. 
\end{example}

\begin{rem}\label{En_Mon}
    One may similarly consider $n$-fold \textit{monoid} objects in an $\oo$-category $\cC$, called \tdef{$\EE_n$-objects}:
    \[
        \Mon_{\mdef{\EE_n}}(\cC) = 
        \overset{n\text{  times}}{\overbrace{\Mon(\Mon(\dots\Mon}}(\cC)\dots)).
    \]
    These are governed by their own, somewhat more complicated,  $\oo$-categorical algebraic theories, which we shall revisit later. For $\EE_n$-spaces, the ones whose $\pi_0$ is a group are exactly the $n$-fold loop spaces. It is also true that the limiting case is just commutative monoids,
    \[
        \CMon(\cC) \simeq \invlim_{n \in \NN}\ \Mon_{\EE_n}(\cC) \qin \Cat_\oo.
    \]
    This is already a fixed point for the process
    \[
        \Mon(\CMon(\cC)) \simeq \CMon(\cC). 
    \]
\end{rem}

\begin{rem}\label{Eta_Comm}
    There is a subtle trap in interpreting the aforementioned \textit{contractibility} of the space of $n$-fold products of an infinite loop space.
    For $n=2$, the composite of the higher commutator with the fold map
    \[
        S^3 \oto{\ \beta_2 \ } 
        S^2 \vee S^2 \oto{\ \nabla \ }
        S^2
    \]
    is equal to  $2\eta$, where $\eta \in \pi_3(S^2)$ is the \textit{Hopf map}. 
    Composition with the fold map amounts to restriction along the diagonal $x=y$. On the diagonal, however, the \textit{half-turn} exchanging the two inputs already closes up to a loop at $m(x,x)$. After shifting to the identity, this loop is corepresented by $\eta$ \textit{itself}. 
    Furthermore, although $2\eta$ vanishes under the suspension map
    \[
        \ZZ \simeq \pi_3(S^2) \oto{\ \Sigma \ } \pi _4(S^3) \simeq \ZZ/2,
    \]
    the class $\eta$ \textit{survives}. Since $\pi_4(S^3)$ is already in the \textit{stable range}, it survives \textit{all the way up}.  
    This reveals an unintuitive yet important feature of homotopy-coherent commutativity. The square $x^2$ of any element $x$ is identified with itself in two essentially different ways: tautologically, and through the commutativity between the two copies of $x$. Because commutativity is \textit{structure} rather than a mere property, these two identifications need not agree.\footnote{A reader who found the puzzle in \Cref{Coherence_Puzzle} challenging may want to revisit it now.}      
\end{rem}

\subsection{Rotating loops}

As with the loop space model for noncommutative groups, the infinite loop space model for commutative groups connects homotopy-coherent algebra to classical homotopical and geometric notions. To begin with, the underlying space of every homotopy-coherent commutative group $G$ is an $n$-fold loop space for all $n\ge0$,
\[
    G \simeq \Omega^n B^n G := \Map_*(S^n,B^n G) \qin \Spc_*.
\]
As such, it admits a canonical action of the orthogonal group $O(n)$ via its natural (pointed) action on $S^n \simeq (\RR^n)^+$. These actions are compatible and thus assemble to an action of the \textit{stable orthogonal group}, 
\[
    O(\oo) = \colim_{n \in \NN}\ O(n).
\]
The action of $O(1) \simeq \ZZ/2$ is simply the involution on $G$ given by taking an element to its \textit{inverse}. There is a certain mind-bending fact about the group $O(\oo)$ arising from contemplating the following variant of topological $K$-theory:

\begin{defn}\label{O_Comm}
    Replacing complex vector spaces by real ones gives connective 
    \tdef{real topological $K$-theory},
    \[
        \mdef{\ko} := 
        M_\RR^{\gp} \simeq 
        \ZZ \times BO(\oo) \qin \CGrp(\Spc).
    \]
\end{defn}
Looping inside $\CGrp(\Spc)$ yields 
\[
    \Omega \ko \simeq 
    \Omega BO(\oo) \simeq 
    O(\oo),
\]
and thereby equips $O(\oo)$ with a homotopy-coherent \textit{commutative} group structure. This is puzzling as matrix multiplication is certainly not commutative. The apparent escape is that the commutative operation inherited from $\ko$ is not matrix multiplication at all, but the one induced by the \textit{block sum}
\[
    O(n)\times O(m) \intoo O(n+m) \intoo O(\oo).
\]
Yet the two operations \textit{interchange} and hence \textit{coincide} by a homotopy-coherent version of the Eckmann-Hilton argument (see \Cref{En_Mon}). Thus, matrix multiplication on $O(\oo)$ \textit{is}, in a concrete and natural way, \textit{homotopy-coherently} commutative.

The left adjoint of $\Omega^\oo$, whose image consists of \tdef{suspension spectra}, is denoted by
\[
    \mdef{\Sigma^\oo} \colon \Spc_* \too \Sp_{\ge 0},
\]
and bears a similar relation to the suspension functor $\Sigma$, as $\Omega^\oo$ to the loops functor $\Omega$. 
Their composition has the following geometric model:
\[
    \Omega^\oo\Sigma^\oo X \simeq 
    \colim_{k\in \NN}\ \Omega^k \Sigma^k X \qin \Spc_*.
\]

\begin{example}\label{Stable_Stems}
    The \textit{sphere spectrum} is therefore given by 
    \[
        \Sph \simeq \Sigma^\oo S^0 \qin \Sp_{\ge0}
    \]
    and its underlying space is 
    \[
        \Omega^\oo \Sph \simeq 
        \colim_{k\in \NN}\ \Omega^k S^k \qin \Spc_*,
    \]
    explaining the sense in which it is ``built from spheres''. In particular, its homotopy groups are the \textit{stable homotopy groups of spheres}\footnote{This sequence actually stabilizes at $k=n+2$ by the \textit{Freudenthal suspension theorem}.}
    \[
        \mdef{\pi_n^{s}} := 
        \colim_{k\in \NN}\ \pi_{n+k}S^k
        \qin \Ab.
    \]
    These groups appear throughout topology and geometry, notably in (framed) cobordism and the study of exotic spheres. It is no small matter that such an inherently geometric structure emerges from the symmetries of finite sets. Much is known and unknown about these groups.\footnote{For the current state of knowledge, see \cite{isaksen2020stable}.}
\end{example}

\subsection{A minor point}

The need to pass back and forth between pointed and unpointed constructions is a constant source of irritation in algebraic topology. From the algebraic perspective this can be explained as follows. The usual forgetful functor from monoids to sets factors through \textit{pointed} sets by remembering the unit element,
\[
    \Mon \too \Set_* \too \Set.
\]
The left adjoint of the first functor takes the free monoid on a pointed set and identifies its basepoint with the unit. Since every pointed set is of the form $X_+=X\sqcup{\mathrm{pt}}$, this simply computes the free monoid on $X$. The same factorization exists in the homotopy-coherent setting:
\[
    \Mon(\Spc) \oto{\ U^\pt_\Mon\ } 
    \Spc_* \oto{\ U_\pt\ } 
    \Spc.
\]
The left adjoint of $U_{\pt}$ still adjoins a disjoint basepoint, $X\mapsto X_+$, but is no longer \textit{essentially surjective}. Accordingly, the left adjoint $F_{\Mon}^{\pt}$ is now more interesting. Starting from
\[
    F_{\Mon}(X) \simeq \bigsqcup_{n\in \NN} X^n
\]
one has to identify coherently all the $n$-tuples containing $k$ basepoints with the corresponding $(n-k)$-tuples obtained by deleting those coordinates.\footnote{This is classically known as the \textit{James construction}.} The upshot is that if $X$ is \textit{connected} then $F_\Mon^\pt(X)$ is also connected and so automatically a \textit{group}. Thus, the pointed setting offers substantially more flexibility in homotopy-coherent algebra than in ordinary algebra.
 
These considerations carry over, with appropriate changes, to the \textit{commutative} setting. The free homotopy-coherent commutative group on an (unpointed) space is given by the composition
\[
    \Sigma^\oo_{\mdef{+}} \colon 
    \Spc \oto{\ (-)_+\ } 
    \Spc_* \oto{\ \Sigma^\oo \ }
    \Sp_{\ge 0}.
\]
As it takes a space $X$ to the constant $X$-shaped colimit on $\Sph$ it might also be denoted by\footnote{Or $X \otimes \Sph$, alluding to the symmetric monoidal structure on (connective) spectra.}
\[
    \mdef{\Sph\langle}X\mdef{\rangle} := \Sigma^\oo_+X \qin \Sp_{\ge 0}.
\]
In this setting, the difference between the pointed and the unpointed constructions is rather small. Any point $\pt \to X$ has a unique retraction $X\to \pt$, which in the additive $\oo$-category $\Sp_{\ge 0}$ \textit{splits} 
\[
    \Sigma^\oo_+X \simeq \Sigma^\oo X\oplus \Sph \qin \Sp_{\ge 0}.
\]

\section{Spectra}

As their name betrays, \textit{connective spectra} are a special case of the more general notion of a \textit{spectrum}. Although this lies somewhat off our main path, being not governed by an algebraic theory, it would be remiss not to say something about it, given its prominence and how well it complements the picture we have developed thus far. We shall not, however, attempt to give a full `introduction to stable homotopy theory' as it would take us too far afield. 

\subsection{Desuspension of belief}

 Like many fundamental notions in mathematics, spectra admit several descriptions that arise from quite different perspectives and are not obviously equivalent.\footnote{And, like many of them, it bears a remarkably generic, overused, and uninformative name.} The most natural of these from our perspective is to simply remove the connectivity constraints in the definition of a connective spectrum.

\begin{defn}
    The $\infty$-category of \tdef{spectra} is given by the limit of the following tower of $\oo$-categories:\footnote{This follows closely the classical notion of an \textit{$\Omega$-spectrum}, only it is an $\infty$-category out of the box, with mapping spaces, (co)limits, and all.} 
    \[
        \mdef{\Sp} := \invlim (\dots \oto{\ \Omega\ } \Spc_* \oto{\ \Omega\ } \Spc_* \oto{\ \Omega\ } \Spc_*)
        \qin \Cat_\oo.
    \] 
\end{defn}

We will identify $\Sp_{\ge0}$ with the full subcategory of $\Sp$ in the obvious way. The fully faithful inclusion $\Sp_{\ge 0} \into \Sp$ admits a \textit{right} adjoint called \tdef{connective cover}
\[
    \mdef{\tau_{\ge 0}} \colon \Sp \too \Sp_{\ge 0},
\]
which levelwise replaces $X_n$ with $\tau_{\ge n}X_n$.
The zeroth-space functor $\Omega^\oo \colon \Sp \to \Spc_*$ factors as a composite
\[
    \Sp \oto{\ \tau_{\ge 0}\ } 
    \Sp_{\ge 0} \oto{\ \Omega^\oo \ }
    \Spc_*
\]
and is \textit{not} conservative as it forgets all the information \textit{below zero}. The kernel consists of the \tdef{strictly coconnective} spectra. 
To restore conservativity one has to consider the \textit{$\ZZ$-graded} homotopy groups. 

\begin{defn}
    For $X \in \Sp$ and $n\in \ZZ$ we set
    \[
        \mdef{\pi_n} X := \pi_{n+k} (X_k) 
        \quad\text{for}\quad 
        k\gg 0. 
    \]    
\end{defn}
The greatest advantage of working with $\Sp$ over $\Sp_{\ge 0}$ is that the functor $\Omega$ on it is \textit{invertible}. Indeed, we have\footnote{Observe, on the other hand, that even when $X$ is connective, $\Omega X$ need not be.}
\[
    (\Omega X)_n \simeq\Omega X_n 
    \qquad\text{and}\qquad 
    (\Omega^{-1}X)_n \simeq X_{n+1}.
\]
The inverse of $\Omega$ is its adjoint $\Sigma$ and accordingly $\Sigma^{\mdef{-1}} := \Omega$ is known as \tdef{desuspension}. It is convenient to extend this notation to\footnote{For $n \in \NN$ we have $\Sigma^\oo \Sigma^n = \Sigma^n\Sigma^\oo$ as $\Sigma^\oo$ is a left adjoint and hence preserves colimits.} 
\[
    \Sigma^{\mdef{\oo + n}} := 
    \Sigma^n\Sigma^{\oo} 
    \quad\text{for all}\quad n \in \ZZ.
\]

\begin{rem}
    Some point-set models of spectra are based on the notion of a \tdef{prespectrum}, which is a sequence of pointed spaces $X_n$ and maps $\Sigma X_n\to X_{n+1}$. From our perspective, the spectrum represented by such data is the sequential colimit
    \[
        \colim \big(
            \Sigma^\oo X_0 \too 
            \Sigma^{\oo-1} X_1 \too 
            \Sigma^{\oo-2} X_2 \too 
            \dots
        \big)
        \qin \Sp.
    \]
    The formal fact underlying this approach is that \textit{every} spectrum can be constructed in such a way. Namely, it can be approximated by desuspensions of suspension spectra. 
\end{rem}

The passage from connective spectra $\Sp_{\ge 0}$ to all spectra $\Sp$ is analogous to the passage from the connective derived category $\sD_{\ge 0}(R)$ to the unbounded derived category $\sD(R)$. The upshot is that $\Sp$, like  $\sD(R)$, is \textit{stable}. 

\begin{defn}
    An $\oo$-category $\cC$ is \tdef{stable} if it admits a zero object, finite (co)limits, and the loop--suspension adjunction is an equivalence\footnote{Other than the obvious \cite[\S 1]{HA}, we recommend \cite{harpaz2013introduction} as an introduction to the notion of stability.}  
    \[
        \Sigma\colon \cC \isoadj \cC \noloc \Omega.
    \]
\end{defn}

Stable $\oo$-categories are an $\oo$-categorical analogue of ordinary \textit{abelian categories}. 
\begin{thm}
    An $\oo$-category $\cC$ is \textit{stable} if and only if it is  \textit{additive}, every morphism admits a fiber and a cofiber, and it satisfies a strong form of the \textit{first isomorphism theorem}: a null sequence
    \[
        Z \too Y \too X \qin \cC
    \]
    is a \textit{cofiber} sequence if and only if is a \textit{fiber} sequence.
\end{thm}
Such a (co)fiber sequence can be extended indefinitely by (co)fibers in both directions:\footnote{Taking $\pi_0$ one gets a long exact sequence in \textit{homotopy groups}.}
\[
    \dots \too \Sigma^{-1} Y \too \Sigma^{-1} X \too 
    Z \too Y \too X \too 
    \Sigma Z \too \Sigma Y \too \dots
    \qin \cC.
\]
In particular, \textit{extensions} of $X$ by $Z$ are classified by \textit{maps} $X \to \Sigma Z$. Since 
\[
    \Map_\cC(X,Z) \simeq \Omega^n\Map_\cC(X,\Sigma^n Z),
\]
these shifted mapping spaces assemble naturally into a \textit{spectrum},
\[
    \mdef{\hom}_\cC(X,Z)_n := \Map_\cC(X,\Sigma^n Z) \qin \Spc_*.
\]
Thus, $\Sp$ plays for stable $\oo$-categories the role that $\CMon(\Spc)$ plays for \textit{semiadditive} ones and $\CGrp(\Spc)$ for \textit{additive} ones: every stable $\oo$-category is canonically enriched in spectra.
The \textit{negative} homotopy groups of the mapping spectrum 
\[
    \mdef{\Ext^n}(X,Z) := 
    \pi_{-n} \hom(X,Z)
    \quad\text{for}\quad n \in \NN
\]
are a natural generalization of the familiar $\Ext^n_R$-groups for the $\oo$-category $\sD(R)$; they encode extension information.  

\begin{rem}
    To dispel any confusion that may have arisen, we reiterate that the ordinary category of abelian groups is \textit{fully faithfully} embedded in the $\oo$-category of spectra via
    \[
        \Ab \simeq 
        \CGrp(\Set) \intoo 
        \CGrp(\Spc) \simeq 
        \Sp_{\ge 0} \intoo
        \Sp.
    \]
    The $\Ext$-groups appear as the \textit{negative} homotopy groups of the mapping \textit{spectrum} $\hom_{\Sp}(A,B)$, but are invisible to the mapping \textit{space}, which is discrete:
    \[
        \Map_{\Sp}(A,B) \simeq
        \Omega^\oo\hom_{\Sp}(A,B) \simeq 
        \hom_{\Ab}(A,B).
    \]
\end{rem}

\begin{example}
    For a pair of abelian groups $A$ and $B$ viewed as objects of $\sD(\ZZ)$ the corresponding groups $\Ext^n_\ZZ(A,B)$ are zero for $n\ge 2$ and $\Ext^1_\ZZ(A,B)$ classifies ordinary abelian group extensions 
    \[
        0 \too B \too E \too A \too 0.
    \]
    As objects of $\Sp$, the group $\Ext_\Sph^1(A,B)$ classifies analogous extensions but in \textit{spectra}. Since $A$ and $B$ are discrete, 
    the long exact sequence of homotopy groups shows that $E$ is discrete as well. We deduce that 
    \[
        \Ext_\Sph^1(A,B) \ \simeq\ \Ext_\ZZ^1(A,B) \qin \Ab.
    \]
    However, the groups $\Ext_\Sph^n(A,B)$ can be \textit{nonzero} for arbitrarily large $n$. In fact, the graded algebra 
    \[
        \sA := \bigoplus_{n\in \NN}\Ext^n_\Sph(\FF_p,\FF_p)
    \]
    is precisely the classical \textit{Steenrod algebra} of stable $\FF_p$-cohomology operations.\footnote{A modern exposition can be found in Lurie's \href{https://ocw.mit.edu/courses/18-917-topics-in-algebraic-topology-the-sullivan-conjecture-fall-2007/pages/lecture-notes/}{notes} on the Sullivan conjecture.} 
\end{example}

\begin{rem}
    The homotopy category $h\cC$ of a stable $\oo$-category $\cC$ has a canonical structure of a \textit{triangulated category}. Roughly speaking, the idea is to remember the cofiber sequences of $\cC$ as additional structure on $h\cC$, in the form of diagrams \(X \to Y \to Z \to \Sigma X\) known as \textit{distinguished triangles}. One can do quite a lot with this 1-categorical framework, including constructing and computing spectral sequences, but it has its limitations, as it does not retain the higher coherence data.\footnote{See \cite[Theorem 1.1.2.14 and Remark 1.1.2.15]{HA}.} 
\end{rem}

\subsection{Cohomology theories}

Negative homotopy groups also appear when spectra serve as representing objects for \textit{cohomology theories}. For a space $X$ and a spectrum $E$ we denote by $E^X$ the constant $X$-shaped \textit{limit} on $E$. 

\begin{defn}
    For a space $X$ and a spectrum $E$, the (unreduced) 
    \tdef{$E$-cohomology} groups of $X$ are given by\footnote{This is the origin of the historical mistake known as ``cohomological grading''.}   
    \[
        \mdef{E^n}(X) := 
        \pi_{-n} (E^X)
        \quad\text{for all}\quad n\in \ZZ.
    \]
    The \textit{reduced} version $\mdef{\tilde{E}}^n(X)$ for a \textit{pointed} space $X$ is obtained by taking the negative homotopy groups of the fiber of the map 
    \[
        E^X \too E^\pt = E \qin \Sp.
    \]
\end{defn}

The usual (singular) cohomology is recovered as a special case.  
\begin{example}
    For an ordinary abelian group $A$ viewed as a (discrete connective) spectrum,
    \[
        \pi_{-n}(A^X) \simeq 
        \pi_{-n}\hom(\Sph\angles{X}, A) \simeq  
        \pi_0\Map(\Sph\angles{X}, \Sigma^nA) \simeq
    \]
    \[
        \pi_0\Map(X, \Omega^\oo\Sigma^nA) \simeq 
        \pi_0\Map(X, B^n A) \simeq
        H^n(X ; A).
    \]    
    Allowing the $X$-indexed coefficient system of abelian groups to \textit{vary} similarly recovers cohomology with \textit{local coefficients}.
\end{example}

By the long exact sequence of homotopy groups for a (co)fiber sequence of spectra, the induced sequence of functors of (ordinary) categories
\[
    \tilde{E}^n \colon h\Spc_*^\op \too \Ab 
    \qquad\text{for}\qquad
    n \in \ZZ,
\]
is equipped with canonical natural isomorphisms 
\[
    \tilde{E}^{n+1}(\Sigma X) \simeq 
    \Tilde{E}^n(X) \qin \Ab
\]
and satisfies all the Eilenberg--Steenrod axioms, except for the dimension axiom,\footnote{$\tilde{E}^n(S^0)=0$ for $n\neq 0$.} making it a \textit{generalized cohomology theory}.
The celebrated \textit{Brown representability} theorem provides a converse.

\begin{thm}[Brown]
    Every generalized cohomology theory $\{E^n\}_{n\in \ZZ}$ arises in the way above from a, unique up to equivalence, spectrum $E$. 
\end{thm}

Historically, generalized cohomology theories were one of the principal ways in which spectra came to be understood: as functorial homotopy invariants of spaces which could be assembled from local geometric pieces by excision and the resulting long exact sequences.
From the point of view developed here, however, this is somewhat backwards. The $\oo$-category of spectra emerges naturally from contemplating homotopy-coherent commutative groups. The associated cohomology theory is obtained by probing the spectrum with spaces, recording only the homotopy groups of the result. In that sense, regarding a spectrum primarily as a cohomology theory is quite an incidental way of looking at it. Case in point, ordinary abelian groups appear from this perspective in the rather roundabout guise of \textit{singular cohomology} with coefficients.

\begin{rem}
    There is also a real mathematical deficiency in the cohomology theory perspective, which is not a mere question of psychology. While Brown representability ensures that cohomology theories capture correctly the \textit{objects} (ie spectra), they fail to capture the right notion of a \textit{morphism}, even up to homotopy. There exist nonzero maps of spectra that induce the zero natural transformation on the associated cohomology theories. In fact, there is such a map from $\KU$ to~$\Sigma \ZZ$. This is closely related to (though not quite the same as) the phenomenon of \textit{phantom maps}.  
\end{rem}

None of this is meant to suggest that the cohomology theory perspective has become obsolete. For certain spectra~$E$, classes in $E^n(X)$ admit profound geometric, categorical, or even analytic interpretations that are not immediate from $E$ itself as a homotopy-theoretic object. Indeed, some of the most beautiful mathematical stories arise precisely where these two viewpoints meet: the intrinsic structure of spectra, governed by homotopy-coherent commutative algebra and algebraic geometry, and the concrete incarnations of them as cohomology theories elsewhere in mathematics.\footnote{One of the best such stories is discussed in \cite{hopkins1995topological}.}

In a different direction, representability can also manufacture spectra from properties of functors, even when no coherent construction of the desired spectrum is apparent.

\begin{example}
    The collection of functors 
    \[
        \mdef{I_{\QQ/Z}^n}(X) \:= 
        \hom_{\Ab}(\pi_n \Sigma^\oo X, \QQ/\ZZ)
    \]
    is a generalized cohomology theory, owing to the \textit{divisibility} (and hence injectivity) of the abelian group $\QQ/\ZZ$. It is therefore represented by a \textit{coconnective} spectrum $\mdef{I_{\QQ/\ZZ}}$ known as the \tdef{Brown-Comenetz} spectrum. 
\end{example}
The operation taking a spectrum $X$ to 
\[
    I_{\QQ/Z} X := \hom_{\Sp}(X,I_{\QQ/\ZZ}) \qin \Sp
\]
is a spectral form of the classical Pontryagin duality, called, unsurprisingly, \tdef{Brown-Comenetz duality}. Apart from that, $I_{\QQ/\ZZ}$ is the uncontested source of \textit{counterexamples} in stable homotopy theory. 

\subsection{Excisive functors}

In addition to a \textit{cohomology} theory, a spectrum determines also a \textit{homology} theory. For a space $X$ and a spectrum $E$, the constant $X$-shaped \textit{colimit} on $E$ is denoted by
\[
    \mdef{E\langle}X\mdef{\rangle} := 
    \colim _X E \qin \Sp.
\]
When $X$ is pointed we consider the reduced version of the above defined by the cofiber sequence 
\[
    E \too E\angles{X} \too \tilde{E}\angles{X} \qin \Sp.
\]

\begin{defn}
    For a (pointed) space $X$ and a spectrum $E$, the unreduced and reduced  \tdef{$E$-homology} groups of $X$, for $n\in \ZZ$, are given by
    \[
        \mdef{E_n}(X) := 
        \pi_{n} E\angles{X} 
        \qquad\text{and}\qquad
        \mdef{\tilde{E}_n}(X) := \pi_n\tilde{E}\angles{X}.
    \] 
\end{defn}
The homology functors $\tilde E_n$ satisfy properties roughly dual to those of the cohomology functors $\tilde E^n$, but there is a twist. Instead of considering the homotopy groups of the spectrum $\tilde E\langle X\rangle$, one may consider its entire underlying \textit{space} $\Omega^\infty\tilde E\langle X\rangle$. This construction can be thought of as a coherent refinement of the associated homology theory and turns out to capture $E$ completely.

\begin{defn}
    Let $\Spc_*^{\mdef{\fin}}$ be the $\oo$-category of pointed \textit{finite} spaces. A functor 
    \[
        F\colon \Spc_*^\fin \too \Spc
    \]
    is called \tdef{reduced} if it preserves the terminal object $\pt$ and \tdef{excisive} if it takes pushout squares to pullback squares.\footnote{Note that it does \textit{not} take pushouts to pullbacks (which doesn't even make sense).}
\end{defn}

\begin{thm}
    Taking a spectrum $E$ to the functor $\Omega^\oo \tilde{E}\angles{-}$ assembles into an equivalence of $\oo$-categories
    \[
        \Sp \iso
        \Fun_*^{{\rm exc}}(\Spc_*^\fin, \Spc)
        \Gray{\ \sseq\ \Fun(\Spc_*^\fin,\Spc)},
    \]
    where the target is the full subcategory on reduced excisive functors. 
\end{thm}

\begin{proof}[Proof idea.]
    For a spectrum $E$, the functor 
    \[
        F_E(X) := \Omega^\oo(\tilde{E}\angles{X}) \qin \Spc
    \]
    is reduced by construction. Excisiveness follows from combining the following three facts: $\tilde{E}\angles{-}$ preserves pushouts (being a left adjoint), in $\Sp$ pushout squares are pullback squares (as a consequence of stability), and the functor $\Omega^\oo$ preserves pullbacks (being a right adjoint).  

    Conversely, given a reduced excisive functor $F\colon \Spc_*^\fin \to \Spc$, we define the spaces 
    \[
        E_n := F(S^n) \qquad\text{for}\qquad n \in \NN.
    \]
    Since $F$ is \textit{excisive}, we get for all $n\ge 0$ a pullback square 
    \[\begin{tikzcd}
    	{F(S^n)} & {F(\pt)} \\
    	{F(\pt)} & {F(S^{n+1}).}
    	\arrow[from=1-1, to=1-2]
    	\arrow[from=1-1, to=2-1]
    	\arrow["\lrcorner"{anchor=center, pos=0.125}, draw=none, from=1-1, to=2-2]
    	\arrow[from=1-2, to=2-2]
    	\arrow[from=2-1, to=2-2]
    \end{tikzcd}\]
    But $F$ is also \textit{reduced}, so $F(\pt) = \pt$ and we get an identification 
    \[
        E_n = F(S^n) \simeq \Omega F(S^{n+1})=\Omega E_{n+1},
    \]
    turning the sequence $\{E_n\}$ into a \textit{spectrum}. The two constructions are mutually inverse. 
\end{proof}

While this does not quite present spectra as models of an algebraic theory, it comes very close in spirit: it describes them as diagrams of spaces of a prescribed shape satisfying prescribed conditions. In particular, we can consider \textit{spectrum objects} in arbitrary $\oo$-categories.

\begin{defn}
    Let $\cC$ be an $\oo$-category with finite limits.
    The $\oo$-category of \tdef{spectrum objects} in $\cC$ is given by
    \[
        \mdef{\Sp}(\cC)  := \Fun_*^{\rm exc}(\Spc_*^\fin,\cC).
    \]
\end{defn}
Evaluation at $S^0$ gives a forgetful functor
\[
    U_{\mdef{\Sp}} \colon 
    \Sp(\cC) \too 
    \cC,
\]
which is an equivalence if and only if $\cC$ is \textit{stable}. 

A key feature of this model is that the inclusion of excisive functors into the $\infty$-category of \textit{all} functors
\[
    F\colon 
    \Spc_*^{\fin}\too 
    \Spc
\]
admits a \textit{left} adjoint, giving the \emph{excisive approximation} $\mdef{P_1}F$, which serves as a \emph{linearization} of $F$. For reduced $F$, it can be computed by the formula
\[
    P_1F(X)\simeq \colim_{n \in \NN}\ \Omega^nF(\Sigma^nX) \qin\Spc.
\]
More generally, one defines the notion of an \textit{$n$-excisive} functor that can be thought of as a \textit{polynomial} of degree at most $n$. Every functor $F$ as above has a \textit{Taylor series} in the form of a tower of $n$-excisive approximations\footnote{It is natural to include $P_0F$, which is constant on $F(\pt)$.} 
\[
    F \too 
    \dots \too
    P_n F \too 
    \dots \too 
    P_2 F \too 
    P_1 F \Gray{\too
    P_0 F.}
\]
This tower need not converge to $F$, but it often does for sufficiently connective spaces (ie, there is a ``radius of convergence'' around $\pt$). This is the beginning of Goodwillie's \textit{calculus of functors}.\footnote{See \cite{arone2020goodwillie} for an exposition.} The first indication that this analogy should be taken seriously is the following:

\begin{prop}
    For every $n \ge 1$ and a reduced functor 
    \[
        F \colon \Spc_*^\fin \too \Spc,
    \]
    there exists a spectrum $D_n \in \Sp$ with an action\footnote{This is a \textit{naive} action in the sense of \Cref{Naive_Action}.} of $\Sigma_n$ fitting into a fiber sequence
    \[
        \Omega^\oo \big(\tilde{D}_n\angles{X^{\wedge n}}/\!\!/\Sigma_n\big) \too 
        P_nF(X) \too 
        P_{n-1}F(X) \qin \Spc.
    \]
\end{prop}
This should be compared with the ordinary \textit{calculus of functions}. For a smooth function $f\colon \RR \to \RR$, the difference between successive Taylor polynomials has the form 
\[
    \frac{d_nx^n}{n!} = p_n(x) - p_{n-1}(x) \quad\text{for some}\quad d_n \in \RR.
\]
The $\Sigma_n$-spectrum $D_n$ is thus called the $n$-th \tdef{derivative} of $F$ and denoted by $\mdef{\partial^n}F$.

\begin{example}
    The \textit{identity} (ie forgetful) functor $\Spc_*^\fin \to \Spc$ is reduced but \textit{not} excisive. Its excisive approximation $P_1\Id$ is the functor $\Omega^\oo\Sigma^\oo$. Thus, the first derivative is given by $\partial^1\Id \simeq \Sph$. The higher derivatives can be expressed in terms of the \textit{partition complex} associated with the poset $\Pi_n$ of partitions of $\{1,\dots,n\}$ ordered by refinement. The $\Sigma_n$-spectra $\partial^n\Id$ constitute the underlying symmetric sequence of the \textit{spectral Lie operad}\footnote{See \cite[\S3-4]{arone2020goodwillie}.}. 
\end{example}

\subsection{Periodicity}

The dual theory to Postnikov towers allows one to express every spectrum $E$ as a colimit of its \textit{bounded below} covers:
\[  
    E \simeq \colim_{n\in \NN}\ \tau_{\ge -n} E 
    \quad\text{where}\quad
    \tau_{\mdef{\ge -n}} E := \Sigma^{-n}(\tau_{\ge 0}  \Sigma^n E) \qin \Sp_{\ge -n}.
\]
Combined with the equivalences
\[
    \Sigma^n \colon \Sp_{\ge -n} \iso \Sp_{\ge 0},
\]
one might get the impression that spectra as a whole are not much more than a bookkeeping device for the ``fringe effects'' appearing in (co)fiber sequences of \textit{connective} spectra. Such an impression would be misleading: genuinely nonconnective spectra can exhibit phenomena with no connective counterpart. The most striking example is \textit{periodicity}. 

To illustrate this point, recall the connective topological $K$-theory spectrum $\ku$. We have defined it as the group completion of the commutative monoid of $\CC$-vector spaces, but it is instructive to understand it also as an \textit{infinite loop space}, which we can do with the help of the celebrated \textit{Bott periodicity} theorem.

\begin{thm}[Bott]
    There is an equivalence of pointed spaces
    \[
        \Omega U(\oo) \ \simeq\  \ZZ \times BU(\oo).
    \]
\end{thm}
The name `periodicity' stems from the immediate consequence $\Omega^2 U(\oo) \simeq U(\oo)$ and its corollary,
\[
    \pi_n \ku \simeq \begin{cases}
        \ZZ & n \text{ even} \\
        0 & n \text{ odd}. 
    \end{cases}
\]
\begin{proof}
    Because the unitary group lies at the intersection of topology, geometry, analysis, and algebra, Bott periodicity admits an unusually diverse collection of 
    \href{https://mathoverflow.net/questions/8800/proofs-of-bott-periodicity/9151?utm_source=chatgpt.com}{proofs}.
    The classical geometric one studies the identity component of the loop space 
    \[
        \Omega U(\oo) \simeq \colim \Omega U(n)
    \]
    through the geometry of geodesics and the energy functional, ultimately identifying it with $BU(\oo)$. We shall, on the other hand, march bravely in the opposite direction and analyze the \textit{bar construction} on the homotopy-coherent group $\ku \simeq \ZZ\times BU(\oo)$ identifying it with $U(\oo)$.
    Since the composite $\Omega B$ implements \textit{group completion}, we get
    \[
        B\ku \simeq 
        B M_\CC^\gp \simeq 
        B\Omega B M_\CC \simeq
        B M_\CC.
    \]
    Hence, it suffices to compute the bar construction on the homotopy-coherent commutative \textit{monoid}~$M_\CC$. That is, we need to show that
    \[
        B\big( \bigsqcup_{n\in \NN} BU(n)\big) \simeq U(\oo). 
    \]
    For this we introduce a strict topological model for the simplicial space $\BB M_\CC$. Let
    \[
        \mdef{X} \colon \Delta^\op \too \Top
    \]
    be the simplicial \textit{topological space} with
    \[
        X_n := \big\{(U_1,\dots,U_n) \mid U_i\le \CC^\oo, \ \dim(U_i)<\oo, \ U_i \perp U_j \text{ for } i\neq j\big\}.
    \]
    In words, $X_n$ consists of $n$-tuples of finite dimensional pairwise orthogonal subspaces of $\CC^\oo$. The degeneracy maps insert the zero subspace and the face maps replace two adjacent subspaces by their direct sum or discard the leading or trailing one. The contractibility of the infinite \textit{Stiefel spaces} $V_d(\CC^\oo)$ ensures that $X_\bullet$, naturally topologized, indeed models $\BB M_\CC$. Thus, the colimit of $\BB M_\CC$ in $\Spc$ is modeled by the \textit{topological geometric realization} $|X_\bullet| \in \Top$. A point in the latter is represented by a tuple of subspaces $U_1,\dots,U_n \le \CC^\oo$ as above together with a point
    \[
        0\le t_1\le \dots \le t_n \le 1 \qin |\Delta^n|.
    \]
    To such data one can associate a \textit{unitary operator} in $U(\oo)$ acting on $\CC^\infty$ with the scalar 
    \[
        \lambda_k:=\exp(2\pi i t_k) 
        \quad\text{on}\quad
        U_k
    \]
    and as the identity on their orthogonal complement. Conversely, the \textit{spectral theorem} expresses every element of $U(\oo)$ \textit{uniquely} in this form, after discarding the 1-eigenspace and merging summands with equal eigenvalues; these are precisely the simplicial gluing relations. We thus get a  \textit{bijective} map
    \[
        |X_\bullet| \iso U(\oo).
    \]
     This map is in fact a \textit{homeomorphism} and a fortiori a homotopy equivalence.\footnote{The point-set topological details we swept under the carpet can be found in \cite{harris1980bott}.} 
\end{proof}

As a corollary, we get an explicit description of the deloopings of $\ku$:
\[
    B^n \ku \simeq \begin{cases}
        \tau_{\ge n}(\ZZ \times BU(\oo)) & n \text{ even} \\
        \tau_{\ge n}U(\oo) & n \text{ odd}.
    \end{cases}
\]
Written like this, the $n$-connective cover seems utterly unnatural and practically begs to be removed.

\begin{defn}\label{Def_KU}
    The \tdef{nonconnective (complex) topological $K$-theory} spectrum is given by
    \[
        \mdef{\KU}_n := \begin{cases}
            \ZZ\times BU(\oo) & n \text{ even} \\
            U(\oo) & n \text{ odd}
        \end{cases}
    \]
    with equivalences $\KU_n \simeq \Omega \KU_{n+1}$ alternating between the canonical ones and Bott periodicity. 
\end{defn}

The spectrum $\KU$ is genuinely \textit{2-periodic}; it satisfies $\Omega^2 \KU \simeq \KU$ internally to the $\oo$-category $\Sp$. While it is, of course, possible to filter $\KU$ by bounded below spectra 
\[
    \KU \simeq \colim \tau_{\ge -n}\KU \qin \Sp,
\]
this periodicity would be lost. This is not merely an offense against aesthetics though. As we shall later see, the \textit{tensor product} of $\CC$-vector spaces endows $\ku$ and $\KU$ with a \textit{multiplicative structure}. On the level of~$\pi_*$, the canonical map $\ku \to \KU$, identifying $\ku$ with $\tau_{\ge 0}\KU$, induces a localization map of graded rings
\[
    \pi_*\ku = 
    \ZZ[\beta] \too 
    \ZZ[\beta^{\pm1}] = 
    \pi_*\KU 
    \quad\text{with}\quad
    |\mdef{\beta}| = 2.
\]
The generator $\beta \in \pi_2\ku$ is called the \tdef{Bott element}. The existence of an invertible element in nonzero degree forces the spectrum to be unbounded in both directions.
Thus, while the filtration by connective covers above recovers the underlying spectrum of $\KU$, it cannot capture its multiplicative structure. 

\begin{rem}
    There is a similar construction of the nonconnective \textit{real} topological $K$-theory spectrum $\mdef{\KO}$ using the \textit{$8$-fold} periodicity $\tau_{\ge 0}\Omega^8\ko \simeq \ko$. It also coincides with $\KU^{hC_2}$ for the action of $C_2$ on $\KU$ by complex conjugation. Furthermore, it exhibits $\KU$ as a \textit{Galois extension}\footnote{In the sense of \cite{RognesGal}.} of $\KO$, which is not true for the connective covers.     
\end{rem}

Such \textit{periodic ring spectra} occupy a central place in stable homotopy theory. For every prime $p$ and \tdef{height} $n\geq 1$, there is a \tdef{Morava $K$-theory} ring spectrum $\mdef{K(n)}$ having\footnote{$K(1)$ is just a slightly leaner version of $\KU/p$. As a rule, $p$ is traditionally omitted from the notation.}
\[
    \pi_*K(n)\simeq 
    \mathbb F_p[v_n^{\pm1}]
    \quad\text{with}\quad |\mdef{v}_n|=2(p^n-1).
\]
Together with 
\[
    K(\mdef{0}):=\QQ
    \qquad\text{and}\qquad
    K(\mdef{\oo}) := \FF_p,
\]
they serve, in a precise sense, as the \textit{prime fields} of $\Sp$. Accordingly, stable homotopy theory can be analyzed prime by prime with respect to these ``primes'', in close analogy with localization (and completion) in ordinary commutative algebra.
This is known as the \textit{chromatic picture}. Much as a prism resolves white light into its constituent wavelengths, the Morava $K$-theories resolve stable homotopy into periodic phenomena of different heights.


\section{Commutative rings}

Having treated homotopy-coherent (commutative) monoids and groups, the next natural species to consider is that of (commutative) \textit{rings}. The reader should be warned that algebraic theories are not quite the default framework for doing so, and that the resulting theory has certain limitations. 
For example, it produces (commutative) connective ring spectra such as $\ku$, but does not readily capture \textit{nonconnective} ones such as $\KU$.
The more robust approach is to consider (commutative) \textit{algebras} in spectra for the \textit{tensor-product} symmetric monoidal structure, as we shall do in the next chapter. 
Nevertheless, algebraic theories have something interesting to say here, and the detour is worth the reader's time.\footnote{The definitive source for this material is \cite{cnossen2025homotopical} and the references cited there.}

\subsection{Bispans}
We focus on the commutative case and, as we did with groups, do not yet impose additive inverses. In other words, we shall consider \textit{commutative semirings}. 
The ordinary algebraic theory of commutative semirings consists of polynomial semirings over $\NN$ in finitely many variables where a morphism  
\[
    \NN[y_1,\dots,y_m] \too \NN[x_1,\dots,x_n]
\]
is given by a sequence of polynomials
\[
    f_1,\dots,f_m \qin \NN[x_1,\dots,x_n].
\]
Each such polynomial is a \textit{sum} of \textit{products} of the variables $x_1,\dots,x_n$. For the homotopy-coherent version, one needs to keep track of the \textit{symmetries} of the sets on which these sums and products are indexed. Moreover, in addition to the individual commutativity and associativity of the sum and product, one has to build in homotopy-coherent \textit{distributivity}. 

There is a \textit{span-like} technology that manages to achieve this rather neatly. And as in the case of commutative monoids, it is clearer when one replaces the numbered sequences of variables $x_1,\dots,x_n$ and $y_1,\dots,y_m$  with general finite sets $X$ and $Y$. 
\begin{defn}
    A \tdef{bispan} from $X$ to $Y$ is a diagram in $\Fin$ of the form
    \[\begin{tikzcd}
    	& V & M & \\
    	X &&& Y.
    	\arrow[from=1-2, to=1-3]
    	\arrow[from=1-2, to=2-1]
    	\arrow[from=1-3, to=2-4]
    \end{tikzcd}\]    
\end{defn}

It should be thought of as a commutative semiring homomorphism, from the free such on $Y$ variables to the free such on $X$ variables, given by the polynomials (with hopefully self-explanatory notation):
\[
    y\longmapsto \sum_{m\in M_y} \prod_{v\in V_m} x(v).
\]
The letter `$M$' stands for the \textit{monomials} to be summed and `$V$' for the \textit{variables} to be multiplied in each monomial. The symmetries of such expressions are encoded by \textit{isomorphisms} of bispans, which are commutative diagrams of the form
\[\begin{tikzcd}[row sep = 1.2em]
	& V & M & \\
	X &&& Y \\
	& {V'} & {M'}
	\arrow[from=1-2, to=1-3]
	\arrow[from=1-2, to=2-1]
	\arrow["\wr"', from=1-2, to=3-2]
	\arrow[from=1-3, to=2-4]
	\arrow["\wr"', from=1-3, to=3-3]
	\arrow[from=3-2, to=2-1]
	\arrow[from=3-2, to=3-3]
	\arrow[from=3-3, to=2-4]
\end{tikzcd}\]

\begin{example}
    For $X = Y = \pt$, a (flattened) bispan diagram in $\Fin$ 
    \[
        \pt \froom V \too M \too \pt
    \]
    represents a single polynomial in a single variable with \textit{finite-set exponents}
    \[
        y\longmapsto \sum_{m \in M} x^{V_m}.
    \]
    An automorphism of such a polynomial is a permutation of the summands $\sigma\colon M\iso M$ together with identifications of the corresponding exponent sets $V_m\iso V_{\sigma(m)}$ for each $m\in M$.
\end{example}

The real question is of course how to define the composition of two such bispans,
\[\begin{tikzcd}
	& V & M && U & N & \\
	X &&& Y &&& Z
	\arrow["\Pi", color={rgb,255:red,179;green,179;blue,179}, from=1-2, to=1-3]
	\arrow[from=1-2, to=1-3]
	\arrow["x"', color={rgb,255:red,179;green,179;blue,179}, from=1-2, to=2-1]
	\arrow[from=1-2, to=2-1]
	\arrow["\Sigma", color={rgb,255:red,179;green,179;blue,179}, from=1-3, to=2-4]
	\arrow[from=1-3, to=2-4]
	\arrow["\Pi", color={rgb,255:red,179;green,179;blue,179}, from=1-5, to=1-6]
	\arrow[from=1-5, to=1-6]
	\arrow["y"', color={rgb,255:red,179;green,179;blue,179}, from=1-5, to=2-4]
	\arrow[from=1-5, to=2-4]
	\arrow["\Sigma", color={rgb,255:red,179;green,179;blue,179}, from=1-6, to=2-7]
	\arrow[from=1-6, to=2-7]
\end{tikzcd}\]
which represents the composite operation
\[
    \sum_{n\in N_z} \prod_{u\in U_{n}}
    \sum_{m\in M_{y(u)}} \prod_{v\in V_m} x(v).
\]
The culprit is the middle $\Pi\Sigma$, a \textit{product of sums}, which needs to be rewritten as a \textit{sum of products}. Isolating the heart of the matter, suppose we are given a pair of maps of finite sets
\[
    X \oto{\quad s \quad} A \oto{\quad p \quad} Y
\]
representing the operation
\[
    y\mapsto 
    \prod_{a\in A_y}\Big(\sum_{x\in X_a} x\Big).
\]
After expanding, the monomials are indexed by all the possible ways to choose an $x \in X_a$ for every $a \in A_y$. In other words, they are indexed by the set of \textit{sections} of the map $s\colon X \to A$ over the subset $A_y \sseq A$, 
\[
    \Gamma_y(s) := 
    \{f\colon A_y \to X \mid s\circ f = \Id\}.
\]
This yields the following \textit{distributivity equation}:
\[
    \prod_{a\in A_y}\sum_{x\in X_a} x \simeq
    \sum_{f\in \Gamma_y(s)} \prod_{a \in A_y} x(f(a)).
\]
The disjoint union over $Y$ of all these fiberwise sections 
\[
    \pi \colon \coprod_{y\in Y} \Gamma_y(s) \too Y
\]
is obtained by applying to $s\colon X\to A$, viewed as an object of $\Fin_{/A}$, the right adjoint $p_*$ of the basechange functor $p^*$, 
\[
    p_* \colon 
    \Fin_{/A} \adj 
    \Fin_{/Y}
    \noloc p^*.
\]
Diagrammatically, the right-hand side of the distributivity equation is thus encoded by the (top) bispan 
\[\begin{tikzcd}[sep = 3em]
	& \textcolor{rgb,255:red,31;green,122;blue,31}{{p^*p_*X}} & \textcolor{rgb,255:red,31;green,122;blue,31}{{p_*X}} & \\
	\textcolor{rgb,255:red,31;green,122;blue,31}{X} & A & Y & \textcolor{rgb,255:red,31;green,122;blue,31}{Y}
	\arrow[draw={rgb,255:red,31;green,122;blue,31}, from=1-2, to=1-3]
	\arrow[draw={rgb,255:red,31;green,122;blue,31}, from=1-2, to=2-1]
	\arrow[from=1-2, to=2-2]
	\arrow["\lrcorner"{anchor=center, pos=0.125}, draw=none, from=1-2, to=2-3]
	\arrow["\pi"', from=1-3, to=2-3]
	\arrow[draw={rgb,255:red,31;green,122;blue,31}, from=1-3, to=2-4]
	\arrow["s", from=2-1, to=2-2]
	\arrow["p", from=2-2, to=2-3]
	\arrow[equals, from=2-3, to=2-4]
\end{tikzcd}\]


where the left diagonal is the adjunction counit.

\begin{example}
    The case $Y=\pt$ is particularly instructive.
    Thinking of $s\colon X\to A$ as an $A$-indexed family
    $\{X_a\}_{a\in A}$ of finite sets, a monomial in the
    expansion of the ``product of sums'' is specified by
    a choice of an element $x_a\in X_a$ for each $a\in A$.
    These choices form precisely the set
    $p_*X=\prod_{a\in A}X_a$.
    Thus, the bispan
    \[
        X \froom
        p^*p_*X \too
        p_*X \too
        \pt
    \]
    takes the concrete form
    \[
        X \froom
        A \times \prod_{a\in A} X_a
        \oto{\Gray{\ \Pi\ }}
        \prod_{a\in A} X_a
        \oto{\Gray{\ \Sigma\ }}
        \pt.
    \]
    The first map is the adjunction counit, given by evaluation at the \(a\)-coordinate.
\end{example}

Putting everything back together, the composition of bispans is described by the rather baroque diagram:
\[\begin{tikzcd}[column sep = 3.2em]
	&&& \textcolor{rgb,255:red,31;green,122;blue,31}{W} & {p^*p_*Q} & \textcolor{rgb,255:red,31;green,122;blue,31}{{p_*Q}} & \\
	&&& Q \\
	& \textcolor{rgb,255:red,31;green,122;blue,31}{V} & \textcolor{rgb,255:red,31;green,122;blue,31}{\ M\ } && \textcolor{rgb,255:red,31;green,122;blue,31}{U} & \textcolor{rgb,255:red,31;green,122;blue,31}{N} \\
	\textcolor{rgb,255:red,31;green,122;blue,31}{X} &&& \textcolor{rgb,255:red,31;green,122;blue,31}{Y} &&& \textcolor{rgb,255:red,31;green,122;blue,31}{Z}
	\arrow[from=1-4, to=1-5]
	\arrow[draw={rgb,255:red,31;green,122;blue,31}, curve={height=-18pt}, dashed, from=1-4, to=1-6]
	\arrow[from=1-4, to=3-2]
	\arrow["\lrcorner"{anchor=center, pos=0.125}, draw=none, from=1-4, to=3-5]
	\arrow[draw={rgb,255:red,31;green,122;blue,31}, curve={height=24pt}, dashed, from=1-4, to=4-1]
	\arrow[from=1-5, to=1-6]
	\arrow[from=1-5, to=2-4]
	\arrow[from=1-5, to=3-5]
	\arrow["\lrcorner"{anchor=center, pos=0.125}, draw=none, from=1-5, to=3-6]
	\arrow["\pi"', from=1-6, to=3-6]
	\arrow[draw={rgb,255:red,31;green,122;blue,31}, curve={height=-12pt}, dashed, from=1-6, to=4-7]
	\arrow[from=2-4, to=3-3]
	\arrow[from=2-4, to=3-5]
	\arrow["\lrcorner"{anchor=center, pos=0.125, rotate=-45}, draw=none, from=2-4, to=4-4]
	\arrow[color={rgb,255:red,31;green,122;blue,31}, from=3-2, to=3-3]
	\arrow[color={rgb,255:red,31;green,122;blue,31}, from=3-2, to=4-1]
	\arrow[color={rgb,255:red,31;green,122;blue,31}, from=3-3, to=4-4]
	\arrow["p", color={rgb,255:red,31;green,122;blue,31}, from=3-5, to=3-6]
	\arrow[color={rgb,255:red,31;green,122;blue,31}, from=3-5, to=4-4]
	\arrow[color={rgb,255:red,31;green,122;blue,31}, from=3-6, to=4-7]
\end{tikzcd}\]

Encoding polynomial composition --- and especially the distributivity hidden within it --- by diagrammatic maneuvers is already something of an achievement. But a recipe for composing \textit{two} bispans is not yet a fully coherent associative composition law. It still remains to construct the \textit{$\oo$-category} of bispans. 

For that we observe that a bispan can be thought of as a (special kind of a) \textit{span of spans},
\[\begin{tikzcd}[row sep = 1em]
	&& M && \\
	& V && M \\
	X &&&& Y.
	\arrow[from=2-2, to=1-3]
	\arrow[from=2-2, to=3-1]
	\arrow[equals, from=2-4, to=1-3]
	\arrow[from=2-4, to=3-5]
\end{tikzcd}\]
This observation might seem accidental were it not for the following miracle: composition of bispans in $\Fin$ is precisely what composition in $\Span(\Span(\Fin))$ would be. Strictly speaking, the latter category does not quite exist, since composing arbitrary spans of spans requires pullbacks that need not exist in $\Span(\Fin)$. For those arising from bispans, however, the usual recipe is well-defined and recovers exactly the composition constructed by hand above.\footnote{A proof can be found in \cite[Proposition 2.5.9]{elmanto2023distributivity}.}
\newcommand{\Bispan}{\mathrm{Bispan}}
\newcommand{\fwd}{\mathrm{fwd}}
\begin{defn}
    The $(2,1)$-category of bispans of finite sets is given by
    \[
        \mdef{\Bispan}(\Fin) := 
        \Span(\Span(\Fin),\all,\fwd),
    \]
    where `all' means backward legs are allowed to be arbitrary maps in $\Span(\Fin)$ while `fwd' means the forward legs are allowed to be only the \textit{forward maps} in $\Span(\Fin)$,
    \[
        X \ofrom{\quad =\quad } 
        X \oto{\quad \ \quad}
        Y.
    \]
\end{defn}

There is nothing special about $\Fin$ in the above definition. 
The general construction of \textit{span} $\oo$-categories (allowing restrictions on backward and forward legs) readily produces the $\oo$-category of \textit{bispans} by suitable iteration as long as the original $\oo$-category admits nicely behaved pullbacks.\footnote{The basechange functor needs to admit a \textit{right} adjoint.} 

\subsection{Commutative (semi)rings}

It can be easily verified that finite \textit{products} in $\Bispan(\Fin)$ exist and are computed by levelwise \textit{coproduct} in $\Fin$. In particular, every finite set is the product of its points in $\Bispan(\Fin)$. We think of $\Bispan(\Fin)^\op$ as the algebraic theory of homotopy-coherent commutative semirings. 

\begin{defn}
    The $\oo$-category of \tdef{commutative semirings} in an $\oo$-category $\cC$ is given by\footnote{We use `rig' as another term for \textit{semiring}, which is a ri\textbf{n}g with no \textbf{n}egatives.}
    \[
        \mdef{\CRig}(\cC) := \Fun^\times(\Bispan(\Fin), \cC).
    \]
\end{defn}

The underlying space of a commutative semiring admits \textit{two} commutative monoid structures: the additive and the multiplicative. These correspond to two product-preserving inclusions
\[
    \Span(\Fin) \too \Bispan(\Fin)
\]
taking a span 
\[
    X\froom A \too Y \qin \Fin
\]
to one of the following bispans
\[\begin{tikzcd}
	& A & A &&&& A & Y & \\
	X &&& Y & {,} & X &&& Y.
	\arrow[equals, from=1-2, to=1-3]
	\arrow[from=1-2, to=2-1]
	\arrow[from=1-3, to=2-4]
	\arrow["\Sigma", color={rgb,255:red,179;green,179;blue,179}, draw=none, from=1-3, to=2-4]
	\arrow[from=1-7, to=1-8]
	\arrow["\Pi", color={rgb,255:red,179;green,179;blue,179}, draw=none, from=1-7, to=1-8]
	\arrow[from=1-7, to=2-6]
	\arrow[equals, from=1-8, to=2-9]
\end{tikzcd}\]
The reader is invited to verify that these embeddings respect composition. Restriction along these functors corresponds to the two forgetful functors 
\[
    \CRig(\cC) \too \CMon(\cC)
\]
remembering either only the \textit{additive} commutative monoid structure or only the \textit{multiplicative} one. This allows one to define commutative \textit{rings}. 

\begin{defn}
    The $\oo$-category of \tdef{commutative rings} in an $\oo$-category $\cC$ is given by the full subcategory
    \[
        \mdef{\CRing}(\cC) \ \sseq \ \CRig(\cC)
    \]
    on those semirings whose underlying additive monoid is a \textit{group}. Objects of $\CRing(\Spc)$ will be referred to as (connective) \tdef{commutative ring spectra}.\footnote{The equivalence of this definition with the standard one is \cite[Corollary B]{cnossen2025homotopical}.
    \textit{Nonconnective} ring spectra will have to wait until the next chapter.} 
\end{defn}

For an ordinary category $\cC$, all these notions recover the usual ones. In particular, $\CRig(\Set)$ and $\CRing(\Set)$ are the ordinary categories $\CRig$ and $\CRing$ of commutative (semi)rings. This follows from the equivalence of categories
\[
    h\Bispan(\Fin) \simeq \CRig^{\ff,\op}.
\]
Moreover, the truncation functor 
\[
    \Bispan(\Fin)\too h\Bispan(\Fin)
\]
is product-preserving and thus induces a functor on \textit{models}, which preserves the grouplike  objects
\[\begin{tikzcd}
	\textcolor{rgb,255:red,179;green,179;blue,179}{{\Fun^\times(h\Bispan(\Fin),\Spc)}} & \textcolor{rgb,255:red,179;green,179;blue,179}{{\Fun^\times(\Bispan(\Fin),\Spc)}} \\
	{\Anim(\CRing)} & {\CRing(\Spc)}
	\arrow[color={rgb,255:red,179;green,179;blue,179}, from=1-1, to=1-2]
	\arrow[color={rgb,255:red,179;green,179;blue,179}, hook, from=2-1, to=1-1]
	\arrow[from=2-1, to=2-2]
	\arrow[color={rgb,255:red,179;green,179;blue,179}, hook, from=2-2, to=1-2]
\end{tikzcd}\]

This functor is \textit{not} an equivalence. 
\begin{rem}
    This situation reminds one of the distinction between 
    \[
        \Anim(\Ab) \simeq \sD_{\ge 0}(\ZZ)
        \qquad\text{and}\qquad
        \CGrp(\Spc)\simeq \Sp_{\ge 0},
    \]
    where the former is equivalent to the $\oo$-category of (homotopy-coherent) \textit{modules} over the (discrete) commutative ring $\ZZ$ in $\Sp_{\ge 0}$. 
    In fact, it suggests a \textit{third} candidate for what homotopy-coherent commutative rings might be: the homotopy-coherent \tdef{commutative $\ZZ$-algebras}. Namely, the under category 
    \[
        \CRing_{\mdef{\ZZ}}(\Spc) := \CRing(\Spc)_{\ZZ/}.
    \]
    This one is equivalent to \textit{commutative algebras} in $\sD_{\ge0}(\ZZ)$ with respect to the derived tensor product~$\otimes_{\ZZ}^{\mathbb{L}}$.  
    The difference between the three options can be loosely attributed to the level of ``strictness'' in the commutativity of addition and multiplication.
    One is advised to \textit{unask} which is the correct homotopy-coherent version of commutative rings. 
\end{rem}

\subsection{Examples}

The initial (ordinary) commutative semiring is $\NN$ with its usual addition and multiplication. In the homotopy-coherent world we have the following matural analogue:
\begin{prop}
    The \textit{initial object} of $\CRig(\Spc)$ has underlying space $\MM$. The addition and multiplication maps 
    \[
        +\colon \MM \times \MM \too \MM \qquad,\qquad
        \cdot\ \colon \MM \times \MM \too \MM
    \]
    are given respectively by the disjoint union $\sqcup$ and cartesian product $\times$ of finite sets. 
\end{prop}

We have secretly used this semiring structure of $\MM$ when defining the category $\Mat(\MM)$ before recasting it in terms of $\Span(\Fin)$.

\begin{proof}
    The initial object in the category
    \[
        \Bispan(\Fin)^\op \simeq 
        \CRig(\Spc)^{\ff} \sseq 
        \CRig(\Spc)
    \]
    is represented by $\es$ --- the free commutative semiring on the empty set of variables. To compute its \textit{underlying space}, we need to compute the space of maps from the free commutative semiring on one element. Namely, 
    \[
        U_{\CRig}F_{\CRig}(\es) \simeq 
        \Map_{\Bispan(\Fin)}(\es,\pt) \qin \Spc.
    \]
    This is the groupoid of bispans of the form
    \[
        \es \ofrom{\quad = \quad} 
        \es \oto{\quad \ \quad}
        X \oto{\quad \ \quad}
        \pt,
    \]
    which is exactly the groupoid $\MM$ of finite sets. Given two such elements $X, Y\in \Fin$, their sum is computed by the composition
    \[\begin{tikzcd}[column sep = 2em]
    	&&& \textcolor{rgb,255:red,31;green,122;blue,31}{\es} & {X\sqcup Y} & \textcolor{rgb,255:red,31;green,31;blue,122}{{X\sqcup Y}} & \\
    	&&& {X\sqcup Y} \\
    	& \textcolor{rgb,255:red,31;green,122;blue,31}{\es} & \textcolor{rgb,255:red,31;green,122;blue,31}{{X\sqcup Y}} && \textcolor{rgb,255:red,31;green,122;blue,31}{{\pt \sqcup\pt}} & \textcolor{rgb,255:red,31;green,122;blue,31}{{\pt \sqcup \pt}} \\
    	\textcolor{rgb,255:red,31;green,122;blue,31}{\es} &&& \textcolor{rgb,255:red,31;green,122;blue,31}{{\pt \sqcup\pt}} &&& \textcolor{rgb,255:red,31;green,122;blue,31}{\pt}
    	\arrow[from=1-4, to=1-5]
    	\arrow[draw={rgb,255:red,31;green,122;blue,31}, curve={height=-18pt}, dashed, from=1-4, to=1-6]
    	\arrow[from=1-4, to=3-2]
    	\arrow["\lrcorner"{anchor=center, pos=0.125}, draw=none, from=1-4, to=3-5]
    	\arrow[draw={rgb,255:red,31;green,122;blue,31}, curve={height=24pt}, dashed, from=1-4, to=4-1]
    	\arrow[equals, from=1-5, to=1-6]
    	\arrow[equals, from=1-5, to=2-4]
    	\arrow[from=1-5, to=3-5]
    	\arrow["\lrcorner"{anchor=center, pos=0.125}, draw=none, from=1-5, to=3-6]
    	\arrow[from=1-6, to=3-6]
    	\arrow[draw={rgb,255:red,31;green,122;blue,31}, curve={height=-12pt}, dashed, from=1-6, to=4-7]
    	\arrow[from=2-4, to=3-3]
    	\arrow[from=2-4, to=3-5]
    	\arrow["\lrcorner"{anchor=center, pos=0.125, rotate=-45}, draw=none, from=2-4, to=4-4]
    	\arrow[draw={rgb,255:red,31;green,122;blue,31}, from=3-2, to=3-3]
    	\arrow[draw={rgb,255:red,31;green,122;blue,31}, from=3-2, to=4-1]
    	\arrow[draw={rgb,255:red,31;green,122;blue,31}, from=3-3, to=4-4]
    	\arrow[draw={rgb,255:red,31;green,122;blue,31}, equals,  from=3-5, to=3-6]
    	\arrow[draw={rgb,255:red,31;green,122;blue,31}, equals,  from=3-5, to=4-4]
    	\arrow[draw={rgb,255:red,31;green,122;blue,31}, from=3-6, to=4-7]
    \end{tikzcd}\]
    whereas their product is computed by the composition
    \[\begin{tikzcd}[column sep = 2em]
    	&&& \textcolor{rgb,255:red,31;green,122;blue,31}{\es} & {(X\times Y)^{\sqcup2}} & \textcolor{rgb,255:red,31;green,31;blue,122}{{X\times Y}} & \\
    	&&& {\quad X\sqcup Y \quad } \\
    	& \textcolor{rgb,255:red,31;green,122;blue,31}{\es} & \textcolor{rgb,255:red,31;green,122;blue,31}{{X\sqcup Y}} && \textcolor{rgb,255:red,31;green,122;blue,31}{{\pt \sqcup\pt}} & \textcolor{rgb,255:red,31;green,122;blue,31}{\pt} \\
    	\textcolor{rgb,255:red,31;green,122;blue,31}{\es} &&& \textcolor{rgb,255:red,31;green,122;blue,31}{{\pt \sqcup\pt}} &&& \textcolor{rgb,255:red,31;green,122;blue,31}{\pt}
    	\arrow[from=1-4, to=1-5]
    	\arrow[draw={rgb,255:red,31;green,122;blue,31}, curve={height=-18pt}, dashed, from=1-4, to=1-6]
    	\arrow[from=1-4, to=3-2]
    	\arrow["\lrcorner"{anchor=center, pos=0.125}, draw=none, from=1-4, to=3-5]
    	\arrow[draw={rgb,255:red,31;green,122;blue,31}, curve={height=24pt}, dashed, from=1-4, to=4-1]
    	\arrow[from=1-5, to=1-6]
    	\arrow[from=1-5, to=2-4]
    	\arrow[from=1-5, to=3-5]
    	\arrow["\lrcorner"{anchor=center, pos=0.125}, draw=none, from=1-5, to=3-6]
    	\arrow[from=1-6, to=3-6]
    	\arrow[draw={rgb,255:red,31;green,122;blue,31}, curve={height=-12pt}, dashed, from=1-6, to=4-7]
    	\arrow[equals, from=2-4, to=3-3]
    	\arrow[from=2-4, to=3-5]
    	\arrow["\lrcorner"{anchor=center, pos=0.125, rotate=-45}, draw=none, from=2-4, to=4-4]
    	\arrow[draw={rgb,255:red,31;green,122;blue,31}, from=3-2, to=3-3]
    	\arrow[draw={rgb,255:red,31;green,122;blue,31}, from=3-2, to=4-1]
    	\arrow[draw={rgb,255:red,31;green,122;blue,31}, from=3-3, to=4-4]
    	\arrow[color={rgb,255:red,31;green,122;blue,31}, from=3-5, to=3-6]
    	\arrow[draw={rgb,255:red,31;green,122;blue,31}, equals, from=3-5, to=4-4]
    	\arrow[draw={rgb,255:red,31;green,122;blue,31}, equals, from=3-6, to=4-7]
    \end{tikzcd}\]
\end{proof}

\begin{rem}
    The two bispans used in the proof to compute the addition and multiplication in the commutative semiring $\MM$, 
    \[\begin{tikzcd}[column sep=1.5em]
    	& {\pt\sqcup\pt} && {\pt\sqcup\pt} &&&& {\pt\sqcup\pt} && \pt & \\
    	{\pt\sqcup\pt} &&&& \pt & {,} & {\pt\sqcup\pt} &&&& \pt
    	\arrow[equals, from=1-2, to=1-4]
    	\arrow[equals, from=1-2, to=2-1]
    	\arrow[from=1-4, to=2-5]
    	\arrow["\Sigma", color={rgb,255:red,179;green,179;blue,179}, draw=none, from=1-4, to=2-5]
    	\arrow[from=1-8, to=1-10]
    	\arrow["\Pi", color={rgb,255:red,179;green,179;blue,179}, from=1-8, to=1-10]
    	\arrow[equals, from=1-8, to=2-7]
    	\arrow[equals, from=1-10, to=2-11]
    \end{tikzcd}\]
    correspond under the Yoneda embedding
    \[
        \Bispan(\Fin)^\op \intoo \CRig(\Spc)
    \]
    to the two commutative semiring homomorphisms $\MM[t] \to \MM[x,y]$ given by $t \mapsto x+y$ and $t\mapsto xy$ respectively.
\end{rem}

More generally, one can consider commutative semirings in $\Cat_\oo$. Namely, $\oo$-categories $\cC$ with two symmetric monoidal structures satisfying coherent distributivity. If the additive one is \textit{cocartesian}, i.e., given by the coproduct $\sqcup$, and the multiplicative one is denoted by $\otimes$ there is a canonical and natural distributivity map
\[
    \bigsqcup_{a\in A}(X\otimes Y_a)\too
    X\otimes\big(\bigsqcup_{a\in A}Y_a\big).
\] 
If this map happens to be always an isomorphism, $\cC$ promotes canonically to a commutative semiring in $\Cat_\oo$. In particular, this works when $\cC$ is cartesian closed and the multiplicative structure is \textit{cartesian}, such as $\cC = \Fin$. Passing to the underlying space we obtain a homotopy-coherent commutative semiring whose entire distributivity coherence comes for free.\footnote{See \cite{elmanto2023distributivity} for this and more.}

\begin{example}
    For a \textit{commutative} ring $R$, the groupoid $M_R$ of finite free $R$-modules from \Cref{Ex_Proj_CMon} is a commutative semiring with addition given by $\oplus$ and multiplication by the relative tensor product $\otimes_R$. 
\end{example}

The \textit{group completion} of $M_R$ is the algebraic $K$-theory (connective) spectrum $K(R)$ (see \Cref{Ex_Alg_K_Thy}). It would be great if for \textit{commutative} $R$, we could transfer the multiplication $\otimes_R$ of $M_R$ through the group completion making $K(R)$ a (connective) commutative ring spectrum. By the adjoint functor theorem, it is formal that the inclusion
\[
    \CRing(\Spc) \intoo \CRig(\Spc)
\]
admits a \textit{left adjoint}. What is not formal though is that this left adjoint is given by group completing the underlying additive commutative monoid. Although the lack of additive inverses is precisely the defect this left adjoint is meant to repair, one could imagine the multiplicative structure interfering with the process rather than being silently complicit in it. This complicity is guaranteed in ordinary algebra by \textit{distributivity}, and the same turns out to hold homotopy-coherently. But the proof becomes much more transparent in the language of symmetric monoidal categories, which is why we defer it to the next chapter. 

\begin{prop}
    The following square, with vertical maps forgetting multiplication, is horizontally left adjointable: 
    \[\begin{tikzcd}[sep = large]
    	{\CRing(\Spc)} & {\CRig(\Spc)} \\
    	{\CGrp(\Spc)} & {\CMon(\Spc)}
    	\arrow[hook, from=1-1, to=1-2]
    	\arrow[from=1-1, to=2-1]
    	\arrow[curve={height=10pt}, dashed, from=1-2, to=1-1]
    	\arrow[from=1-2, to=2-2]
    	\arrow[hook, from=2-1, to=2-2]
    	\arrow[curve={height=10pt}, dashed, from=2-2, to=2-1]
    \end{tikzcd}\]
\end{prop}

Aside from $K(R)$ for a commutative ring $R$, we shall consider one additional example of a commutative ring spectrum.  

\begin{cor}
    The initial object of $\CRing(\Spc)$ has underlying connective spectrum $\Sph$. 
\end{cor}
\begin{proof}
    The left adjoint 
    \[
        \CRig(\Spc) \too \CRing(\Spc)
    \]
    preserves initial objects. In the source the initial object is $\MM$ and hence the underlying additive group of the initial object in the target is $\MM^\gp \simeq \Sph$.
\end{proof}

The \textit{multiplicative} commutative monoid of $\Sph$, and especially its subgroup of \textit{units} $\Sph^\times$, are very interesting objects as well.\footnote{As a connective spectrum, $\Sph^\times$ is classically denoted by $\gl_1(\Sph)$.} 
The action of $O(\oo)$ on the sphere spectrum viewed as an infinite loop space $\Sph = \colim \Omega^nS^n$ arises via a homotopy-coherent \textit{commutative} group homomorphism known as the \tdef{$J$-homomorphism},
\[
    \mdef{J} \colon O(\infty) \too \Sph^\times \qin \CGrp(\Spc).
\]
To recover its usual formulation, observe that $\Sph^\times \sseq \Sph$ is the union of the components indexed by $\pm 1$. Its identity component is therefore $\Sph^\times_1 = \Sph^{}_1$. 
Translation by $-1$ in the additive group $\Sph$ gives an equivalence of pointed spaces 
\[
    \Sph^\times_1 \iso \Sph_0^{} \qquad\text{by}\qquad x \longmapsto x-1.
\]
Hence, for $n \ge 1$, the map $J$ induces group homomorphisms
\[
    \pi_n O(\oo) \too \pi_n(\Sph^\times) \simeq \pi_n(\Sph).
\]
The equivalence $\Sph^\times_1 \simeq \Sph^{}_0$, however, identifies only the underlying \textit{pointed spaces}, not their \textit{group structures}. 
Thus, although the classical $J$-homomorphism is expressed in terms of the homotopy groups of $\Sph$, its actual target before taking homotopy groups is the multiplicative group $\Sph^\times$ and not the additive group~$\Sph$.

\begin{rem}
    A better formulation still is in terms of the \textit{deloopings} of the source and target. There is a functor 
    \[
        M_\RR := 
        \bigsqcup_{n\in \NN} BO(n) \oto{\ (-)^+\ }
        \Spc_* \oto{\ \Sigma^\oo\ }
        \Sp
    \]
    taking a finite dimensional $\RR$-vector space $V$ to the suspension spectrum of the one-point compactification $V^+$, which we denote by $\Sph^V$. This construction takes a direct sum of vector spaces to the \textit{tensor product} of spectra
    \[
        \Sph^{V\oplus W} \simeq 
        \Sph^V \otimes \Sph^W \qin \Sp.
    \]
    Moreover, the map lands in the subgroup $\mdef{\pic}(\Sp)$ of \textit{tensor-invertible} spectra, and thus factors through the \textit{group completion} of the source
    \[
        \ZZ \times BO(\oo) \simeq \ko \too \pic(\Sp) \qin \Sp_{\ge 0}.
    \]
    Looping this map once yields the $J$-homomorphism above. 
\end{rem}
 
\section{Equivariance}

There is a standard procedure for producing further mathematics: prepend the word \textit{equivariant}. Any mathematical object worth studying eventually demands an analysis of its symmetries. Yet, the real complexity lies in the attempt to axiomatize the structure of  \textit{fixed-point data} for a group action, and especially its presumed interaction with addition and multiplication operations. Our discussion of this theme will only scratch the surface. It will culminate with the definitions of the very basic objects of study in equivariant stable homotopy theory, but will stop short of defining the all-important \textit{geometric fixed points}.\footnote{Which, in the finest tradition of mathematical terminology, are neither geometric nor fixed points.}

\subsection{Naive $G$-actions}

An action of a finite group $G$ is a rather simple kind of algebraic structure governed by the algebraic theory~$\Fin_G^\ff$ of \textit{finite free} $G$-sets. 
\begin{prop}\label{Naive_Action}
    Let $\cC$ be an $\oo$-category and $G$ a finite group. There is an equivalence of $\oo$-categories,
    \[
        \Fun^\times(\Fin_G^{\ff,\op} , \cC) \simeq 
        \Fun(BG, \cC).
    \]
\end{prop}
\begin{proof}
    This is due to $\Fin_G^\ff$ being itself freely generated under \textit{finite coproducts} by the single free $G$-orbit $G/e$ whose automorphism group is $G$.\footnote{In fact, it is canonically $G^\op$, but these are isomorphic.} 
\end{proof}

This is the notion of an object with a (naive) \tdef{$G$-action}.
While not uninteresting, from an abstract point of view there is little more to say. 

\begin{example}
    The $\oo$-category $\Cat_\oo$ of $\oo$-categories has $C_2$ as its automorphism group. The nontrivial element acts as the \textit{covariant} involution\footnote{See, for example, \cite[\S 4.4]{lurie2009infinity}.}  
    \[
        (-)^\op\colon \Cat_\oo \iso \Cat_\oo.
    \]
    An object of $\Cat_\oo^{hC_2}$ is an $\oo$-category with \tdef{duality}. This is roughly an $\oo$-category $\cC$ with an equivalence 
    \[
        D \colon \cC \iso \cC^\op \qin \Cat_\oo,
    \]
    such that $D\circ D^\op \simeq \Id$, but also with all the higher homotopies exhibiting $D$ as a coherent \textit{anti-involution}.\footnote{Even for an \textit{ordinary} category $\cC$, there is one further condition.} One natural example of this is \textit{symmetric monoidal duality}, such as the dual of a finite dimensional vector space. When $\cC$ is in addition \textit{stable} and $D$ is \textit{exact}, this is one possible (though not the most general) input for \textit{Hermitian $K$-theory} and \textit{$L$-theory}.\footnote{A detailed account of this can be found in \cite{heine2016infinity}.} 
\end{example}

\subsection{Declared fixed points}

To make things more interesting, one could use the \textit{entire} category $\Fin_G$ as an algebraic theory. This one contains $G$-orbits of all types. In fact, it is similarly freely generated under finite coproducts from the full subcategory $\cO_G \sseq \Fin_G$ on \textit{orbits} $G/H$. 

\begin{defn}
    Let $\cC$ be an $\oo$-category and $G$ a finite group. The $\oo$-category of \tdef{$G$-equivariant} objects in $\cC$ is given by
    \[
        \cC_{\mdef{G}} :=
        \Fun(\cO_G^\op, \cC) 
        \Gray{\simeq 
        \Fun^\times(\Fin_G^\op, \cC).}
    \]
\end{defn}

There is a natural passage in both directions between $G$-equivariant objects and objects with a naive $G$-action. The full subcategory on the free orbit $G/e \in \cO_G^\op$ is equivalent to $BG$. Restriction and (the fully faithful) right Kan extension along its inclusion constitute an adjunction
\[\begin{tikzcd}[column sep = 2em]
	{\Fun(\cO_G^\op,\cC)} & {\Fun(BG,\cC).}
	\arrow[shift right, from=1-1, to=1-2]
	\arrow[shift right, hook', from=1-2, to=1-1]
\end{tikzcd}\]
One may thus say that a $G$-equivariant object $X$ has an underlying object $X(G/e)$ with a naive $G$-action. And conversely, an object with a $G$-action is equivalently a \textit{naive} $G$-equivariant object,\footnote{Another classical term for that is \textit{Borel}.} 
\[
    X \colon \cO_G^\op \too \cC
    \qquad\text{with}\qquad
    X(G/H) := X^{hH} \qin \cC.
\]
Modeled on this, we regard $X(G/H)$, for an arbitrary $G$-equivariant object $X$, as the object of \tdef{declared $H$-fixed points}.\footnote{Also called \textit{categorical fixed points} in the equivariant homotopy literature.} Functoriality in $\cO_G^\op$ makes these objects interact as actual fixed-point data would: subgroup inclusions induce \textit{restriction} maps, while elements of $G$ induce \textit{conjugation} isomorphisms.

\begin{rem}
    The unit of the adjunction
    \[
        \Fun(\cO_G^\op, \cC) \adj \Fun(BG, \cC)
    \]
    consists of natural maps 
    \[
        X(G/H) \too X^{hH} \qin \cC.
    \]
    Thus, the declared fixed points \textit{refine} the homotopy fixed points of the underlying naive $G$-object. Moreover, a $G$-equivariant $X$ is \textit{naive} precisely when all these comparison maps are \textit{isomorphisms}.
\end{rem}

A $G$-equivariant object may begin its life as a naive one arising from an ordinary $G$-action in $\cC$. But after applying a functor $F\colon \cC \to \cD$ levelwise, the resulting $G$-equivariant object in $\cD$ need no longer be of this form, since $F$ may fail to preserve fixed points. This is often precisely how nontrivial $G$-equivariant objects arise. 

\begin{example}
        Let $X$ be a \textit{topological space} with a continuous $G$-action. The actual \textit{strict} fixed points $X^{H}\sseq X$ for the various subgroups $H \le G$ need not coincide with the \textit{homotopy} fixed points~$X^{hH}$. Nevertheless, by pushing the naive $G$-equivariant topological space $X$ along the canonical functor $\mathrm{Top} \to \Spc$ we get a $G$-equivariant \textit{space}\footnote{In fact, $X^H$ is a common notation for declared $H$-fixed points for arbitrary $G$-equivariant objects $X$.} 
    \[
        X \colon \cO_G^\op \too \Spc \qquad\text{with}\qquad
        X(G/H) := X^{H} \qin \Spc.
    \]
    For example, there is a natural $G$-equivariant space $\mdef{E}G$ characterized by
    \[
        EG(G/H) \simeq \begin{cases}
            \pt & H = e \\
            \es & else.
        \end{cases}
    \]
    Its underlying naive $G$-space is just $\pt$ with the trivial $G$-action, but unlike it, it has ``no fixed points''. 
    This construction brings strict equivariant structures within reach of homotopy-theoretic methods. For example, assuming $X$ is sufficiently nice, we have the following \textit{homotopy} colimit formula for the \textit{strict} orbits $X/G$ in terms of the \textit{strict} fixed points $X^H$,
    \[
        X/G \simeq \colim_{G/H\in \cO_G^\op}\ X^H 
        \qin \Spc,
    \]
    which in turn gives rise to a \textit{homology spectral sequence}, 
    \[
        E^2_{n,s} = H_s(\cO_G^\op, H_{n-s}(X^{(-)})) \implies 
        H_n(X/G).
    \]
\end{example}

Not all $G$-equivariant objects arise in this way. In other cases, the declared fixed points reflect nontrivial structure intrinsic to the situation itself.

\begin{example}
    The $C_2$-action on $\Cat_\st$ by $\op$ admits a refinement to a \textit{$C_2$-equivariant} $\oo$-category with declared $C_2$-fixed points the $\oo$-category $\Cat_\oo^{\rm p}$ of \textit{Poincar\'{e}} categories. 
    These are stable $\oo$-categories $\cC$ with a quadratic functor 
    \[
        Q\colon \cC^\op \too \Sp,
    \]
    whose bilinear cross-effect is nondegenerate.
    This is a more refined input for Hermitian $K$-theory than just a stable $\oo$-category with duality.\footnote{For the full story start here \cite{calmes2020hermitian}. It's a long ride.}  
\end{example}

\begin{rem}
    The microcosm principle dictates that the true home of a $G$-equivariant object is a $G$-equivariant \textit{$\oo$-category} $\cC \colon \cO_G^\op \to \Cat_\oo$. In this setting, $X(G/H)$ lives in $\cC(G/H)$. Namely, fixed points with respect to different \textit{subgroups} live in different \textit{categories}. This, not at all perverse, level of generality is the subject of \textit{parametrized} category theory.\footnote{For a systematic development start with \cite{martini2021colimits}.}
\end{rem}

\subsection{Transfer maps}

All this was just a warm-up. The real action begins when there is also an \textit{additive} structure in the picture. For an ordinary commutative monoid $M$ with an action of a \textit{finite} group $G$, the various fixed point submonoids are related not only by restriction and conjugation, but also by \textit{transfers}. 
For instance, any $x \in M$ determines a $G$-fixed point by summing over its orbit:
\[
    \sum_{g \in G} gx \qin M^G.
\]
This construction depends only on the $G$-orbit of $x$ and hence  defines a map 
\[
    M_G \too M^G. 
\]
Analogues of such maps exist also in the homotopy-coherent setting. 

\begin{defn}
    Let $\cC$ be an $\oo$-category and $G$ a finite group. For a commutative monoid $M\in \CMon(\cC)$ with a $G$-action, there is a canonical \tdef{transfer} map\footnote{Or \textit{norm}, or \textit{trace}, depending on the context and author.}
    \[
        \mdef{\tr}_G \colon M_{hG} \too M^{hG} \qin \cC,
    \]
    which roughly speaking ``sums over the orbit''.
\end{defn}
The construction relies on the \textit{semiadditivity} property of the $\oo$-category $\CMon(\cC)$. Since finite products and coproducts coincide, the \textit{induced} and \textit{coinduced} $G$-representations on $M$ in $\CMon(\cC)$ are canonically identified with the $G$-fold biproduct $\bigoplus_G M$ equipped with the summand-permutation action. We can now compose the coinduction unit map with the (nonequivariant) fold map:
\[\begin{tikzcd}[row sep=tiny]
	M & {\bigoplus_G M} & M \\
	x & {(gx)_{g\in G}} & {\sum_{g\in G} gx.}
	\arrow["\eta", from=1-1, to=1-2]
	\arrow["\nabla", from=1-2, to=1-3]
	\arrow[maps to, from=2-1, to=2-2]
	\arrow[maps to, from=2-2, to=2-3]
\end{tikzcd}\]
This is indeed ``summation along the orbit'' though not yet factored through the (co)invariants. To achieve \textit{that}, one has to be a bit clever and interpret $\bigoplus_G M$ as the (co)induction of $M$ along the \textit{diagonal} $G \to G\times G$. Under this perspective, $\eta$ and $\nabla$ become $G\times G$-equivariant.  The left copy of $G$ acts trivially on the source $M$ and ordinarily on the target, while the right copy acts ordinarily on the source and trivially on the target $M$. This equivariance precisely induces the desired factorization.\footnote{A rigorous construction using Kan extension calculus appears in \cite[\S 6.1.6]{HA}, vastly generalized in \cite[\S 4]{AmbiKn}.}

The transfer map $\tr_G$ has variants where one sums only over the orbits of a subgroup $H\le G$ and starts already with a fixed point for a smaller subgroup $K\le H$. All these compose and conjugate\footnote{The factorization of the transfer through the \textit{coinvariants} is closely related to the compatibility with conjugation. We shall not dwell on this point.} in a natural way and organize into a \textit{covariant} functor
\[
    X^{(-)} \colon \cO_G \too \cC .
\]
The most intricate part of the story is the interaction between these \textit{transfers} and the \textit{restrictions}, governed by the Mackey \textit{double-coset formula}. Stating it will require some notation. For a subgroup $K \le H$ we denote the induced restriction and transfer maps by
\[
    \tr_K^H \colon 
    X^{hK} \adj 
    X^{hH}
    \noloc \res_K^H,
\]
and for $g\in G$, we denote the induced \textit{conjugation} map by
\[
    X^{hH} \isob[\ \ c_g \ \ ] X^{h(gHg^{-1})}
\]

\begin{prop}\label{Double_Coset_Formula}
    For $M \colon BG \to \CMon(\cC)$ and subgroups 
    $K_1, K_2\le H \le G$,
    the composite of restriction and transfer
    \[
        \res_{K_2}^H\circ \tr_{K_1}^H \colon
        M^{hK_1} \too M^{hK_2}
    \]
    is homotopic to the sum
    \[
        \sum_{[g]\in K_2\backslash H/K_1}
        \tr_{K_2\cap gK_1g^{-1}}^{K_2}
        \circ c_g\circ
        \res_{K_1\cap g^{-1}K_2g}^{K_1}.
    \]
\end{prop}

In light of this, any collection of \textit{declared fixed point} commutative monoids $M(G/H)$ aiming to impersonate the fixed-point data of a $G$-action on some $M\in \CMon(\cC)$ would have to carry this elaborate web of maps and relations --- not to mention the \textit{coherence} homotopies for repeated applications of the above formula. Faced with the daunting task of formalizing all this, the reader may be relieved to learn that it can be packaged neatly in one go. It all comes down to identifying the \textit{combinatorial origin} of double-cosets. 

\begin{prop}
    Let $G$ be a finite group with subgroups $K_1,K_2 \le H \le G$. There is a pullback diagram in $\Fin_G$,
    \[\begin{tikzcd}
    	{\bigsqcup_{[g]\in K_2\backslash H/K_1} G/(K_1\cap g^{-1}K_2g)} & {G/K_2} \\
    	{G/K_1} & {G/H.}
    	\arrow[from=1-1, to=1-2]
    	\arrow[from=1-1, to=2-1]
    	\arrow["\lrcorner"{anchor=center, pos=0.125}, draw=none, from=1-1, to=2-2]
    	\arrow[from=1-2, to=2-2]
    	\arrow[from=2-1, to=2-2]
    \end{tikzcd}\]
\end{prop}

With this in hand, the following definition is hopefully not entirely surprising. 

\begin{defn}
    Let $G$ be a finite group. The $\oo$-category of \tdef{genuine} $G$-commutative monoids in $\cC$ is given by
    \[
        \CMon_{\mdef{G}}(\cC) := 
        \Fun^\times(\Span(\Fin_G),\cC).
    \]
\end{defn}
Given such a finite product preserving functor $M$, the restriction along the inclusions
\[
    \Fin_G \intoo \Span(\Fin_G) \longhookleftarrow \Fin_G^\op
\]
produces functors $\Fin_G \to \cC$ and $\Fin_G^\op \to \cC$ taking finite coproducts and products respectively to products in $\cC$. These in turn correspond to functors $\cO_G \to \cC$ and $\cO_G^\op \to \cC$ with a common value $M(G/H)$ on objects. These are the transfer and restriction functorialities (including conjugation). Finally, the composition of spans, defined via pullback, is precisely the homotopy-coherent Mackey double-coset formula. 

Thus, all one has to do in order to add $G$-equivariance \textit{additively} into commutative monoids is to replace the algebraic theory $\Span(\Fin)$ with $\Span(\Fin_G)$. What might not be immediately clear is that for a genuine $G$-commutative monoid $M$ all the declared fixed point objects $M(G/H)$ are themselves commutative monoids. There is a unique finite coproduct preserving functor $\Fin \to \Fin_G$ taking $\pt$ to $G/H$. This functor preserves pullbacks and thus induces a product preserving  functor
\[
    \Span(\Fin) \too \Span(\Fin_G),
\]
restriction along which endows $M(G/H)$ with a commutative monoid structure. 

\begin{example}
    A novice in equivariant homotopy theory is well advised to spend about 90\% of their time contemplating the case $G=C_2$. A genuine $C_2$-equivariant commutative monoid amounts to the following diagram of commutative monoids:
    \[\begin{tikzcd}
    	{M(C_2/C_2)} && {M(C_2/e)}
    	\arrow["\res"', shift right, from=1-1, to=1-3]
    	\arrow["\tr"', shift right, from=1-3, to=1-1]
    	\arrow["{C_2}", from=1-3, to=1-3, loop, in=55, out=125, distance=10mm]
    \end{tikzcd}\]
    such that $\res$ and $\tr$ are $C_2$-invariant, and the composite
    \[
        M(C_2/e)_{hC_2} \oto{\ \tr\ }
        M(C_2/C_2) \oto{\ \res \ }
        M(C_2/e)^{hC_2}
    \]
    is the (honest) transfer map. The endomorphism 
    \[
        T:=\tr\circ\res \qquad\text{of}\qquad M(C_2/C_2)
    \]
    need \textit{not} be multiplication by $2$, as one might naively expect from summing the two $C_2$-translates of a ``fixed point.''  It can be shown though to satisfy the weaker relation $T^2 = 2 T$. Indeed, while $C_2/e$ is different from two copies of $C_2/C_2$, we do have
    \[
        C_2/e \times C_2/e \ \simeq\ 
        C_2/e \sqcup C_2/e \qin \Fin_{C_2}.
    \]
\end{example}
\newcommand{\Res}{\mathrm{Res}}

\begin{example}
    Mackey's double-coset formula is most familiar in the context of \textit{representation theory}. To an inclusion $K\le H$ of finite subgroups corresponds a \textit{restriction--induction} adjunction of representations,
    \[
        \Ind_K^H \colon 
        \Vect_\CC^{BK} \adj 
        \Vect_\CC^{BH} 
        \noloc \Res^H_K.
    \]
    And for $K_1,K_2 \le H$ one can similarly express the composition $\Res^H_{K_2} \circ \Ind_{K_1}^H$ by a formula analogous to that of \Cref{Double_Coset_Formula}. This is indeed a special case. For a finite group $G$, all this structure organizes into a genuine $G$-equivariant commutative monoid in \textit{categories},
    \[
        \Vect_\CC^{(-)/\!\!/G} \colon 
        \Span(\Fin_G) \too \Cat,
    \]
    taking an orbit $G/H$ to the category of $H$-representations $\Vect_\CC^{BH}$. A special feature of this example is that the transfer maps (ie induction functors) are not an additional structure, but come as \textit{adjoints} of the restriction functors, leveraging the 2-categorical structure of $\Cat$.\footnote{The double coset formula follows from further having the adjointability (aka Beck-Chevalley) condition.} This is another instance of a central theme of categorification: structure becomes property.
\end{example}

\subsection{Genuine equivariant spectra}

If one wishes to restrict to genuine $G$-equivariant commutative \textit{groups}, one has to impose invertibility on all declared fixed point monoids. 

\begin{defn}
    Let $G$ be a finite group. The $\oo$-category of \tdef{genuine $G$-commutative groups} in $\cC$ is given by the full subcategory
    \[
        \mdef{\CGrp}_{G}(\cC) \ \sseq\ \CMon_{G}(\cC)
    \]
    on objects $M$ with $M(G/H)$ a group for all $H\le G$.
\end{defn}

In light of the equivalence $\CGrp(\Spc) \simeq \Sp_{\ge0}$, we might call $\CGrp_{G}(\Spc)$ \textit{connective genuine $G$-equivariant spectra}. There is a more direct description of this $\oo$-category for which one must contemplate for a moment the $\oo$-category 
\[
    \CMon_G(\CMon(\Spc)) \simeq 
    \CMon(\CMon_G(\Spc)) \simeq 
    \CMon_G(\Spc).
\]
By semiadditivity of $\CMon_G(\Spc)$ the two commutative monoid structures coincide. In particular, there is no difference between imposing invertibility \textit{inside} and \textit{outside}.

\begin{prop}
    There is an equivalence of $\oo$-categories
    \[
        \CGrp_G(\Spc) \simeq 
        \CMon_G(\Sp_{\ge 0}). 
    \]
\end{prop}

This insight can be further exploited to extend the notion of genuine $G$-equivariant objects to \textit{nonconnective} spectra. This is the central object of study in  \textit{stable equivariant homotopy theory}. 

\begin{defn}
    Let $G$ be a finite group. The $\oo$-category of \tdef{genuine $G$-spectra} is given by
    \[
        \Sp_{\mdef{G}} := \CMon_G(\Sp).
    \]
\end{defn}

\begin{rem}
    This approach to genuine equivariant spectra via \textit{spectral Mackey functors} was not the original perspective. Historically, $\Sp_G$ was constructed by stabilizing  $G$-equivariant pointed spaces with respect to all of the operations
    \[
        X \longmapsto S^V\wedge X =: \Sigma^{\mdef{V}}X 
        \qin \Fun(\sO_G^\op , \Spc_*),
    \]
    where $S^V = V^+$ is the one-point compactification of a finite dimensional $G$-representation $V$, known as a \tdef{representation sphere}. A pro the comparison can be found in \cite[Appendix A]{clausen2020descent} (as well as in the further references cited there).
\end{rem}

Just as one can forget from spectra to spaces, one can forget from genuine equivariant $G$-spectra to equivariant $G$-spaces,
\[
    \Omega^\oo_{\mdef{G}} \colon \Sp_G \too \Spc_G,
\]
by restricting along $\sO_G^\op \into \Span(\Fin_G)$ and applying $\Omega^\oo$ levelwise.\footnote{Here we focus on the unpointed version, but there is an obvious pointed version as well.} The left adjoint of $\Omega^\oo_G$ is denoted by
\[
    \Sigma^\oo_{\mdef{G},+} \colon \Spc_G \too \Sp_G.
\]
With this setup one can associate to every genuine $G$-equivariant spectrum $E$ a \textit{$G$-equivariant cohomology theory}. We make use of the following fact:
\begin{prop}
    The $\oo$-category $\Sp_G$ is stable.
\end{prop}
\begin{proof}
    The $\oo$-category $\Sp$ is stable and hence the functor category $\Fun(\Span(\Fin_G),\Sp)$ is stable as (co)limits in functor $\oo$-categories are computed pointwise. The full subcategory 
    \[
        \Gray{\Sp_G :=}
        \Fun^\times(\Span(\Fin_G),\Sp) \sseq
        \Fun(\Span(\Fin_G),\Sp)
    \]
    is closed under limits and colimits. Indeed, since finite products in $\Sp$ are also finite coproducts, arbitrary pointwise limits and colimits of finite product preserving functors remain finite product preserving. Thus $\Sp_G$ is closed under limits and colimits in the stable $\oo$-category $\Fun(\Span(\Fin_G),\Sp)$, and is therefore itself stable.
\end{proof}

This means that the mapping spaces in $\Sp_G$ canonically refine to \textit{mapping spectra}. For every $X \in \Spc_G$, we associate the sequence of groups 
\[
    E^n_{\mdef{G}}(X) := \pi_{-n}\hom(\Sigma^\oo_{G,+}X, E)
    \quad\text{for}\quad n \in \ZZ.
\]

\begin{example}
    The most famous example of a $G$-equivariant cohomology theory is \textit{$G$-equivariant topological $K$-theory}. Its degree-zero part assigns to a finite $G$-equivariant space $X$ the group 
    \[
        K^0_{\mdef{G}}(X) = \pi_0(\Vect_{\CC,\mdef{G}}(X))^\gp \qin \Ab,
    \]
    where $\Vect_{\CC,G}(X)$ is the space of finite dimensional $G$-equivariant $\CC$-vector bundles on $X$. In particular, $K_G^0(\pt)$ is the usual \textit{representation ring} of $G$. This equivariant cohomology theory is represented by a genuine $G$-equivariant spectrum $\KU_{\mdef{G}}$ constructed by a similar process to $\KU$, starting with the topological commutative monoid of finite dimensional $G$-representations. 
\end{example}

We conclude this discussion by introducing one fundamental example of a genuine $G$-equivariant spectrum --- the \textit{$G$-equivariant sphere spectrum}. 
The ordinary sphere spectrum $\Sph$ is the free homotopy-coherent commutative group on one element. In spectral notation, this reads as $\Sph \simeq \Sigma_+^\oo(\pt)$. 
In the equivariant world, the \tdef{$G$-equivariant point} is the naive $G$-equivariant space $\pt_{\mdef{G}}$ associated to $\pt$ with the (necessarily) trivial $G$-action. 

\begin{defn}
    The \tdef{$G$-equivariant sphere spectrum} is given by
    \[
        \Sph_{\mdef{G}} := \Sigma_{G,+}^\oo(\pt_G) \qin \Sp_G.
    \]
\end{defn}
It is not immediately clear what the declared \(H\)-fixed-point spectra of $\Sph_G$ are, since we have not written down a formula for $\Sigma_{+,G}^\oo$. As in the nonequivariant case, we can instead pass through commutative monoids and then group-complete. The relevant left adjoints give the following factorization of $\Sigma_{+,G}^\oo$:
\[\begin{tikzcd}[sep = 4em]
	{\Spc_G} & {\CMon_G(\Spc)} & {\CGrp_G(\Spc)} & {\Sp_G,}
	\arrow["{F_{\CMon,G}}", from=1-1, to=1-2]
	\arrow["{(-)^{\gp}}", from=1-2, to=1-3]
	\arrow[hook, from=1-3, to=1-4]
\end{tikzcd}\]
where $G$-equivariant group completion is obtained by applying ordinary group completion levelwise to the declared $H$-fixed points.
We let the \tdef{$G$-equivariant matural numbers} be 
\[
    \MM_{\mdef{G}} := F_{\CMon,G}(\pt_G)  \qin \CMon_G(\Spc).
\]
Composing the adjunctions, we find that  $\Sph_G$ is simply the group completion of $\MM_G$. The chief virtue of the latter is its \textit{naivety}. 

\begin{prop}
    For every subgroup $H \le G$ we have a canonical isomorphism
    \[
        \MM_G(G/H) \ \simeq\  \MM^{BH} \qin \CMon(\Spc).
    \]
    Namely, the underlying $G$-equivariant space of $\MM_G$ is just $\MM$ with the naive trivial $G$-action.\footnote{It nevertheless admits the additional \textit{transfer} maps as part of the genuine structure.}
\end{prop}
\begin{proof}[Proof sketch]
    Consider the following commutative square, where the bottom horizontal map is~$F_{\CMon,G}$:
    \[\begin{tikzcd}
    	& {\Fin_G} & {\Span(\Fin_G)^\op} & \\
    	\textcolor{rgb,255:red,179;green,179;blue,179}{{\Spc_G}} & {\Fun^\times(\Fin_G^\op,\Spc)} & {\Fun^\times(\Span(\Fin_G),\Spc)} & \textcolor{rgb,255:red,179;green,179;blue,179}{{\CMon_G(\Spc)}}
    	\arrow[from=1-2, to=1-3]
    	\arrow[hook, from=1-2, to=2-2]
    	\arrow[hook, from=1-3, to=2-3]
    	\arrow[color={rgb,255:red,179;green,179;blue,179}, equals, from=2-2, to=2-1]
    	\arrow[from=2-2, to=2-3]
    	\arrow[color={rgb,255:red,179;green,179;blue,179}, equals, from=2-3, to=2-4]
    \end{tikzcd}\]
    Since $\pt_G$ is the Yoneda image of the trivial orbit $G/G$, we get 
    \[
        \MM_G(G/H) \simeq 
        \Map_{\Span(\Fin_G)}(G/G,G/H) \simeq 
        ((\Fin_G)_{/(G/H)})_\isom \simeq 
        (\Fin_H)_\isom \simeq
        \MM^{BH}.
    \]
    The penultimate equivalence takes a $G$-equivariant map $X \to G/H$ to its fiber over $eH \in G/H$, which is an $H$-set. 
\end{proof}

Naivety does not survive group completion, though we can still analyze the answer. On the level of $\pi_0$ the result is known as the \textit{Burnside ring},
\[
    \pi_0\Sph_G(G/H) \simeq \pi_0(\MM^{\mathit{BH}})^\gp =: \sA(H).
\]
As an abelian group, $\sA(H)$ is freely generated by the \textit{isomorphism} classes of $H$-orbits $H/K$, which are in bijection with the \textit{conjugacy} classes of subgroups $K \le H$.\footnote{Multiplication is induced from the \textit{product} of $H$-sets.}

To express the spectrum $\Sph_G(G/H)$ itself, we also need to recall that the automorphisms of $H/K$ in $\Fin_H$ are given by the \textit{Weyl group},
\[
    \mdef{W\!}_H(K) := N_H(K)/K.
\]
With this we can give a formula for the declared fixed points of $\Sph_G$, known as the \textit{tom Dieck splitting}.

\begin{thm}[tom Dieck]
    For every $H \le G$ we have an isomorphism
    \[
        \Sph_G(G/H) \ \simeq\ \bigoplus_{(K)}\ \Sph\angles{BW_H(K)} \qin \Sp,
    \]
    where the sum is over conjugacy classes of subgroups $K \le H$.
\end{thm}
\begin{proof}[Proof sketch]
    The connected component of an orbit $H/K$ in $\MM^{\mathit{BH}}$ is equivalent to the classifying space of its automorphisms $BW_{H}(K)$. As every finite $H$-set decomposes uniquely as a finite coproduct of orbits of this form, $\MM^{\mathit{BH}}$ is just the \textit{free} homotopy-coherent commutative monoid on the disjoint union of all the corresponding classifying spaces\footnote{A more rigorous argument might use \Cref{F_CMon_Cart}.}
    \[
        X = \bigsqcup_{(K)}BW_H(K) \qin \Spc.
    \]
    Assuming this, we immediately deduce the desired isomorphism
    \[
        \Sph_G(G/H) \simeq 
        (\MM^{\mathit{BH}})^\gp \simeq 
        F_{\CMon}(X)^\gp \simeq 
        F_{\CGrp}(X) \simeq 
        \bigoplus_{(K)}\ \Sph\angles{BW_H(K)} \qin \Sp.
    \]
\end{proof}
The commutative monoid $\MM^{BG}$ is also the constant $BG$-shaped \textit{limit} on $\MM$. Group completion does not preserve this limit; nevertheless, there is a canonical assembly map 
\[
    \Sph_G(G/G) \simeq 
    (\MM^{BG})^\gp \too 
    (\MM^\gp)^{BG} \simeq
    \Sph^{BG}
    \qin \Sp.
\]
We have just computed the left-hand side in terms of \textit{colimits} of copies of the sphere spectrum. The right-hand side, however, is quite mysterious: it is the \textit{cohomotopy} of $BG$. One of the wildest facts in stable homotopy theory is that this map is almost an equivalence. We define the augmentation ideal $\mdef{I_G}$ to be the kernel of the augmentation map
\[
    \pi_0 \Sph_G(G/G) = \sA(G) \oto{\ \mdef{\varepsilon}\ } \ZZ.
\]

\begin{thm}[Segal conjecture]\!\!\!\footnote{This is the standard name despite the fact it has long been proven.}
    The canonical map above induces an equivalence upon $I_G$-adic completion of the source\footnote{When $G$ is a $p$-group, it is absorbed by $p$-completion.}
    \[
        \Sph_G(G/G)_{I_G}^\wedge \iso \Sph^{BG} 
        \qin \Sp.
    \]
    In particular,
    \[
        \pi_0 (\Sph^{BG}) \simeq \sA(G)_{I_G}^\wedge.
    \]
\end{thm}
The immediate corollary that $\Sph^{BG}$ is \textit{connective} is by itself bewildering. The naive approach to computing $\pi_*(\Sph^{BG})$ would be to apply $(-)^{BG}$ to the Postnikov filtration of $\Sph$, which gives rise to a spectral sequence
\[
    E_2^{st} = H^s(BG, \pi_t^{st}) \implies 
    \pi_{t-s}(\Sph^{BG}). 
\]
The $E_2$-page contains infinitely much complicated junk below the diagonal; yet by $E_\oo$, it somehow manages to cancel itself out through an intricate conspiracy of the differentials. 

\subsection{Multiplicative norms}

We have thus far considered only the \textit{additive} aspects of equivariant structures, but one might as well have \textit{multiplication} thrown into the mix. Consider an ordinary commutative ring $R$ with a $C_2$-action. Namely, a commutative ring automorphism 
\[
    \sigma \colon R \iso R \qquad\text{with}\qquad \sigma^2 = \Id.
\]
From an element $x\in R$ we can produce a $C_2$-fixed point either by \textit{summing} or by \textit{multiplying} over its $C_2$-translates,  
\[
    \tr(x) = x+\sigma (x) \qquad\text{and}\qquad 
    \mdef{\Nm}(x) = x \cdot \sigma(x).
\]
The first construction is the transfer discussed above for the underlying \textit{additive} commutative monoid, while the second is the~\tdef{norm}, which is the transfer for the underlying \textit{multiplicative} commutative monoid. The essential point, however, is not merely that the transfer and norm coexist. The norm is multiplicative rather than additive, and its failure to preserve addition is itself controlled by transfers:
\[
    \Nm(x+y) = \Nm(x) + \tr(x\sigma(y))  + \Nm(y).
\]
Since additive equivariant structure is governed by spans of finite $G$-sets, while compatible addition and multiplication are governed (nonequivariantly) by bispans of finite sets, one is naturally led to consider bispans of finite $G$-sets.

\begin{defn}
    For a finite group $G$, the $\oo$-category of  \tdef{genuine $G$-equivariant commutative semirings} in an $\oo$-category $\cC$ is given by\footnote{The interested reader is encouraged to locate the above formula for $\Nm(x+y)$ in the bispan calculus of finite $C_2$-sets.}
    \[
        \CRig_{\mdef{G}}(\cC) := \Fun^\times(\Bispan(\Fin_G), \cC).
    \]
    The full subcategory of genuine $G$-equivariant commutative \textit{rings} is defined in the obvious way and denoted by $\CRing_G(\cC)$. 
\end{defn}

In the case of $\cC = \Spc$, this notion coincides with the more classical one of (connective) \textit{normed $G$-ring spectra} defined in terms of tensor products in spectra.\footnote{This is \cite[Theorem B]{cnossen2024normed}.} 
Importantly, though, not only are these not $G$-equivariant objects in commutative ring spectra, but they are also not commutative algebras in genuine $G$-equivariant spectra. To capture this notion one has to venture deeper into the equivariant world --- to \textit{$G$-symmetric monoidal categories}. Namely, one has to promote $\Sp_G$ itself to an object 
\[
    \underline{\Sp}_G \ \in \ \CMon_G(\Cat_\oo) \qquad\text{with}\qquad \underline{\Sp}_G(G/H) = \Sp_H.
\]
This requires, in particular, constructing for every spectrum $X$, a genuine $G$-equivariant spectrum $N^G_eX$ with underlying naive $G$-spectrum $\bigotimes_G X$. Such operations were constructed in \cite{HillHopkinsRavenel2016} and played a crucial role in the solution of the \textit{Kervaire invariant one} problem,\footnote{The last remaining case, dimension $126$, was resolved positively in \cite{lin2024last}.} a landmark achievement of stable equivariant homotopy theory.

\section{The view from here}

Before turning to the next chapter, it is worth pausing to look back at the paths our characters have taken so far and, with the benefit of hindsight, trace the arcs of their development and place them in a broader context. In ordinary mathematics, algebraic structures are commonly presented as a list of operations of varying arity, subject to a list of identities among their compositions. The inciting incident of our story was the realization that this definition is inadequate in the homotopy-coherent world: as identities become (chosen) homotopies, one is left to ponder the higher relations among them. The systematic solution is not to enumerate these higher relations one by one, but rather to draw out the shape of the diagram that governs the structure as a whole.  

There is genuinely no single right answer here, but there is a canonical one called \textit{animation}. The full subcategory spanned by the free objects on finitely many generators encodes all the operations of all arities of the structure, together with their full ``composition table''. This is its \textit{algebraic theory}, which can be realized in any $\oo$-category by taking product preserving functors from its opposite. The equations built into the composition table then unfold into homotopies, higher homotopies between them etc.\! producing an infinite hierarchy of coherences. In this process the higher structure is not imposed by hand, but is already latent in the classical theory and emerges once the target is given room for homotopy. This recipe produces homotopy-coherent unital and nonunital monoids, (connective) derived categories of rings, and various other useful things.

Animation, however, is by no means the only option. 
And while for an arbitrary algebraic structure there is no objective ranking on its homotopy-coherent refinements, for certain structures we have rather concrete expectations. Homotopy-coherent \textit{associativity} is one such case. Ordinary associativity does not by itself determine its preferred homotopy-coherent refinement, but monoids occupy a sufficiently special place in category theory to force the answer: just as monoids are one-object categories, homotopy-coherent monoids ought to be precisely one-object \(\infty\)-categories. Animated monoids pass this test, though not tautologically. 

The case of homotopy-coherent \textit{commutativity} is another, quite different, example. We were led beyond $\Mat(\NN)$, the ordinary algebraic theory of commutative monoids, because its models in spaces failed to accommodate natural examples such as $\MM$. Passing to $\Mat(\MM)$ solved this problem, but did not yet explain why it captures the right notion of a homotopy-coherent commutative monoid in general. The missing justification comes from the categorical source of commutative monoid structure --- \textit{semiadditivity}. This property is defined for $\oo$-categories in essentially the same way as for ordinary categories and $\Mat(\MM)$ is precisely the free semiadditive $\oo$-category on one object, while $\Mat(\NN)$ is only its homotopy category. Thus, the same categorical property that singles out the classical theory also singles out its homotopy-coherent refinement, forcing higher coherence into the algebraic theory itself rather than only into its values. In more concrete terms, the operations themselves acquire symmetries. Indeed, the $n$-fold sum of an element of a commutative monoid with itself admits a natural, and generally nontrivial, action of $\Sigma_n$.

The passage from $\Mat(\MM)$ to the equivalent $\oo$-category $\Span(\Fin)$ is not merely a clever application of general categorical machinery (as the reader may have been led to believe), but in fact represents a marked shift in perspective. Morphisms in $\Span(\Fin)$ factor uniquely as a composition of two special types of maps\footnote{What is called a \textit{factorization system}.} --- the backward legs, and the forward ones --- doing two different jobs. The backward maps $Y \ofrom{f} X \oto{=} X$, constituting a copy of $\Fin^\op$, encode \textit{restrictions} $f^*\colon M^Y \to M^X$, while the forward maps $X \ofrom{=} X\oto{g} Y$ encode \textit{summation} along the fibers $g_!\colon M^X \to M^Y$. Put this way, it is natural to contemplate more general ``spaces'' than finite sets where one can restrict and sum along more general classes of maps: $\pi$-finite, proper submersions, finite etal\'{e}, etc. This is in a sense what spans \textit{are}: the abstract calculus of restrictions and transfers interacting via basechange. Bispans extend this picture by adding an additional \textit{multiplicative} transfer, which distributes over the additive one. This too generalizes to other homotopical, geometric and algebro-geometric situations.  

The equivariant theory fits this broader interpretation rather well. A genuine $G$-equivariant commutative monoid $M$ assigns to finite $G$-sets $X$ objects which we may demagogically denote by $M^X$ related by restrictions and transfers. However, this is not an algebraic structure on a single underlying object; each $M^{G/H}$ is an atomic unit, which may or may not be determined by $M^{G/e}$ with its naive $G$-action. This kind of flexibility can be quite useful. The general restriction--transfer calculus applies in situations where $M^X$ is not literally the $X$-indexed power of a fixed object, and need not even have a fixed type as $X$ varies.

Besides being a major theme in homotopy theory in its own right, equivariant homotopy theory is useful for seemingly nonequivariant questions as well. The two examples discussed above --- the Segal conjecture and the Kervaire invariant problem --- illustrate this. Another natural beneficiary of equivariant techniques is algebraic $K$-theory. Suppose $L/F$ is a finite Galois extension of fields\footnote{Or commutative rings, or commutative ring \textit{spectra}.} with Galois group $G$. While on the level of commutative monoids we have $M_F^{} \simeq M_L^{hG}$, after group completion, the assembly map 
\[
    K(F) = K(L^G)\too K(L)^{hG} \qin \Sp
\]
is rarely an equivalence. This failure of \textit{Galois descent} for algebraic $K$-theory is one of the main hurdles for computing it. Nevertheless, it can also be systematically studied and put to work. The spectra $K(L^H)$ assemble into a genuine $G$-spectrum with restriction and transfer maps induced
by extension and restriction of scalars respectively.
One can get surprising mileage out of this highly structured object.\footnote{A good place to read about some of it is \cite{mathew2016descent}.}

Once one further imposes \textit{invertibility} on the  monoid structure, a close connection to ``geometry'' emerges. Homotopy-coherent groups are loop spaces, which as an $\oo$-category are simply pointed connected spaces. Adding increasing levels of commutativity amounts to adding further deloopings, which in the limit identifies homotopy-coherent commutative groups with infinite loop spaces. Removing the connectivity restrictions on the deloopings produces the $\oo$-category of \textit{spectra}, which is the natural home for homotopy-coherent commutative algebra. 

The relationship between group structures and loop spaces can be somewhat demystified by recalling that spaces are themselves homotopy-coherent groupoids and one-object groupoids are groups just like one-object categories are monoids. There is in fact a ``categorical delooping'' procedure for monoids as well, but it requires one to leave the world of ($\oo$,1)-categories. The endomorphisms of an object in an $(\oo,2)$-category form a monoidal $(\oo,1)$-category. One can further consider the endomorphisms of the unit $\one$ in the latter. A pair of such endomorphisms $f,g\colon \one \to \one$ can be multiplied in \textit{two} ways: by composition $g\circ f$, and by tensor product $g\otimes f$. Moreover, these two multiplications \textit{interchange} and so provide an $\EE_2$-monoid structure on~$\End(\one)$. 
This pattern continues for all \(n\) by iterating the endomorphism object construction. In particular, \(\EE_n\)-monoids in spaces can be identified with pointed \(n\)-connective \((\infty,n)\)-categories.\footnote{See \cite[Corollary 6.3.11]{gepner2015enriched}.} Commutative monoids therefore admit compatible categorical deloopings of arbitrary order. By removing the connectivity constraint and allowing arbitrary $(\oo,\oo)$-categories one arrives at the notion of a \textit{categorical spectrum}.\footnote{See \cite{masuda2024algebra}.}

There is another interesting and more hidden arc of the \textit{bar construction}. In classical homotopy theory it is one of the fundamental constructions relating category theory, topology, and homological algebra. 
In our story,  it made its first appearance merely as an economical model for a homotopy-coherent monoid structure. However, when the monoid in question was assumed to be a group, and hence a loop space $\Omega X$, the bar construction was reinterpreted as the \textit{\v{C}ech resolution} of the basepoint inclusion~$\pt \to X$. Furthermore, as the way back to $X$ is given by realizing the bar construction $X \simeq |\BB \Omega X|$ (ie taking its colimit), we obtained a model for the \textit{group completion} map $M \to \Omega |\BB M|$. 
For a homotopy-coherent \textit{commutative} monoid $M$, the bar construction~$BM = |\BB M|$ was once more reinterpreted, this time as the internal \textit{suspension} in $\CMon(\Spc)$ and accordingly, the commutative group completion became the unit of the loop--suspension adjunction. For a homotopy-coherent commutative \textit{group}, this unit map is an isomorphism and so the bar construction, now suspension, becomes \textit{fully faithful}, which in spectra is further an \textit{autoequivalence}. These different guises of the bar construction are given by the categorical world around it, which assigns the same formal construction a different role each time. 

Finally, throughout the chapter we have witnessed the interplay between \textit{structure} and \textit{properties}. The properties, of increasing strength, an $\oo$-category might enjoy
\[
    \text{plain} \too
    \text{pointed} \too
    \text{semiadditive} \too
    \text{additive} \too
    \text{stable}
\]
are reflected in the canonical extra structure its \textit{mapping spaces} acquire
\[
    \Spc \xrightarrow{\qquad}
    \Spc_* \xrightarrow{\qquad} 
    \CMon(\Spc) \xrightarrow{\qquad}
    \CGrp(\Spc) \xrightarrow{\qquad}
    \Sp.
\]
Conversely, each of these properties is itself equivalent to the \textit{objects} of the category admitting the corresponding structure in a unique way.
Such universal sources of homotopy-coherent structure are useful since exhibiting it in terms of explicit formulae is usually rather difficult. The uniqueness part is similarly handy to compare a priori different approaches for constructing the same structure.\footnote{See \cite{gepner2016universality} for a fuller discussion of this idea.}

Our study of algebraic theories highlighted their versatility in encoding various kinds of algebraic structures and the ways in which they interact with one another. However, it also pointed to some of their limitations and deficiencies. Take for example the notion of a commutative ring. The algebraic theory approach allows one to define this notion as a structure on a set whose $n$-ary operations are encoded by polynomials. But one might also wish to think of a commutative ring as an abelian group, with respect to addition, endowed with additional multiplicative structure. There is, however, one snag: for a ring $R$, the distributivity condition dictates that the multiplication map $R \times R \to R$ is \textit{bilinear} and so not even a morphism in the category of abelian groups. To represent it by a \textit{linear} map, one has to replace the \textit{cartesian product} by the \textit{tensor product} $R\otimes R \to R$. With this modification the definition of a (commutative) ring is essentially the same as that of a (commutative) monoid.

This rearrangement relegates some of the structure to the ambient category. Namely, the distributivity of the ring multiplication is absorbed into the definition of the tensor product of abelian groups. 
From a broader perspective, an algebraic theory encodes certain $n$-ary operations, but one can be more subversive in interpreting this notion. 
One can replace the cartesian product $\times$ by a more general ``product'' $\otimes$, so that an $n$-ary operation on $X$ takes the form $X^{\otimes n}\to X$ rather than $X^n\to X$.
For this to make sense, $\otimes$ has to share enough of the properties of $\times$, like homotopy-coherent associativity, commutativity, and unitality. In other words, it needs to be a symmetric monoidal structure, which is none other than a commutative monoid in categories.\footnote{It does not, however, need to share \textit{all} of the properties of a cartesian product, such as the existence of a diagonal $X \to X\times X$. Thus, while monoids generalize to algebras, \textit{groups} do not have a similar interpretation in an arbitrary symmetric monoidal category.} 
Thus, doing algebra \textit{in} categories has brought us to doing algebra \textit{on} categories.
And there is more to this than the study of (commutative) algebra objects. One can also study invertible, dualizable, and idempotent objects and the resulting theories of traces, localizations etc.  

Importantly, the tensor product of a general (symmetric) monoidal category is not determined by the underlying category, but is an externally imposed structure. In particular, it can do more (or, rather, other things) than represent bilinear maps. With an appropriate choice one obtains notions such as graded rings, differential graded rings, and ultimately (nonconnective) ring spectra. More radically still, even the meaning of \textit{commutativity} depends on the chosen symmetry. The Koszul sign rule for graded abelian groups, for instance, produces graded-commutative rings, in which elements of odd degree \textit{anticommute}. 
At a further remove, one can dispense with commutativity altogether, and the monoidal product itself may have an entirely different provenance.
The endomorphism category of a category, for example, is monoidal under composition. Its algebras are precisely \textit{monads}, which, as we shall explain (much) later, provide yet another route to encoding algebraic structures.

\bibliographystyle{alpha}
\phantomsection\addcontentsline{toc}{section}{\refname}
\bibliography{four}



\end{document}